\documentclass[11pt]{article}
\usepackage{amsmath,amssymb,amscd,enumerate,eucal}
\numberwithin{equation}{section} 
\newcommand{\PP}{\mathcal{P}}
\newcommand{\Gll}{\mathsf{Gl}}
\newcommand{\R}{\mathbb{R}}

\newcommand{\Z}{\mathbb{Z}}

\newcommand{\Ad}{\text{Ad}}
\newcommand{\Aut}{\text{Aut}}

\newcommand{\HH}{\mathbb{H}}
\newcommand{\YMe}{\mathcal{YM}_{\epsilon}}
\newcommand{\YM}{\mathcal{YM}}
\newcommand{\q}{\mathsf{q}}
\newcommand{\A}{\mathsf{a}}
\newcommand{\loc}{\text{loc}}
\newcommand{\ttilde}{\tilde{\intercal}}
\newcommand{\tperp}{\tilde{\perp}}
\newcommand{\Je}{\mathcal{J}_{\epsilon}}
\newcommand{\tr}{\mathsf{Tr}}

\begin{document}
\setlength{\baselineskip}{17pt}

\title{Small coupling limit and multiple solutions to the Dirichlet
Problem for Yang Mills connections in $4$ dimensions - Part I}
\author{Takeshi Isobe \thanks {%
Tokyo Institute of Technology; email:
isobe.t.ab@m.titech.ac.jp} $\,$ and Antonella Marini \thanks {%
University of L'Aquila / Yeshiva University; email:
marini@dm.univaq.it / marini@yu.edu}}


\date{}
\maketitle

\begin{abstract}
{\footnotesize \baselineskip 4mm} In this paper (Part I) and its
sequels (Part II and Part III), we analyze the structure of
the space of solutions to the $\epsilon$-Dirichlet problem for the Yang-Mills
equations on the $4$-dimensional disk, for small values of the
coupling constant $\epsilon$. These are in one-to-one correspondence with solutions
to the Dirichlet problem for the Yang Mills
equations, for small boundary data $\epsilon A_0$. We prove the existence of multiple
solutions, and, in particular, non minimal ones, and establish a
Morse Theory for this non-compact variational problem. In part I, we
describe the problem, state the main theorems and do the first part
of the proof. This consists in transforming the problem into a
finite dimensional problem, by seeking solutions that are
approximated by the connected sum of a minimal solution with an
instanton, plus a correction term due to the boundary. An auxiliary
equation is introduced that allows us to solve the problem
orthogonally to the tangent space to the space of approximate
solutions. In Part II, the finite dimensional problem is solved via
the Ljusternik-Schirelman theory, and the existence proofs are
completed. In Part III, we prove that the space of gauge equivalence
classes of Sobolev connections with prescribed boundary value is a
smooth manifold, as well as some technical lemmas essential to the
proofs of Part I. The methods employed still work when $B^4$ is
replaced by a general compact manifold with boundary, and $SU(2)$ is
replaced by any compact Lie group.
\end{abstract}

\section{Introduction and statement of the main results}
A solution to the Yang Mills equations is a critical point for the
Yang Mills functional defined on the space of connections. These
equations are particularly interesting in $4$ dimensions, since in
this case, the Yang Mills equations are not only invariant under the
infinite dimensional automorphism group of the bundle, namely, the
gauge group, but are also invariant under the group of conformal
transformations over the base manifold. Since the latter is
non-compact, the associated variational  problem is non-compact
(i.e., it never satisfies the Palais-Smale condition) even when
quotiented out by the automorphism group of the bundle. Finding
critical points and establishing a Morse theory for such non-compact
variational problems is one of the most challenging problems in
nonlinear functional analysis.
See~\cite{Aubin},~\cite{Brezis},~\cite{Struwe} and references
therein for an interesting list of non-compact variational problems
and their applications. For the existence of solutions to the Yang
Mills equations on closed manifolds, not necessarily
action-minimizing,
see~\cite{Bor},~\cite{Parker},~\cite{Sadun},~\cite{SS},
~\cite{SSU},~\cite{Taubes1},~\cite{Taubes2},~\cite{Taubes3},~\cite{Wang}.

We establish a relation between the small coupling limit problem and
the problem of existence of multiple critical points for the Yang
Mills functional with small Dirichlet boundary conditions on a
compact 4-dimensional manifold with boundary. We thoroughly analyze
the simplest case, namely, the $SU(2)$-Yang Mills problem on the
$4$-dimensional disk. Our approach is based on a perturbation method
as developed by the first author in~\cite{I2},~\cite{I3},~\cite{I4}
for the $H$-surface equations. See also~\cite{AM} and references
therein for applications of perturbation methods to popular
non-compact variational problems, and, in particular, to Yamabe-like
equations.

In order to describe the problem, we need to establish first some
basic notation. We denote by $B^4$ the open unit disk in $\R^4.$ For
$\epsilon >0$, we define the Lie algebra
$(\mathfrak{su}(2)_{\epsilon},[\cdot,\cdot]_{\epsilon})
:=(\mathfrak{su}(2),\epsilon[\cdot,\cdot])$, where $[\cdot,\cdot]$
is the ordinary Lie bracket on $\mathfrak{su}(2)$. Precisely,
$\mathfrak{su}(2)_{\epsilon}=\mathfrak{su}(2)$ as vector spaces, but
with Lie bracket $[X,Y]_{\epsilon}:=\epsilon [X,Y]$ for
$X,Y\in\mathfrak{su}(2)$. The map
$\phi_{\epsilon}:\mathfrak{su}(2)_{\epsilon}\to\mathfrak{su}(2),$
defined by $\phi_{\epsilon}(X):=\epsilon X$ is a Lie algebra
isomorphism (namely $\phi_{\epsilon} [X,Y]_\epsilon = \epsilon^2
[X,Y] = [\phi_\epsilon X,\phi_\epsilon Y]),$ thus induces a Lie
group isomorphism $\Phi_{\epsilon}:SU(2)_{\epsilon}\to SU(2)$. In
the context of $SU(2)_{\epsilon}$-principal bundles,  the covariant
differentiation associated to a connection $A$ is locally
${d_A}^\epsilon := d + [A, \cdot]_\epsilon := d + \epsilon [A,
\cdot],$ where $\epsilon$ is understood as the coupling constant
that appears in much of the physics literature (cf., for example,
\cite{Ry}, or Ch. 15 of \cite{Wei}).

\medskip
Let now $A_0$ be a smooth connection on an
$SU(2)_{\epsilon}$-principal bundle over $\partial B^4$. Since such
bundles are trivial, $A_0$ lives on the trivial bundle, that is
$A_0\in C^{\infty}(T^{\ast}\partial
B^4\otimes\mathfrak{su}(2)_\epsilon).$ For $A\in L^2_1(T^*B^4\otimes
\mathfrak{su}_\epsilon(2)),$ the $SU(2)_{\epsilon}$-Yang Mills
functional is given by
\begin{equation}
\label{1} \YMe(A)=\int_{B^4}|{F_A}^\epsilon|^2\,dx\,,
\end{equation}
where
${F_A}^\epsilon=dA+\frac{1}{2}[A,A]_{\epsilon}:=dA+\frac{\epsilon}{2}[A,A]$
is the curvature of the connection $A$, and the
$SU(2)_{\epsilon}$-Yang Mills Dirichlet problem in exam, obtained
via a variational method, is given by
$$
\bigl(\mathcal{D}_{\epsilon}\bigr)\quad\qquad\left\{\begin{array}{ll}
{d_A^{\ast}}^\epsilon {F_A}^\epsilon=0\quad&\mbox{in } B^4\\
\iota^{\ast}A\sim A_0\quad&\mbox{at } \partial B^4\,.
\end{array}\right.
$$
Here, $\iota:\partial B^4\to\overline{B}^4$ is the inclusion,
$\iota^{\ast}A\sim A_0$ on $\partial B^4$ means that $\iota^{\ast}A$
is gauge equivalent to $A_0$ over $\partial B^4$ via a gauge
transformation that extends smoothly to the interior,
${d_A^{\ast}}^\epsilon:=\ast
d\ast+\ast[A,\ast\cdot]_{\epsilon}:=\ast
d\ast+\epsilon\ast[A,\ast\cdot]$, where $\ast$ is the Hodge star
operator with respect to the flat metric on $\R^4$. Notice that,
since $\mathfrak{su}(2)_{\epsilon}=\mathfrak{su}(2)$ as sets, $A_0$
may be regarded as an $\mathfrak{su}(2)$-valued $1$-form on
$\partial B^4$, i.e. $A_0\in C^{\infty}(T^{\ast}\partial
B^4\otimes\mathfrak{su}(2))$. This is a canonical identification
between $\mathfrak{su}(2)_{\epsilon}$-valued 1-forms and
$\mathfrak{su}(2)$-valued $1$-forms. A different identification is
given by $\phi_{\epsilon}$. We will use the former identification
throughout this paper, unless we explicitly write otherwise.


\medskip
By the direct method of the calculus of variations,
Marini~\cite{Marini} obtained the first solution to
($\mathcal{D}_{\epsilon}$), that is, an absolute Yang Mills action
minimizing solution, which we shall call \emph{small solution} and
will denote by $\underline{A}_{\epsilon}$. (Note that the
\emph{small solution} is, in general,  not unique, so we choose one
of these solutions for each $\epsilon>0$). Moreover, it is known
(cf.~\cite{IM}) that the space of connections with boundary value
$A_0$, denoted by $\mathcal{A}(A_0)$, has infinitely many connected
components indexed by $\Z$:
$\mathcal{A}(A_0)=\bigsqcup_{k=-\infty}^{\infty}\mathcal{A}_k(A_0)$,
where
$\mathcal{A}_k(A_0)=\{A\in\mathcal{A}(A_0):\mathsf{c}_2(A)=k\}$,
with
$\mathsf{c}_2(A)=\frac{\epsilon^2}{8\pi^2}\int_{B^4}\tr({F_A}^{\epsilon}\wedge
{F_A}^{\epsilon})-\frac{\epsilon^2}{8\pi^2}\int_{B^4}\tr({F_{\underline{A}_{\epsilon}}}^{\epsilon}\wedge
{F_{\underline{A}_{\epsilon}}}^{\epsilon})$ (the relative 2nd Chern
number with respect to $\underline{A}_{\epsilon}$, where
$\underline{A}_\epsilon$ is a fixed absolute minimizer).
In~\cite{IM}, the problem of finding a minimum in each component
$\mathcal{A}_k(A_0)$ is thoroughly solved, yielding many so-called
\emph{large solutions}. In particular, it is proved that, for
non-flat boundary values $A_0$, there exists a minimum at least in
one of the components $\mathcal{A}_{\pm1}(A_0)$. Since all the
solutions known to the Dirichlet problem for Yang Mills are minima
(minimizers for the action restricted to the connected components
$\mathcal{A}_k(A_0)$), it is left open for investigation the
interesting problem whether there exist non-minimal solutions and,
in general, whether the solution found in~\cite{IM} is unique in
each component. (Notice that the results in \cite{Marini, IM} cited
above, proven for the $\mathfrak{su}(2)$-Yang Mills functional, that
is for the standard Dirichlet problem $(\mathcal{D}_1)$ ($\epsilon =
1$), automatically extend to $\YMe(A)).$
This problem can also be related to 
the quantization of Yang Mills theory.

Since the uniqueness result for flat boundary values has been
established in ~\cite{I}, we henceforth assume that $A_0$ is
non-flat, and investigate the existence of non-minimal solutions
and, more in general, seek non-uniqueness results in
$\mathcal{A}_{+1}(A_0)$, (or, by the same arguments, in
$\mathcal{A}_{-1}(A_0))$, since we know that an absolute minimum
exists in at least one of these components. By our method, we find
multiple solutions and non-minimal ones, as stated in Theorems 1-3,
in $\mathcal{A}_{+1}(A_0)$ for small values of the coupling constant
$\epsilon>0$ 
\footnote{One may also investigate these issues in
different connected components, i.e. for $k\neq \pm
1$, and 
in particular, in the connected component $\mathcal{A}_0(A_0)$,
where an absolute minimum $\underline{A}_{\epsilon}$ is already
known to exist. The techniques used in this paper extend to all the
cases (cf. Conjecture 8.1 in \cite{IM2}). Nevertheless, these
non-uniqueness problems via the perturbation approach are
technically harder for $k\neq \pm 1$.}.  It is important to point
out that, by the argumentation in $\S 2.2$, the isomorphism
$\phi_\epsilon$ establishes a correspondence between solutions to
$\bigl(\mathcal{D}_{\epsilon}\bigr)$ and solutions to
$(\mathcal{D}_1)$ (the standard Dirichlet problem) with boundary
value $\epsilon A_0$.

\medskip
In order to state our main theorems, we need to introduce further
notation
 used throughout this paper.
 We denote by $\HH$ the algebra of quaternions, i.e., $\HH$ is the
associative algebra over $\R$ generated by $i,j,k$, which satisfy
$i^2=j^2=k^2=-1$, $ij=-ji=k$, $jk=-kj=i$ and $ki=-ik=j$. Thus,
$x\in\HH$ is written as $x=x^0+x^1i+x^2j+x^3k$, $x_i\in\R$ ($1\le
i\le 4$). The real and the imaginary components of $x$ are
$\text{Re}\,x:=x^0$ and $\text{Im}\,x:=x^1i+x^2j+x^3k$,
respectively. The inner product of $x,y\in\HH$ is $(x,y)
:=\text{Re}\, (x\overline{y})$, where
$\overline{y}=y^0-y^1i-y^2j-y^3k$. The Lie algebra,
$\text{Im}\,\HH$, of imaginary quaternions, with Lie bracket $[x,y]
:= xy - yx$, is isomorphic to $\mathfrak{su}(2)$, and the Lie group
$Sp(1)$ of unit quaternions is isomorphic to $SU(2)$. An isomorphism
between $\text{Im}\,\HH$ and $\text{su}(2)$ is given explicitly by
$\text{Im}\,\HH\ni x^1i+x^2j+x^3k\mapsto \begin{pmatrix} x^1i&x^2+x^3i\\
-x^2+x^3i&-x^1i\end{pmatrix}\in\mathfrak{su}(2)$. We endow
$\mathfrak{su}(2)$ with the inner product $(X,Y)=-\tr(XY)$, which
translates in terms of quaternions into $(X,Y)=2(x,y)$, where
$X,Y\in\mathfrak{su}(2)$ correspond to $x,y\in\text{Im}\,\HH$ via
the above isomorphism. The pointwise inner product on
$\mathfrak{su}(2)$-valued forms on $B^4$ is defined via the inner
product on $\mathfrak{su}(2)$ and the standard metric on $B^4$. In
the following, for any $\text{Im}\,\HH$-valued $q$-form $\omega$, we
denote by $\omega_1,\omega_2,\omega_3$ its real-valued components,
i.e., $\omega=\omega_1i+\omega_2j+\omega_3k$, where the $q$-forms
$\omega_1$, $\omega_2$, $\omega_3$ are real-valued. For $p\in B^4$,
we define $h_p$ as the $\text{Im}\,\HH$-valued $1$-form which is the
unique solution of
$$\left\{\begin{array}{ll}
\Delta h_p=0\quad&\mbox{in } B^4\\
h_p=\text{Im}\frac{(\overline{x}-\overline{p})dx}{|x-p|^4}\quad&\mbox{at
} \partial B^4. \end{array}\right.$$ Here, all the components of
$h_p$ (not only the tangential ones) are prescribed at the boundary,
so the one-form $h_p$ is harmonic component-wise, with assigned
Dirichlet boundary data. This one-form, as well as the function
$F(p)$ and the matrix $M(A_0, p)$ defined below play a crucial role
in this paper. We define the function
$$F(p):=\int_{B^4}|(dh_p)^-|^2\,dx\;,\; p\in B^4\;,$$
\noindent and, for a given boundary value $A_0$, and $p\in B^4$, we
define the $3\times 3$-matrix
$$M(A_0,p):=\bigl(m_{ij}(A_0,p)\bigr)\;,\; 1\le i,j\le 3\,,$$
with
$$m_{ij}(A_0,p):=\int_{B^4}\bigl((d\underline{A}_{0,j})^-,(dh_{p,i})^-\bigr) \,dx\;,\; 1\le i,j\le 3\,,$$
where $\underline{A}_0$ is a solution to
$$\bigl(\mathcal{D}_0\bigr)\qquad\quad\left\{\begin{array}{ll}
d^{\ast}dA=0\quad&\mbox{in } B^4\\
\iota^{\ast}A\sim A_0\quad&\mbox{on } \partial B^4\;,
\end{array}\right.$$
and, for a given $2$-form $\omega$, we denote by $\omega^-$ the
anti-self dual component of $\omega$ (that is
$\omega^-:=(\omega-\ast\omega)/2$).
\newtheorem{remark}{Remark}[section]
\begin{remark}
\label{R2.1} Note that $d\underline{A}_0$ is uniquely determined by
(the gauge equivalence class of) $A_0$, since it satisfies the
system of equations $d^2\underline{A}_0=0$,
$d^{\ast}d\underline{A}_0=0$, $\iota^{\ast}d\underline{A}_0=d_\tau
A_0$, which has a unique solution by relative Hodge theory. Thus the
matrix above is a well defined matrix-valued function of $A_0$ and
$p$.
\end{remark}
We denote by $\mu_1(A_0,p)\ge\mu_2(A_0,p)\ge\mu_3(A_0,p)\ge 0$ the
eigenvalues of the non-negative symmetric matrix $M(A_0,p)^t
M(A_0,p)$, and state Theorems 1-3 for solutions in
$\mathcal{A}_{+1}(A_0)$. (Analogous results would hold for solutions
in $\mathcal{A}_{-1}(A_0)$ by the same arguments, or by simply
reversing the orientation of $B^4$).

\newtheorem{theorem}{Theorem}
\begin{theorem} Let us define the function $G_1^\pm (p):=\frac{(\sqrt{\mu_1(p)}+ \sqrt{\mu_2(p)}\pm\sqrt{\mu_3(p)})^2}{F(p)}\;,\; p\in B^4,$
and assume that $p_0\in B^4$ satisfies either of the following
hypotheses (1),(2):
\begin{enumerate}[(1)]
\item $\det M(A_0,p_0)>0$  and $p_0$ is an isolated local maximum point of $G_1^+(p);$
\item $\det M(A_0,p_0)<0$ and $p_0$ is an isolated local maximum point of $G_1^-(p).$
\end{enumerate}
Then, there exists $\epsilon_0>0$ and  a family of connections
$\{{A}_\epsilon\}$ indexed by $\epsilon\in (0, \epsilon_0]$ with the
following properties: ${A}_{\epsilon}$ is a solution to
$\bigl(\mathcal{D}_{\epsilon}\bigr)$ in $\mathcal{A}_{+1}(A_0)$;
$\epsilon^2|{F_{{A}_{\epsilon}}}^{\epsilon}|^2\,dx\to
8\pi^2\delta_{p_0}$ as $\epsilon\to 0$ in the sense of measures
(i.e. $\epsilon {F_{{A}_{\epsilon}}}^{\epsilon}$ concentrates at
$p_0$ as $\epsilon\to 0$).
%
\end{theorem}
\begin{remark}
Note that the ${{A}_{\epsilon}}$'s above do not satisfy the
Yang Mills equations, but the $SU(2)_\epsilon$-Yang Mills equations
(cf. $\mathcal D_\epsilon$). Nonetheless the thesis of Theorem 1 can
be restated by means of the isomorphism $\Phi_\epsilon$ as follows:

\noindent $<<$There exists $\epsilon_0>0$ and a family of Yang Mills
connections  $\{{A}(\epsilon)\}$, indexed by $\epsilon\in (0,
\epsilon_0]$, with the following properties:
 ${A}(\epsilon)$ is a solution to the Dirichlet problem $\mathcal D_1$ with boundary value $\epsilon A_0$ in
$\mathcal{A}_{+1}(\epsilon A_0)$;
$|F_{{A}(\epsilon)}|^2\,dx\to 8\pi^2\delta_{p_0}$ as
$\epsilon\to 0$ in the sense of measures (i.e. the
$F_{{A}(\epsilon)}$ concentrates at $p_0$ as $\epsilon\to
0)$.$>>$

\noindent Note that there holds the relation
${A}(\epsilon)= \epsilon {A}_\epsilon$ (cf. also
$\S2.2)$. The theses of Theorems 2, 3 can be restated in similar
fashion in terms of the ${A}(\epsilon)$.
\end{remark}
In Lemma 2.1 (1) of \cite{IM2} we show that $F(p)>0$ in $B^4.$ Also,
one of the basic properties of $F(p)$, namely
$\text{dist}(p,\partial B^4)^4F(p)\to C$ for some constant $C>0$ as
$\text{dist}(p,\partial B^4)\to 0$ (cf.  Lemma 2.1 (2) of
\cite{IM2}), implies that the maximum of $G_1^{\pm}$ is always
attained at some point in $B^4$ (provided that $G^\pm_1\neq 0).$
Notice that the solution obtained in~\cite{IM} corresponds to the
global maximum of $G_1^{\pm}$. (Here, we point out that a different
sign convention is used in \cite{IM}: in the main theorem
of~\cite{IM}, ``self dual" should be replaced by ``anti-self dual"
and viceversa, and the glued connection in that proof is in
$\mathcal{A}_{-1}(A_0)$ (not in $\mathcal{A}_{+1}(A_0)$),
accordingly to our current convention).

\begin{theorem}
Let us define the functions
$G_2^\pm(p):=\frac{(\sqrt{\mu_1(p)}-\sqrt{\mu_2(p)}\mp\sqrt{\mu_3(p)})^2}{F(p)};$
\hfill

\noindent
$G_3^-(p):=\frac{(-\sqrt{\mu_1(p)}+\sqrt{\mu_2(p)}+\sqrt{\mu_3(p)})^2}{F(p)};$
$G_1^0(p):=\frac{(\sqrt{\mu_1(p)}+\sqrt{\mu_2(p)})^2}{F(p)};$
$G_2^0(p):=\frac{(\sqrt{\mu_1(p)}-\sqrt{\mu_2(p)})^2}{F(p)}.$

\noindent Assume that $p_0\in B^4$ satisfies one of the following
conditions (1)-(a),(b), (2)-(a),(b),(c), (3)-(a),(b):
\begin{enumerate}[(1)]
\item $\det M(A_0,p_0)>0$ and
\begin{enumerate}
\item $p_0$ is a non-degenerate critical point of $G_1^+(p)$, or
\item $\sqrt{\mu_1(A_0,p_0)}>\sqrt{\mu_2(A_0,p_0)}+\sqrt{\mu_3(A_0,p_0)}$ and $p_0$ is a non-degenerate critical point of $G_2^+(p);$
\end{enumerate}
\item $\det M(A_0,p_0)<0$ and
\begin{enumerate}
\item $\mu_2(A_0,p_0)>\mu_3(A_0,p_0)$ and $p_0$ is a non-degenerate critical point of
$G_1^-(p)$, or
\item $\mu_1(A_0,p_0)>\mu_2(A_0,p_0)>\mu_3(A_0,p_0)$ and $p_0$ is a non-degenerate critical point of
$G_2^-(p)$, or
\item $\mu_1(A_0,p_0)>\mu_2(A_0,p_0)$, $\sqrt{\mu_1(A_0,p_0)}<\sqrt{\mu_2(A_0,p_0)}+\sqrt{\mu_3(A_0,p_0)}$ and
$p_0$ is a non-degenerate critical point of $G_3^-(p);$
\end{enumerate}
\item $\det M(A_0,p_0)=0$ and
\begin{enumerate}
\item $\mu_2(a_0,p_0)>0$ and $p_0$ is a non-degenerate critical point of
$G_1^0(p)$, or
\item $\mu_1(A_0,p_0)>\mu_2(A_0,p_0)>0$ and $p_0$ is a non-degenerate critical point of
$G_2^0(p).$
\end{enumerate}
\end{enumerate}
Then, there exists $\epsilon_0>0$ and  a family of connections
$\{{A}_\epsilon\}$ indexed by $\epsilon\in (0, \epsilon_0]$ with the
following properties: ${A}_{\epsilon}$ is a solution to
$\bigl(\mathcal{D}_{\epsilon}\bigr)$ in $\mathcal{A}_{+1}(A_0)$;
$\epsilon^2|{F_{{A}_{\epsilon}}}^{\epsilon}|^2\,dx\to
8\pi^2\delta_{p_0}$ as $\epsilon\to 0$ in the sense of measures
(i.e. $\epsilon {F_{{A}_{\epsilon}}}^{\epsilon}$ concentrates at
$p_0$ as $\epsilon\to 0$).
\end{theorem}
One of the main properties of $F(p)$, namely $\nabla F(p)\sim
\text{dist}(p,\partial B^4)^{-5}p/|p|$, as $\text{dist}(p,\partial
B^4)\to 0$ (cf. Lemma 2.1 (3) in \cite{IM2}), implies that, if the
hypotheses of the theorem above are satisfied, the derivatives of
the $G_i^{\pm}$'s, $G_i^0$'s  do not vanish at $\partial B^4$. For
each $G_i^{\pm}$, $G_i^0$, there always exists at least one critical
point in $B^4$, a maximum indeed. However, these critical points may
be degenerate.

The next theorem holds without assuming such non-degeneracy.
\begin{theorem}
Assume that there exists $p_0\in B^4$ such that one of the following
holds:

\begin{enumerate}[(1)]
\item   $\det M(A_0,p_0)>0$,
$\mu_1(A_0,p_0)>\mu_2(A_0,p_0)>\mu_3(A_0,p_0)$ and
$\sqrt{\mu_1(A_0,p_0)}>\sqrt{\mu_2(A_0,p_0)}+\sqrt{\mu_3(A_0,p_0)}$;
\item $\det M(A_0,p_0)<0$ and $\mu_1(A_0,p_0)>\mu_2(A_0,p_0)>\mu_3(A_0,p_0)$;
\item $\det M(A_0,p_0)=0$ and $\mu_1(A_0,p_0)>\mu_2(A_0,p_0)>0$;
\end{enumerate}
then, for all sufficiently small $\epsilon>0$, there exist at least
two distinct solutions to $\bigl(\mathcal{D}_{\epsilon}\bigr)$ in
$\mathcal{A}_{+1}(A_0).$ Furthermore, the following alternative
holds: there exists at least one non-minimizing solution, or there
exist infinitely many minimizing solutions. In the hypotheses (2),
if in addition
$\sqrt{\mu_1(A_0,p_0)}<\sqrt{\mu_2(A_0,p_0)}+\sqrt{\mu_3(A_0,p_0)},$
then  there exist at least three distinct solutions, of which at
least two non-minimizing, or there exist infinitely many minimizing
solutions to $\bigl(\mathcal{D}_{\epsilon}\bigr)$ in
$\mathcal{A}_{+1}(A_0)$.
\end{theorem}

\medskip
Let us point out to the reader's attention some important results
obtained in \cite{IM2}, regarding the hypotheses of Theorems 1-3.
Precisely, we construct a family of boundary values yielding
matrices $M(A_0,p_0)$ which realize each of the cases in Theorem 3
for any given point $p_0\in B^4$ (cf. Proposition 8.1 in
\cite{IM2}). We also show that for any boundary value $A_0$, there
exists an arbitrarily small perturbation $\tilde{A}_0$ of $A_0$ such
that $\det M(\tilde{A}_0,p_0)\ne 0$ and
$\mu_1(\tilde{A}_0,p_0)>\mu_2(\tilde{A}_0,p_0)>\mu_3(\tilde{A}_0,p_0)$
(cf. Proposition 8.2 in \cite{IM2}).

\medskip
\noindent In $\S2$ of the present paper, we describe the asymptotic
profile of small and large solutions as $\epsilon\to 0$, thus giving
a heuristic explanation of our method. In $\S3.1$ we construct the
spaces of approximate solutions and introduce the technical notation
used in the estimates that follow. In $\S3.2$ we obtain the
asymptotic expansion of the $\mathfrak{su}(2)_\epsilon$-Yang Mills
functional evaluated on the approximate solutions, and in
$\S3.3-\S3.4$ we estimate the Hessian and the remainder. In $\S3.5$
we introduce and estimate the modified Hessian. In $\S3.6$ we define
the auxiliary equation and solve it.

\section{Asymptotic profile of small and large solutions as $\epsilon\to 0$}
In this section, we analyze the asymptotic behavior as $\epsilon\to
0$ of the family $\{\underline{A}_{\epsilon}\}$ of \emph{small
solutions} to the Dirichlet problems $\bigl(\mathcal{D}_\epsilon
\bigr)$ defined in $\S$1 (i.e. absolute minimizers for the
$SU(2)_\epsilon$-Yang Mills functionals).
 This is a crucial ingredient in the proofs of
Theorems 1-3. We also describe the asymptotic profile as
$\epsilon\to 0$ of the family of \emph{large solutions}
$\{\overline{A}_{\epsilon}\}\subset\mathcal{A}_{+1}(A_0)$ (or
$\mathcal{A}_{-1}(A_0)$) obtained in~\cite{IM}, in order to give a
heuristic argumentation that motivates the procedure we employ to
construct approximate solutions to
$\bigl(\mathcal{D}_{\epsilon}\bigr)$.

\subsection{Asymptotic profile of small solutions}
Let $\bigl(\mathcal{D}_0\bigr)$ be the  Dirichlet problem defined in
$\S$1. The following Proposition holds for the family of
$\{\underline{A}_{\epsilon}\}$, \emph{small solutions} to
$\bigl(\mathcal{D}_\epsilon\bigr)$  .
\newtheorem{proposition}{Proposition}[section]
\begin{proposition}
\label{P2.1} There exists a solution $\underline{A}_0$ to
$\bigl(\mathcal{D}_0\bigr)$ such that
$\underline{A}_{\epsilon}\to\underline{A}_0$ in
$C^{\infty}(\overline{B}^4)$, as $\epsilon\to 0$, in a suitable
gauge. More precisely, for any $k\ge 1$ there exists $C_k>0$ such
that
$\|\underline{A}_{\epsilon}-\underline{A}_0\|_{C^k(\overline{B}^4)}\le
C_k\epsilon$, for small $\epsilon>0$, in a suitable gauge.
\end{proposition}

\textit{Proof.} We define 1-forms $\omega_\epsilon :=
 {\underline  A}_\epsilon - {\underline A}_0.$ These satisfy the following
boundary value
 problems:

\begin{equation}
\label{omega} \left \{
\begin{array}{ll}
d^*d \omega_\epsilon +\epsilon\{ \frac{d^*}{2}[{\underline
A}_\epsilon, {\underline  A}_\epsilon]+* [{\underline  A}_\epsilon,
*d {\underline  A}_\epsilon] +*\frac {\epsilon}
{2}[{\underline  A}_\epsilon,* [{\underline  A}_\epsilon,
{\underline  A}_\epsilon]]\}  =0 & \;\text{on$\;B^4$}\\
i^*\omega_\epsilon= 0& \;\text{on$\;\partial B^4\;.$}
\end{array}
\right .
\end{equation}

One needs the following lemma.
\newtheorem{lemma}{Lemma}[section]
\begin{lemma}
\label{L2.1} There exist constants $B_0$, $B_1$, $C$, $\epsilon_0$
such that
\begin{align}
&\Vert {\underline A}_\epsilon\Vert _{L^2_1(B^4)}\leq C\Vert
{F_{\underline{A}_{\epsilon}}}^{\epsilon}\Vert_{L^2(B^4)}
\leq B_0,\\
&\Vert{\underline  A}_\epsilon\Vert _{L^p_1(B^4)}\leq C\Vert
{F_{\underline{A}_{\epsilon}}}^{\epsilon}\Vert_{L^p(B^4)},\\
&\Vert \omega_\epsilon\Vert _{L^2_1(B^4)}\leq \epsilon B_1,
\end{align}
for $2\leq p<4,$ for all $\epsilon$ such that $0\leq\epsilon\leq
\epsilon_0,$ where $L^p_k$ is the Sobolev space of functions with
$L^p$-integrable partial derivatives up to order $k$ (in the sense
of distributions).
\end{lemma}
\textit{Proof.} We denote by $\mathcal{YM}_\epsilon(A)$ the Yang
Mills functional on $\mathfrak{su}(2)_\epsilon$-connections $A$,
i.e., explicitly, in local coordinates:
\begin{equation} \label{YMepsilon}\mathcal{YM}_\epsilon(A):= \int_{B^4}\vert{F_A}^{\epsilon}\vert := \int_{B^4}\Big\vert dA
+ \frac{1}{2}[A, A]_\epsilon\Big\vert^2 = \int_{B^4}\Big\vert dA +
\frac{\epsilon}{2}[A, A]\Big\vert^2\;.\end{equation}

\noindent We first show that $\mathcal{YM}_\epsilon({\underline
A}_\epsilon)$ is uniformly bounded for $\epsilon$ sufficiently
small. In fact,
\begin{equation}
 \mathcal{YM}_\epsilon({\underline  A}_\epsilon)= m_\epsilon:= \min_{A\in \mathcal{A}}\, \mathcal{YM}_\epsilon
 (A)\,,\,
\text{where}~\mathcal{A}:= \{\text{smooth
$\mathfrak{su}(2)$-connections}\, : i^* A= A_0\}\,.
\end{equation}
Thus,
\begin{align}\mathcal{YM}_\epsilon({\underline
A}_\epsilon)&\leq \mathcal{YM}_\epsilon ({\underline   A}_0) = \Vert
d{\underline A}_0\Vert^2_{L^2(B^4)} +\epsilon\langle d{\underline
A}_0, [{\underline   A}_0, {\underline   A}_0]\rangle+\frac
{\epsilon^2}{4}
\Vert [{\underline A}_0, {\underline   A}_0]\Vert ^2_{L^2(B^4)}\notag\\
&=\mathcal {YM}_0({\underline   A}_0)+\epsilon\langle d{\underline
A}_0, [{\underline   A}_0, {\underline   A}_0]\rangle+\frac
{\epsilon^2}{4} \Vert [{\underline   A}_0, {\underline   A}_0]\Vert
^2_{L^2(B^4)}\leq 2 m_0\;,
\end{align}
for $\epsilon$ sufficiently small, where $m_0:=\mathcal
{YM}_0({\underline   A}_0).$

 Now, let us recall that the absolute Yang Mills minimizer for the Dirichlet problem found
 in~\cite{Marini} is in the good gauge, i.e., it satisfies $d^*A=0$ on $B^4$ and the boundary condition $d^*_\tau A_\tau=0$ at $\partial B^4$, where $\tau$ represents tangential
 directions. These conditions yield the estimate
 \begin{equation}
 \label{estimate}
 \Vert {\underline   A}_\epsilon\Vert _{L^p_1(B^4)}\leq h \Vert d{\underline  A}_\epsilon\Vert _{L^p(B^4)},
 \end{equation}
 for $\epsilon\geq 0,$ and $2\leq p <4,$ with a constant $h$ depending
 only on dimension (not on $\epsilon$).
 For $\epsilon=0,$ this becomes $\Vert {\underline A}_0\Vert _{L^p_1(B^4)}\leq h \Vert {F_{\underline{A}_0}}^0\Vert_{L^p(B^4)},$ i.e., (2.3).
 This estimates also holds on local charts $U$ of type one and type two (boundary and interior neighborhoods, respectively, cf.~\cite{Marini}).
 This yields
 \begin{equation}
 \label{onlocalcharts}
 \Vert {\underline  A}_\epsilon\Vert_{L^2_1(U)}\leq
 h\Vert{F_{\underline{A}_{\epsilon}}}^{\epsilon}\Vert_{L^2(U)} +
 hC^\prime\epsilon\Vert {\underline  A}_\epsilon\Vert_{L^4(U)}\Vert {\underline  A}_\epsilon\Vert_{L^2_1(U)}\;.
 \end{equation}
 Let us now consider covers ${\mathcal C}_j$ of charts $U$ of
 radius $\rho_j$, with $\rho_j\to 0$ as $j\to\infty$, satisfying the following conditions:

 \noindent
 (1) there exists $k$, independent of $j$,
 such that any $k+1$ charts have empty intersection (in particular each cover ${\mathcal C}_j$ is
 finite);

 \noindent
 (2) $\forall \epsilon>0$ given, there exists $j=j(\epsilon)$ such that $hC^\prime\epsilon\Vert
{\underline  A}_\epsilon\Vert_{L^4(U)}<
 \frac{1}{2}$, for every $U\in{\mathcal C}_j$.

 \noindent This can be achieved using the compactness of ${\overline
 B}^4$ and using charts of type one and type two.
Then, by \eqref{onlocalcharts},
\begin{equation}
\Vert {\underline A}_\epsilon\Vert_{L^2_1(U)}\leq 2 h\Vert
{F_{\underline{A}_{\epsilon}}}^{\epsilon}\Vert_{L^2(U)}\;,\;\forall
U\in{\mathcal C}_j\;.
\end{equation}
Thus on $B^4$,
$$\Vert {\underline  A}_\epsilon\Vert _{L^2_1(B^4)}\leq
2 h
k\Vert{F_{\underline{A}_{\epsilon}}}^{\epsilon}\Vert_{L^2(B^4)}\;,$$
with $h,k$ independent of $\epsilon$. This, together with (2.7)
yields (2.2) with $C=2hk$ and $B_0= 2Cm_0.$

\noindent
 For $\epsilon >0 $ and general $p\in [2, 4)$, (2.8) yields
$$\Vert {\underline  A}_\epsilon\Vert _{L^p_1(B^4)}\leq h\Vert{F_{\underline{A}_{\epsilon}}}^{\epsilon}\Vert_{L^p(B^4)} +
 hC^\prime\epsilon\Vert {\underline  A}_\epsilon\Vert_{L^4(B^4)}\Vert {\underline
 A}_\epsilon\Vert_{L^q(B^4)}\leq
 h\Vert{F_{\underline{A}_{\epsilon}}}^{\epsilon}\Vert_{L^p(B^4)} +
 hC^\prime\epsilon\Vert {\underline  A}_\epsilon\Vert_{L^2_1(B^4)}\Vert {\underline  A}_\epsilon\Vert_{L^p_1(B^4)}\;,$$
 for $q=\frac {4p}{4-p}.$
 Thus, applying (2.2), there exists $\epsilon_0>0$ such that (2.3)
 holds in the ``good gauge", for some constant $C$ depending only on dimension, for $2\leq
 p<4$ and for $0\leq\epsilon\leq\epsilon_0$.

 \noindent
 In the same gauge, for $\omega_\epsilon :=
 {\underline  A}_\epsilon - {\underline A}_0$, one also obtains
 \begin{equation}
 \label{estimateomega1}
 \Vert \omega_\epsilon\Vert _{L^2_1(B^4)}\leq h \Vert
 d\omega_\epsilon\Vert_{L^2(B^4)}\;.
 \end{equation}
Then,  system \eqref{omega} and integration by parts yield
\begin{equation}
\label{eqdomega}
 \Vert d\omega_\epsilon\Vert^2_{L^2(B^4)} =-\int_{B^4} \omega_\epsilon\wedge *
 d^*d\omega_\epsilon + \int_{\partial B^4} i^*(\omega_\epsilon\wedge *
 d\omega_\epsilon)= \epsilon \int_{B^4} \omega_\epsilon\wedge *
 \{ \frac{d^*}{2}[{\underline
A}_\epsilon, {\underline  A}_\epsilon]+* [{\underline  A}_\epsilon,
*{F_{\underline{A}_{\epsilon}}}^{\epsilon} \}.
\end{equation}
Sobolev embeddings and H\"older inequalities give
$$
 \Vert d\omega_\epsilon\Vert^2_{L^2(B^4)}\leq \epsilon C^\prime \Vert{\underline
 A}_\epsilon\Vert^2_{L^2_1(B^4)}\Vert{F_{\underline{A}_{\epsilon}}}^{\epsilon}\Vert^2_{L^2(B^4)} \;.$$
 Thus,
 \begin{equation}
\label{estimateomega2} \Vert \omega_\epsilon\Vert_{L^2_1(B^4)}\leq h
\Vert
 d\omega_\epsilon\Vert_{L^2(B^4)}\leq\sqrt\epsilon B^\prime\;,
 \end{equation}
 for some constant $B^\prime,$ if $\epsilon$ is sufficiently small.
 Using estimate \eqref{estimateomega2} into \eqref{eqdomega} and
 estimating again, one obtains (2.4),
for some constant $B_1,$ if $\epsilon$ is sufficiently small.

This completes the proof of the lemma.\hfill$\Box$

\medskip

To prove the proposition, we cover $B^4$ with coordinate patches
$U_1:=\{x\in B^4;\ \vert x\vert^2 <\delta\}$ of ``type one", in the
interior, and $U_2:=\{(x^\prime, x^4): x^\prime\in\partial B^4;\
x^4\geq 0;\ \vert x\vert^2< \delta\}$ of ``type two" near the
boundary. Here, the functions $x^4\to(x^\prime, x^4)$ describe unit
speed geodesics orthogonal to $\partial B^4.$ This way, the metric
$g_{ij}$ satisfies $g_{i4}(x^\prime,0)=0$ and $g_{44}= 1$ in a
neighborhood of the boundary. Doubling $B^4$ via reflection across
the boundary, yields a Lipschitz-bounded metric on the resulting
manifold. We lift this action trivially to the bundle. We show that
$\Vert \omega_\epsilon\Vert_{C^k(U_2)}\leq C_k\epsilon$ holds for
small $\epsilon>0$, in a suitable gauge, all the way up to $\partial
B^4,$ on neighborhoods of type two. (We omit the proof for
neighborhoods of type one).

Locally, on $U_2$, system \eqref{omega} and the good gauge theorem
for boundary neighborhoods (cf. \cite{Marini}) yield
\begin{equation}
\label{omegaloc} \left \{
\begin{array}{ll}
d^*d \omega_\epsilon +\epsilon\{ \frac{d^*}{2}[{\underline
A}_\epsilon, {\underline  A}_\epsilon]+* [{\underline  A}_\epsilon,
*d {\underline  A}_\epsilon] +*\frac {\epsilon}
{2}[{\underline  A}_\epsilon,* [{\underline  A}_\epsilon,
{\underline  A}_\epsilon]]\}  =0 & \;\text{on}\;U_2\\
d^{*_F} \omega_\epsilon =0 & \;\text{on}\;U_2\\
d^{*_F}_\tau (\omega_\epsilon)_\tau =0 & \;\text{on}\;\{x^4=0\}\\
i^*\omega_\epsilon= 0& \;\text{on}\;\{x^4=0\}\;,
\end{array}
\right .
\end{equation}
where $\tau$ denotes tangential components and $*_F$ is the flat
Hodge operator.

\noindent This becomes
\begin{equation}
\label{omegaloc1} \left \{
\begin{array}{ll}
\Lambda \omega_\epsilon:= \Delta_F \omega_\epsilon + \mathcal{E}
\omega_\epsilon +\epsilon
R_\epsilon \omega_\epsilon  =\epsilon G(\underline{A}_\epsilon, d\underline{A}_\epsilon) & \;\text{on}\;U_2\\
d^{*_F} \omega_\epsilon =0 & \;\text{on}\;U_2\\
d^{*_F}_\tau (\omega_\epsilon)_\tau =0 & \;\text{on}\;\{x^4=0\}\\
i^*\omega_\epsilon= 0& \;\text{on}\;\{x^4=0\}\;,
\end{array}
\right .
\end{equation}
where $\mathcal{E}=*_F d(*-*_F)d$ contains only first order
derivatives of the metric and can be made small by dilations,
$R_\epsilon (\cdot) =
*_F \{ \frac{d*}{2}[\cdot, \underline{A}_\epsilon]+ [\cdot,
*d\underline{A}_\epsilon] +\frac {\epsilon}
{2}[\cdot,* [\underline{A}_\epsilon,\underline{A}_\epsilon]]\}$ is
of lower order, and $G(\underline{A}_\epsilon,
d\underline{A}_\epsilon)= -*_F \{ \frac{d*}{2}[\underline{A}_0,
\underline{A}_\epsilon]+ [\underline{A}_0,
*d\underline{A}_\epsilon] +\frac {\epsilon}
{2}[\underline{A}_0,*
[\underline{A}_\epsilon,\underline{A}_\epsilon]]\}$ is uniformly
bounded in $L^2$ for small $\epsilon$ by the previous lemma.

Following the procedure in \cite{Marini}, we reflect $U_2$ across
the boundary $\{x^4=0\}$ and work on the doubled neighborhood
$\tilde U$, after doubling all the operators above via the formula
$r^*\tilde \Lambda = \tilde \Lambda r^*.$ More in detail,
$$\tilde\Lambda(\omega)(x)=\Lambda(\omega\vert_{U^+})\chi_{U^+}(x)
+r^*\Lambda(r^*(\omega\vert_{U^-}))\chi_{U^-}(x),$$ (here
$U^\pm=\{x^4
>0 (<0)\}$), for all 1-forms $\omega$. An easy computation shows that
the double $\tilde\Delta_F$ is $\Delta_F$ on $\tilde U.$

Moreover, $\tilde{\mathcal{E}}$ and $\epsilon \tilde{R_\epsilon}$
are small operators from $L^p_{0,1}(T^{\ast}\tilde
{U}\otimes\mathfrak{su}(2))$ to
$L^p_{-1}(T^{\ast}\tilde{U})\otimes\mathfrak{su}(2))$ and, also,
from $L^p_2(T^{\ast}\tilde {U}) \otimes\mathfrak{su}(2))$ to
$L^p(T^{\ast}\tilde {U}\otimes\mathfrak{su}(2))$ for $p>1$, where
$L^p_{0,1}(T^{\ast}\tilde{U}\otimes\mathfrak{su}(2))$ is the
completion of
$C^{\infty}_0(T^{\ast}\tilde{U}\otimes\mathfrak{su}(2))$ with
respect to the $L^p_1$-norm, see also $\S3.3$.

We take care of the boundary conditions on $\partial \tilde U$ by
introducing a smooth appropriate cut-off function $\phi.$ System
\eqref{omegaloc1} then becomes
\begin{equation}
\label{omegaloc2} \left \{
\begin{array}{ll}
\Delta_F
(\phi{\tilde\omega}_\epsilon)+\mathcal{E}(\phi{\tilde\omega}_\epsilon)+\epsilon
R_\epsilon(\phi{\tilde\omega}_\epsilon)=\epsilon
(\tilde\phi G(\underline{A}_\epsilon, d\underline{A}_\epsilon)+\frac {1}{\epsilon} T_\phi{\tilde\omega}_\epsilon):= \epsilon \alpha_\epsilon & \;\text{on}\;\tilde U\\
\phi ({\tilde\omega}_\epsilon)_\tau =0 & \;\text{on}\;\partial
\tilde U\;,
\end{array}
\right .
\end{equation}
where $\tilde \omega$ and $\tilde G$ are odd extensions of $\omega$
and $G$ (i.e., $r^* \tilde \omega = - \tilde \omega,$ and $r^*
\tilde G = - \tilde G$), and  $T_\phi{\tilde\omega}_\epsilon =
\Lambda (\phi {\tilde \omega}_\epsilon) - \phi \Lambda (
{\tilde\omega}_\epsilon)$ contains only first order derivatives of
${\tilde\omega}_\epsilon.$ With this definition, $\tilde G$ is
uniformly bounded in $L^2,$ thus in $L^p_{-1}$ with $p\leq 4,$ and
so is $\frac {1}{\epsilon} T_\phi{\tilde\omega}_\epsilon$ (by the
estimate (2.4)). So $\alpha_\epsilon$ is uniformly bounded in
$L^p_{-1}$ with $p\leq 4$ for $0\leq\epsilon\leq\epsilon_0.$

It is well known that the system
\begin{displaymath}
\left \{
\begin{array}{ll}
\Delta \omega= \gamma
& \;{\mbox on}\; \Vert x\Vert\leq \delta\\
\omega =0 & \;{\mbox on }\;\Vert x\Vert =\delta\;,
\end{array}
\right .
\end{displaymath}
with $\gamma\in L^p_{-1}$ and $\omega\in L^2_1,$  admits a unique
solution in $L^p_1.$ Let $\mathcal S$ be the solution operator
(bounded). Applying $\mathcal S$ to \eqref{omegaloc2} one obtains
$$I(\phi{\tilde\omega}_\epsilon) + \mathcal S[\mathcal{E}
(\phi{\tilde\omega}_\epsilon) +\epsilon R_\epsilon
(\phi{\tilde\omega}_\epsilon)] =\epsilon \mathcal
S(\alpha_\epsilon).$$ Thus, since $\mathcal{E} +\epsilon
R_\epsilon:L^p_{0,1}\to L^p_{-1}$ is small, we can invert $I+
\mathcal S (\mathcal{E} +\epsilon R_\epsilon)$ and
 obtain
$$\phi{\tilde\omega}_\epsilon = \epsilon [I+ \mathcal S (\mathcal{E} +\epsilon R_\epsilon)]^{-1}\mathcal S \alpha_\epsilon.$$
 Thus
 $$\Vert\phi {\tilde\omega}_\epsilon\Vert_{L^p_1(\tilde U)} \leq \epsilon
 \Vert[I+ \mathcal S (\mathcal{E} +\epsilon R_\epsilon)]^{-1}\mathcal S\Vert \Vert \alpha_\epsilon\Vert_{L^p_{-1}(\tilde
U)}\leq C^\prime\epsilon,$$ yielding $\epsilon^{-1}\omega_\epsilon$
uniformly bounded in $L^p_1(U_2),$
 on a smaller neighborhood $U_2,$ all the way up to the
 boundary $\{x^4=0\}.$
 Iterating the procedure above, recalling that $\underline{A}_\epsilon =\underline{A}_0 +
\omega_\epsilon$ and estimate (2.4), one obtains a system similar to
 \eqref{omegaloc2}, but simpler (this time there is
 no need for the operator $R_\epsilon$), with the right hand side
 uniformly bounded in $L^p$ for any $p<4$, yielding finally
 $$\Vert \omega_\epsilon\Vert_{\mathcal{C}_0}\leq C \Vert \omega_\epsilon\Vert_{L^q_1(U_2)}\leq C C^\prime
 \epsilon:= C_0\epsilon\;,$$
 for some $q>4,$ on a smaller neighborhood $U_2,$ all the way up to and including the
 boundary $\{x^4=0\}.$
To show the analogous result for $\nabla \omega_\epsilon,$ we take
first tangential derivatives in \eqref{omegaloc1} and proceed with
the doubling procedure above (using that $\partial_j \omega_\epsilon
= 0$ at $\{x^4=0\},$ for $j=1,2,3,$ thus the one forms $\widetilde
{\partial_j\omega}$ are continuous). For normal components
$\partial_4 \omega_\epsilon,$ one uses the relations between
tangential and normal derivatives given by the good gauge and the
field equations. Iterating this procedure, after some calculation,
one obtains $\|\underline {A}_{\epsilon}-\underline
{A}_0\|_{C^k(\overline{B^4})}:=\Vert \omega_\epsilon
\Vert_{C^k(\overline{B^4})}\le C_k\epsilon$ for small $\epsilon>0$,
in a suitable gauge. \hfill$\Box$

\subsection{Asymptotic profile of large solutions}

Here we give a heuristic argumentation to motivate our approach to
the existence of new solutions to
$\bigl(\mathcal{D}_\epsilon\bigr)$. We start by observing that the
Lie algebra isomorphism $\phi_{\epsilon}$ (and the Lie group
isomorphism $\Phi_{\epsilon}$) transform the $SU(2)_\epsilon$-Yang
Mills Dirichlet problem $\bigl(\mathcal{D}_{\epsilon}\bigr)$ for
$A$, into the $SU(2)$-Yang Mills Dirichlet problem for
$\phi_{\epsilon}(A):=\epsilon A$
\begin{equation}
\label{transformed}
\bigl(\mathcal{D}(\epsilon)\bigr)\qquad\quad\left\{\begin{array}{ll}
d_A^{\ast}F_A=0\quad&\mbox{in } B^4\\
\iota^{\ast}A\sim \phi_{\epsilon}(A_0):=\epsilon A_0\quad&\mbox{on }
\partial B^4\;,
\end{array}\right.
\end{equation}
where $d_A^{\ast}=\ast d\ast+\ast[A,\ast\cdot]$, and
$F_A=dA+\frac{1}{2}[A,A]$.

Let us set
$\underline{A}(\epsilon):=\phi_{\epsilon}(\underline{A}_{\epsilon})$
and
$\overline{A}(\epsilon):=\phi_{\epsilon}(\overline{A}_{\epsilon})$,
where $\underline{A}_{\epsilon}$ is the absolute minimizing solution
to $\bigl(\mathcal{D}_\epsilon\bigr)$, and $\overline{A}_{\epsilon}$
is the absolute minimizer restricted to the class
$\mathcal{A}_{+1}(A_0)$ (or $\mathcal{A}_{-1}(A_0)$). Then,
$\underline{A}(\epsilon)$ is an absolute minimizing solution to
$\bigl(\mathcal{D}(\epsilon)\bigr)$, and $\overline{A}(\epsilon)$ is
a large solution to $\bigl(\mathcal{D}(\epsilon)\bigr)$, i.e, it
minimizes the Yang Mills functional in $\mathcal{A}_{+1}(\epsilon
A_0)$ (or $\mathcal{A}_{-1}(\epsilon A_0)$). Passing to the limit
$\epsilon\to 0$ formally in $\bigl(\mathcal{D}(\epsilon)\bigr)$, one
obtains the Dirichlet problem for Yang Mills connections with the
trivial boundary value. It is known (cf.~\cite{I}) that this only
admits flat solutions, therefore $\overline{A}(\epsilon)$ cannot
converge strongly since
$\overline{A}(\epsilon)\in\mathcal{A}_{\pm1}(\epsilon A_0)$ for
$\epsilon>0$. Indeed, following the proof in~\cite{IM}, one may
argue that $\int_{B^4}\,|F_{\overline{A}(\epsilon)}|^2\,dx\to
8\pi^2$ as $\epsilon\to 0$, and
$|F_{\overline{A}(\epsilon)}|^2\,dx\to 8\pi^2\delta_p\,dx$ as a
Radon measure for some $p\in \overline{B}^4$. It follows that, in a
suitable gauge, one has asymptotically
$\overline{A}(\epsilon)\approx\underline{A}(\epsilon)\#\text{($\pm1$-instanton
on $S^4$)}$, where $\#$ denotes the connected sum. In terms of
$\overline{A}_{\epsilon}$, one has asymptotically
$\overline{A}_{\epsilon}\approx\underline{A}_{\epsilon}\#\frac{1}{\epsilon}(\text{$\pm1$-instanton
on $S^4$})$.
\begin{remark}
\label{R2.2} The argumentation above would require a little extra
work to be made rigorous. In fact, the spaces
$\mathcal{A}_{\pm1}(A_0)$ do depend on $\underline A_\epsilon$, thus
on $\epsilon$. However, in this paper we construct solutions to
$(\mathcal{D}_\epsilon)$ for `fixed' small positive $\epsilon$ and
we are not concerned with this issue, nor with the issue of
constructing paths of solutions parameterized by $\epsilon$.
\end{remark}

\section{Reduction to a finite dimensional problem}
In this section we construct approximate solutions to
$(\mathcal{D}_\epsilon)$ (for small $\epsilon$) via a gluing
technique, and study the asymptotic expansion of the
$SU(2)_\epsilon$-Yang Mills functional, its gradient and its
Hessian. The approximate solutions depend on a finite-dimensional
parameter (cf. \S3.1). The space tangent to the space of approximate
solutions is a good approximation for the kernel of the Hessian,
thus it constitutes the obstruction to the direct application of the
implicit function theorem. We follow the standard procedure for this
type of problems, consisting of first solving the Yang Mills
equation orthogonally to the kernel of the Hessian by means of the
auxiliary equation introduced (and solved) in \S3.6. Thus, the
problem is transformed into a finite dimensional problem (cf. in
particular Lemma \ref{L3.8} and Proposition \ref{P3.2}).

\noindent We focus on solutions that create a 1-bubble in the limit
as $\epsilon\to 0.$
\subsection{The space of approximate solutions and introduction of the notation}
Motivated by the discussion in $\S 2.2$, we seek approximate
solutions to $\bigl(\mathcal{D}_{\epsilon}\bigr)$ in
$\mathcal{A}_{+1}(A_0)$ of the form:
$A_{\epsilon}=\underline{A}_{\epsilon}\#\frac{1}{\epsilon}\text{(1-bubble)}+a$,
where $a\in C^{\infty}(T^{\ast}B^4\otimes\text{Ad}(P(p,g,\lambda)))$
is small and satisfies $a=0$ on $\partial B^4$. (The bundle
$P(p,g,\lambda)$ will be defined soon). In this section, we
introduce all the technical notation used to prove Theorems 1-3.

We start by describing the main part of the solution
$\underline{A}_{\epsilon}\#\frac{1}{\epsilon}\text{(1-bubble)}$.

For $\lambda>0$, $p\in\R^4$, the $1$-instanton solution
$I_{\lambda,p}$ with center at $p$ and scale $\lambda$ to the Yang
Mills equation on $\R^4$
is defined by
$$I_{\lambda,p}=\left\{\begin{array}{ll}
I^2_{\lambda,p}(x)=\text{Im}\frac{\lambda^2(\overline{x}-\overline{p})dx}{|x-p|^2(\lambda^2+|x-p|^2)} &\quad\mbox{in }U_2(p)=\R^4\setminus\{p\}\\
I^1_{\lambda,p}(x)=\text{Im}\frac{(x-p)d\overline{x}}{\lambda^2+|x-p|^2}
&\quad\mbox{in } U_1(p)=\{x\in\R^4:|x-p|<1\}=B^4 (p)\;,
\end{array}\right .$$
where $x=(x^0,x^1,x^2,x^3)\in\R^4$ is identified with the quaternion
$x=x^0+x^1i+x^2j+x^3k\in\mathbb{H}$ and the transition map is
$g_{12,p}(x)=\frac{x-p}{|x-p|}$ for $x\in U_1(p)\cap U_2(p)$ (notice
that the gluing relation
$I^2_{\lambda,p}=g_{12,p}^{-1}dg_{12,p}+g_{12,p}^{-1}I_{\lambda,p}^1g_{12,p}$
holds in $U_1(p)\cap U_2(p)$).

For $p\in B^4$, $\lambda>0$, we define the $1$-form $h_{\lambda,p}:=
(h_{\lambda,p})_j\,dx^j\in
C^{\infty}(T^{\ast}B^4\otimes\mathfrak{su}(2))$, the components of
which solve the Dirichlet problems
\begin{equation}
\label{3.1} \left\{\begin{array}{ll}
\Delta (h_{\lambda,p})_j=0\quad&\mbox{in } B^4\\
(h_{\lambda,p})_j=(I^2_{\lambda,p})_j\quad&\mbox{at } \partial
B^4\;,
\end{array}\right.
\end{equation}
and set $PI^2_{\lambda,p}:=I^2_{\lambda,p}-h_{\lambda,p}$, the
projection of $I^2_{\lambda,p}$ on the Sobolev space
$L^2_{0,1}(T^{\ast}B^4\otimes\mathfrak{su}(2)):=\{a\in
L^2_1(T^{\ast}B^4\otimes\mathfrak{su}(2)):a=0\text{ on $\partial
B^4$}\}$.

We also define a cut-off function  $\beta(x)=\beta(|x|)\in
C^{\infty}_0(\R^4)$,
 such that $\beta=1$ for $|x|\le 1$, $\beta(x)=0$ for
$|x|\ge 2$ and $0\le\beta(x)\le 1$, and, for $\lambda>0$,
$p\in\R^4$, we define
$\beta_{\lambda,p}(x):=\beta(\lambda^{-1}(x-p))$.

For $d_0$ and $\lambda_0$ small fixed numbers satisfying
$0<2\lambda_0<d_0<1$, we consider the set of parameters
\begin{equation}
\label{par1} \tilde{\PP}(d_0,\lambda_0):= B^4_{1-d_0}\times
SU(2)\times(0,\lambda_0)\,.
\end{equation}
For $\q:=(p,g,\lambda)\in\tilde\PP (d_0,\lambda_0)$, we define the
connections $A(\q)$ as
\begin{equation}
\label{A(q)} A(\q)=\left\{\begin{array}{ll}
(1-\beta_{\lambda,p})\underline{A}_{\epsilon}+\frac{1}{\epsilon}\beta_{\lambda/4,p}\,g\,I^2_{\lambda,p}\,g^{-1}+
\frac{1}{\epsilon}(1-\beta_{\lambda/4,p})gPI^2_{\lambda,p}g^{-1}&\quad\mbox{in } B^4\setminus\{p\}\\
\frac{1}{\epsilon}gI^1_{\lambda,p}g^{-1}&\quad\mbox{in }
B^4_{\lambda/4}(p)\,,
\end{array}\right.
\end{equation}
obtained  by gluing the $1$-instanton to $\underline{A}_{\epsilon}$.
These connections live on the bundles $P(\q)$, defined by the data
\begin{equation}
\label{P(q)} \biggl(B^4_{\lambda/4}(p), B^4\setminus\{p\},
g\,g_{12,p}\,g^{-1}\biggr)\,.
\end{equation}
Notice that the relative 2nd Chern number of $P(\q)$ with respect to
$\underline{A}_{\epsilon}$ is $1$, thus
$A(\q)\in\mathcal{A}_{+1}(A_0)$.

\noindent We observe that the effective parameter space is
\begin{equation}
\label{par2} \PP(d_0,\lambda_0):=\tilde\PP(d_0,\lambda_0)/\{\pm 1\}
=B^4_{1-d_0}\times SO(3)\times(0,\lambda_0),
\end{equation}
since $P(p,-g,\lambda)=P(p,g,\lambda)$ and
$A(p,-g,\lambda)=A(p,g,\lambda)$, therefore from now on we quotient
out with respect to this action and redefine $\q:=(p,[g],\lambda)$,
where $[g]\in SO(3)=SU(2)/\{\pm1\}$.

\noindent We also observe that the bundles $P(\q)$, for
$\q\in\PP(d_0,\lambda_0)$ are all isomorphic, so we  fix ${\q}_0
:=(p_0,[g_0],\lambda_0)\in\PP(d_0,\lambda_0)$ and apply the
convention that everything is pulled back to
$P(p_0,[g_0],\lambda_0)$, via the bundle isomorphisms
$\varphi({\q}):P({\q}_0)\overset{\sim}{\to}P(\q)$.

\noindent We define the map
$\Gll:\PP(d_0,\lambda_0)\to\mathcal{A}_{+1}(A_0)$, $\q\mapsto
A(\q)$,  and the space
$$\mathcal{N}(d_0,\lambda_0):=\Gll(\PP(d_0,\lambda_0)),$$ as the space
of approximate solutions to $\bigl(\mathcal{D}_{\epsilon}\bigr)$ in
$\mathcal{A}_{+1}(A_0)$.


\noindent Let $\mathcal{G}(\q):= \Aut\;P(\q)$ be the space of smooth
gauge transformations of $P(\q)$, that is, the automorphism group of
$P(\q)$, and $\mathcal{G}^p_{k+1}(\q):=L^p_{k+1}(\Aut\;P(\q))$ the
space of $L^p_{k+1}$-gauge transformations. For $A_0\in
C^{\infty}(T^{\ast}\partial B^4\otimes\mathfrak{su}(2))$ and
$\q\in\PP(d_0,\lambda_0)$, we define the following spaces of
connections:
$$\mathcal{A}(A_0;\q)=\{A:\text{$A$ is a smooth connection on $P(\q)$ such that $\iota^{\ast}A\sim A_0$ on $\partial B^4$}\}\;,$$
where $\iota^{\ast}A\sim A_0$ on $\partial B^4$ means that
$\iota^{\ast}A$ is gauge equivalent to $A_0$ over $\partial B^4$,
via a gauge transformation which extends smoothly to $B^4$;
$$\mathcal{A}^p_k(A_0;\q)=\{A:\text{$A$ is a connection of class $L^p_k$ on $P(\q)$ such that $\iota^{\ast}A\sim A_0$ on $\partial B^4$}\}\;,$$
where this time $\iota^{\ast}A\sim A$ on $\partial B^4$ means that
$\iota^{\ast}A$ is gauge equivalent to $A_0$ over $\partial B^4$,
via an $L^p_{k+1-1/p}$-gauge transformation on $P(\q)|_{\partial
B^4}$ which admits an $L^p_{k+1}$ extension to $B^4$. We henceforth
assume that $(k+1)p\ge 4$.

\noindent The spaces $\mathcal{A}(A_0;\q)$ and
$\mathcal{A}^p_k(A_0;\q)$ have connected components labeled by the
integers $\Z$ (cf.~\cite{IM}):
\begin{align*}
\mathcal{A}(A_0;\q)=\bigsqcup_{j\in\Z}\mathcal{A}_j(A_0;\q),\;\text{
with }
\mathcal{A}_j(A_0;\q)=\{A\in\mathcal{A}(A_0;\q):\mathsf{c}_2(A)=j\}\;,\\
\mathcal{A}^p_k(A_0;\q)=\bigsqcup_{j\in\Z}\mathcal{A}^p_{k,j}(A_0;\q),\;
\text{ with }
\mathcal{A}^p_{k,j}(A_0;\q)=\{A\in\mathcal{A}^p_{k,j}(A_0;\q):\mathsf{c}_2(A)=j\}\;,
\end{align*}
where
$\mathsf{c}_2(A):=\frac{\epsilon^2}{8\pi^2}\int_{B^4}\tr({F_A}^{\epsilon}\wedge
{F_A}^{\epsilon})-\frac{\epsilon^2}{8\pi^2}\int_{B^4}\tr({F_{\underline{A}_{\epsilon}}}^{\epsilon}\wedge
{F_{\underline{A}_{\epsilon}}}^{\epsilon})$ is the relative 2nd
Chern class of $A$ with respect to $\underline{A}_{\epsilon}$.

\noindent Since the groups $\mathcal{G}(\q)$,
$\mathcal{G}^p_{k+1}(\q)$, respectively,  act on
$\mathcal{A}(A_0;\q)$, $\mathcal{A}^p_k(A_0;\q)$, preserving these
connected components, we consider the quotient spaces
\begin{align*}
\mathcal{B}(A_0;\q):=\mathcal{A}(A_0;\q)/\mathcal{G}(\q);\quad
\mathcal{B}_{j}(A_0;\q):=\mathcal{A}_j(A_0;\q)/\mathcal{G}(\q)\;,\\
\mathcal{B}^p_k(A_0;\q):=\mathcal{A}^p_k(A_0;\q)/\mathcal{G}^p_{k+1}(\q);\quad
\mathcal{B}^p_{k,j}(A_0;\q)=\mathcal{A}^p_{k,j}(A_0;\q)/\mathcal{G}^p_{k+1}(\q)\;,
\end{align*}
and denote by $[A]$ the class of $A$ in $\mathcal{B}(A_0;\q)$, or in
$\mathcal{B}^p_k(A_0;\q)$.

\noindent We will also make use of the subgroups
$\mathcal{G}^{\ast}(\q)$ of $\mathcal{G}(\q)$, and
$\mathcal{G}^{\ast,p}_{k+1}(\q)$ of $\mathcal{G}^p_{k+1}(\q)$,
consisting of all gauge transformations $g$ such that
$g((1,0,0,0))=\mathbf{1}$, and of corresponding subspaces
$\mathcal{A}^{\ast}(A_0;\q)$ of $\mathcal{A}(A_0;\q)$, and
$\mathcal{A}^{\ast,p}_k(A_0;\q)$ of $\mathcal{A}^p_k(A_0;\q)$,
consisting of all connections satisfying $\iota^{\ast}A\sim A_0$ on
$\partial B^4$ via a gauge transformation $g$ with
$g((1,0,0,0))=\mathbf{1}$ extendible to $B^4$. The subspaces
$\mathcal{A}^{\ast}_j(A_0;\q)$ and
$\mathcal{A}^{\ast,p}_{k,j}(A_0;\q)$ are defined analogously. The
groups $\mathcal{G}^{\ast}(\q)$, $\mathcal{G}^{\ast,p}_{k+1}(\q)$,
act freely on $\mathcal{A}^{\ast}(A_0;\q)$,
$\mathcal{A}^{\ast,p}_k(A_0;\q)$, respectively, and the
corresponding quotients $\mathcal{B}^{\ast}(A_0;\q)$,
$\mathcal{B}^{\ast,p}_k(A_0;\q)$,  are proved to be differentiable
manifolds provided that $(k+1)p>4$ (cf. $\S1.1$ in \cite{IM3}).


\noindent The $\mathcal{YM}_{\epsilon}$-action (cf.~\eqref{1}),
 descends to the quotients $\mathcal{B}(A_0;\q)$,
$\mathcal{B}^{\ast}(A_0;\q)$, $\mathcal{B}^p_{k}(A_0;\q)$,
$\mathcal{B}^{\ast,p}_k(A_0;\q)$, thus solutions to
$\bigl(\mathcal{D}_{\epsilon}\bigr)$ are critical points of
$\mathcal{YM}_{\epsilon}$ on $\mathcal{B}(A_0;\q)$,
$\mathcal{B}^{\ast}(A_0;\q)$, $\mathcal{B}^p_{k}(A_0;\q)$,
$\mathcal{B}^{\ast,p}_k(A_0;\q)$.

As remarked previously, the bundles $P(\q)$ are all isomorphic to
$P({\q}_0)$ for any fixed ${\q}_0\in\PP(d_0,\lambda_0)$, thus, from
now on, we may omit the indication of ${\q}_0$ (or $\q$) from our
notation and simply write $P$, $\mathcal{A}(A_0)$,
$\mathcal{A}^{\ast}(A_0)$, $\mathcal{A}^p_k(A_0)$,
$\mathcal{A}^{\ast,p}_k(A_0)$, $\mathcal{G}$, $\mathcal{G}^{\ast}$,
$\mathcal{G}^p_{k}$, $\mathcal{G}^{\ast,p}_{k+1}$,
$\mathcal{B}^{\ast}(A_0)$ and $\mathcal{B}^{\ast,p}_k(A_0)$ for the
corresponding objects.

\subsection{Asymptotic expansion of $\YMe(A(\q))$}
In this and in the next sections, we show that the connections
$A(\q)$ introduced in $\S3.1$ are indeed good approximate solutions
to the Dirichlet problem $\bigl(\mathcal{D}_{\epsilon}\bigr)$. For
this, we need to consider the following parameter space, a subset of
$ \PP(d_0,\lambda_0)$ (cf. \eqref{par1}, \eqref{par2}). For
$0<D_1<D_2<\infty$, we define
\begin{equation}
\label{3.2}
 \PP(d_0,\lambda_0;D_1,D_2;\epsilon):=\{\q:=(p,[g],\lambda)\in \PP(d_0,\lambda_0):D_1\epsilon<\lambda^2<D_2\epsilon\}.
\end{equation}
Here and throughout the rest of the paper we assume that
$0<2\lambda_0<d_0<1$, $0<D_1<D_2<\infty$, $\epsilon$ and
$\q\in\PP(d_0,\lambda_0;D_1,D_2;\epsilon)$ are fixed.

\medskip
\noindent In the following, we choose $\bigl(\xi_1, \xi_2,
\xi_3\bigr)$, where  $\xi_1=\begin{pmatrix}
0&0&0\\ 0&0&-1\\ 0&1&0\end{pmatrix}$, $\xi_2=\begin{pmatrix} 0&0&1\\
0&0&0\\ -1&0&0\end{pmatrix}$, $\xi_3=\begin{pmatrix} 0&-1&0\\
1&0&0\\ 0&0&0\end{pmatrix}$, as basis for the Lie algebra
 $\mathfrak{so}(3)$. Right translation by $g$ yields a basis for
$T_{[g]}SO(3)$, which we denote by $\bigl(\xi_1[g], \xi_2[g],
\xi_3[g]\bigr)$.

\medskip
\noindent The following proposition holds for the functional
\begin{equation}
\label{3.3} J_{\epsilon}(\q):= \mathcal{YM}_1(\epsilon
A(\q))=\epsilon^2\mathcal{YM}_{\epsilon}(A(\q))\,,
\end{equation}
where the second equality above comes from the Lie algebra
isomorphism $\mathfrak{su}(2)_\epsilon\simeq\mathfrak{su}(2).$
\begin{proposition}
\label{P3.1} In the hypotheses above, for $\q\in
\PP(d_0,\lambda_0;D_1,D_2;\epsilon)$, the functional
$J_{\epsilon}(\q)$ has the following asymptotic expansion:
\begin{align}
\label{J-epsilon}
J_{\epsilon}(\q)=~8\pi^2+\epsilon^2\int_{B^4}|{F_{\underline{A}_{\epsilon}}}^{\epsilon}|^2\,dx
+2\lambda^4\int_{B^4}|(dh_p)^-|^2\,dx
-4\epsilon\lambda^2\int_{B^4}(({F_{\underline{A}_0}}^0)^-,g(dh_p)^-g^{-1})\,dx
+r_1(\q)\,,
\end{align}
where $r_1(\q)=O(\epsilon^3)$ as $\epsilon\to 0$.

\noindent Moreover, for $\mathcal{F}_{\epsilon}$  defined by
$$
\mathcal{F}_{\epsilon}(\q):=2\lambda^4\int_{B^4}|(dh_p)^-|^2\,dx-4\epsilon\lambda^2\int_{B^4}
\bigl(({F_{\underline{A}_0}}^0)^-,g(dh_p)^-g^{-1}\bigr)\,dx,
$$
one has
$$J_{\epsilon}'(\q)=\mathcal{F}_{\epsilon}'(\q)+r_2(\q),$$
where $J_{\epsilon}'(\q)$ and $\mathcal{F}_{\epsilon}'(\q)$ are the
derivatives of $J_{\epsilon}$ and $\mathcal{F}_{\epsilon}$ with
respect to the variable $\q$, and
\begin{align}
&r_2(\q)\Big(\frac{\partial}{\partial
p^i}\Big)=O(\epsilon^{5/2})~(1\le i\le 4),\quad r_2(\q)(\xi_i[g])=O(\epsilon^3|\log\epsilon|)~(1\le i\le 3)\,,\label{3.5}\\
&r_2(\q)\Big(\frac{\partial}{\partial\lambda}\Big)=O(\epsilon^{5/2}|\log\epsilon|)\,,\label{3.6}
\end{align}
for $\q\in \PP(d_0,\lambda_0;D_1,D_2;\epsilon)$ and small
$\epsilon>0$.
\end{proposition}
In order to prove the proposition above, we need an estimate for the
one-form $h_{\lambda,p}$ defined in \S3.1. We recall that
$h_{\lambda,p}$ and $h_p$ are harmonic componentwise on $B^4$,
$h_{\lambda,p}=\text{Im}\frac{\lambda^2\,\overline{x-p}\,dx}{|x-p|^2(\lambda^2+|x-p|^2)}$
at $\partial B^4$,  and
$h_p=\text{Im}\frac{\,\overline{x-p}\,dx}{|x-p|^4}$ at $\partial
B^4$. We have the following:

\begin{lemma} \label{L5.6}  Let $d_0,\lambda\in(0,1)$. For any $k\ge 1$, there exists a constant $C_{k,d_0}$ depending only on $k$ and $d_0$ such that
$$\|h_{\lambda,p}-\lambda^2h_p\|_{C^k(B^4)}\le C_{k,d_0}\lambda^4$$
for any $p\in B_{1-d_0}$.
\end{lemma}
\textit{Proof:} For $x\in\partial B^4$, we have
\begin{equation}
\label{5.9}
h_{\lambda,p}=\text{Im}\biggl(\frac{\lambda^2\,\overline{x-p}\;dx}{|x-p|^4}
-\frac{\lambda^4\,\overline{x-p}\;dx}{|x-p|^4(\lambda^2+|x-p|^2)}\biggr).
\end{equation}
Let us define
$$r_{\lambda,p}(x)=-\frac{\lambda^4\,\overline{x-p}\;dx}{|x-p|^4(\lambda^2+|x-p|^2)}\,,\; x\neq p\,.$$
It is easy to see that there exists a constant $C$ depending only on
$k$ and $d_0$ such that $|\nabla^kr_{\lambda,p}(x)|\le C\lambda^4$
for any $p\in B_{1-d_0}$, $x\in B^4\setminus B_{1-d_0/2}$, and
$\lambda\in(0,1)$.

Let $\varphi$ be a smooth cut-off function such that $\varphi(x)=0$
in $B_{1-d_0/2}$, $\varphi(x)=1$ on $\R^4\setminus B_{1-d_0/4}$ and
$|\nabla\varphi(x)|\le 8d_0^{-1}$. Then there exists another
constant $C$ depending only on $k$ and $d_0$ such that
\begin{equation}
\label{phi r} |\nabla^k(\varphi(x)r_{\lambda,p}(x))|\le C\lambda^4
\end{equation}
for $p\in B_{1-d_0}$, $x\in B^4$ and $\lambda\in (0,1)$. Since
$h_{\lambda,p}-\lambda^2h_p$ is a harmonic function on $B^4$ with
the same boundary value as $\varphi(x)r_{\lambda,p}(x)$ at $\partial
B^4$, the assertion of the lemma follows from \eqref{phi r} and
elliptic estimates (c.f.~\cite{GT}). \hfill$\Box$

\textit{Proof of Proposition \ref{P3.1}:} Throughout these
estimates, the expression $Q(\q)\lesssim f(\epsilon)$ for a quantity
$Q(\q)$ depending on $\q\in\PP(d_0,\lambda_0;D_1,D_2;\epsilon)$ and
a function $f(\epsilon)$ means that there exists a constant $C=
C(d_0; \lambda_0; D_1; D_2)$ such that $Q(\q)\le C f(\epsilon)$ for
all $\q\in\PP(d_0,\lambda_0;D_1,D_2;\epsilon)$ and small positive
$\epsilon$. We write $Q(\q)\simeq f(\epsilon)$ if the inequality
holds both ways.

We first prove the asymptotic expansion \eqref{J-epsilon}. Due to
the definition of $A(\q)$, we decompose the domain $B^4$ into four
subdomains: $B^4\setminus B_{2\lambda}(p)$,
$B_{2\lambda}(p)\setminus B_{\lambda/2}(p)$,
$B_{\lambda/2}(p)\setminus B_{\lambda/4}(p)$ and $B_{\lambda/4}(p)$.

\medskip

\noindent\textbf{(1) Estimate of $\int_{B^4\setminus
B_{2\lambda}(p)}|{F_{A(\q)}}^{\epsilon}|^2\,dx$}

\noindent On $B^4\setminus B_{2\lambda}(p)$,
$A(\q)=\underline{A}_{\epsilon}+\frac{1}{\epsilon}g
PI^2_{\lambda,p}g^{-1}$ and ${F_{A(\q)}}^{\epsilon}={F_{\underline
A_\epsilon}}^{\epsilon}+{F_{\frac{1}{\epsilon}gdPI^2_{\lambda,p}g^{-1}}}^{\epsilon}+[\underline
A_\epsilon , gPI^2_{\lambda,p}g^{-1}].$

\noindent We recall that ${F_{\underline
A_\epsilon}}^{\epsilon}=d\underline
A_\epsilon+\frac{\epsilon}{2}[\underline A_\epsilon,\underline
A_\epsilon],$
 ${F_{\frac{1}{\epsilon}gPI^2_{\lambda,p}g^{-1}}}^{\epsilon}=
d (\frac{1}{\epsilon}gPI^2_{\lambda,p}g^{-1}) + \frac{\epsilon}{2}
[\frac{1}{\epsilon}gPI^2_{\lambda,p}g^{-1},
\frac{1}{\epsilon}gPI^2_{\lambda,p}g^{-1}] = \frac{1}{\epsilon}
F_{gPI^2_{\lambda,p}g^{-1}}$ with $F_{gPI^2_{\lambda,p}g^{-1}} := d
(gPI^2_{\lambda,p}g^{-1}) +\frac{1}{2} [gPI^2_{\lambda,p}g^{-1}, g
PI^2_{\lambda,p}g^{-1}]$.

\noindent Thus,
\begin{align}
\label{YM1}   &\epsilon^2\int_{B\setminus
B_{2\lambda}(p)}|{F_{A(\q)}}^{\epsilon}|^2\,dx=\epsilon^2
\int_{B\setminus B_{2\lambda}(p)}|{F_{\underline
A_\epsilon}}^{\epsilon}|^2 \,dx+\int_{B\setminus
B_{2\lambda}(p)}|F_{\,PI^2_{\lambda,p}}|^2\,dx
+\epsilon^2\int_{B\setminus
B_{2\lambda}(p)}|[\underline A_\epsilon , gPI^2_{\lambda,p}g^{-1}]|^2\,dx\notag\\
  &+2\epsilon\int_{B\setminus
B_{2\lambda}(p)}({F_{\underline
A_\epsilon}}^{\epsilon},F_{gPI^2_{\lambda,p}g^{-1}})\,dx+
 2\epsilon^2 \int_{B\setminus
B_{2\lambda}(p)}({F_{\underline A_\epsilon}}^{\epsilon},[\underline
A_\epsilon,
gPI^2_{\lambda,p}g^{-1}])\,dx\notag\\
  &+2\epsilon\int_{B\setminus
B_{2\lambda}(p)}(F_{g\,PI^2_{\lambda,p}\,g^{-1}},[\underline
A_\epsilon,
 gPI^2_{\lambda,p}g^{-1}])\,dx:= E_1 + \cdots +
 E_6.
\end{align}
We now estimate all terms $E_2$--$E_6$ in \eqref{YM1}.

\medskip

\noindent $\bullet$ \emph{Estimate of $E_2$:}

\noindent On any domain $D\subset B\setminus \{p\}$ one has
\begin{align}
 \int_{D}|F_{\,PI^2_{\lambda,p}}|^2\,dx \,=
\int_{D}|F_{I^2_{\lambda,p}}|^2\,dx\,+
\int_{D}|d_{I^2_{\lambda,p}}h_{\lambda,p}|^2\,dx\,
-\,2\int_{D}(F_{I^2_{\lambda,p}},d_{I^2_{\lambda,p}}h_{\lambda,p})\,dx\notag\\
 +\int_{D}(F_{I^2_{\lambda,p}},[h_{\lambda,p},h_{\lambda,p}])\,dx+
\frac{1}{4}\int_{D}|[h_{\lambda,p},h_{\lambda,p}]|^2\,dx
-\int_{D}(d_{I^2_{\lambda,p}}h_{\lambda,p},[h_{\lambda,p},h_{\lambda,p}])\,dx\notag
\end{align}

\noindent recalling that $PI^2_{\lambda,p}=
I^2_{\lambda,p}-h_{\lambda,p},$
$F_{I^2_{\lambda,p}}=dI^2_{\lambda,p}+\frac{1}{2}[I^2_{\lambda,p},I^2_{\lambda,p}]$
and
$d_{I^2_{\lambda,p}}h_{\lambda,p}=dh_{\lambda,p}+[I^2_{\lambda,p},h_{\lambda,p}]$.

\noindent By Lemma \ref{L5.6} ($\Vert h_{\lambda,p}\Vert_\infty
\lesssim \lambda^2\Vert h_{p}\Vert_\infty$, and $\Vert
dh_{\lambda,p}\Vert_\infty \lesssim\lambda^2\Vert \nabla
h_{p}\Vert_\infty$),
\begin{align}
  \frac{1}{4}\int_{D}|[h_{\lambda,p},h_{\lambda,p}]|^2\,dx&\lesssim\int_{D}|h_{\lambda,p}|^4\,dx
  \lesssim\lambda^8\lesssim\epsilon^4\,,\notag\\
 \biggl|\int_{D}(d_{I^2_{\lambda,p}}h_{\lambda,p},[h_{\lambda,p},h_{\lambda,p}])\,dx\biggr|
&\lesssim\int_{D}|dh_{\lambda,p}||h_{\lambda,p}|^2\,dx+\int_{D}|I^2_{\lambda,p}||h_{\lambda,p}|^3\,dx
\lesssim\lambda^6\lesssim\epsilon^3.\notag
\end{align}

\noindent Thus, for any $D\subset B\setminus \{p\}$, one has
\begin{align}
\label{*}  \int_{D}|F_{\,PI^2_{\lambda,p}}|^2\,dx&=
\int_{D}|F_{I^2_{\lambda,p}}|^2\,dx+
\int_{D}|d_{I^2_{\lambda,p}}h_{\lambda,p}|^2\,dx\notag\\
&\quad-2\int_{D}(F_{I^2_{\lambda,p}},d_{I^2_{\lambda,p}}h_{\lambda,p})\,dx+
\int_{D}(F_{I^2_{\lambda,p}},[h_{\lambda,p},h_{\lambda,p}])\,dx+
O(\epsilon^3)
\end{align}

\noindent $\bullet$ \emph{Estimates of $E_3-E_5$:}

By Proposition \ref{P2.1} ($\Vert \underline A_\epsilon\Vert
\lesssim \Vert \underline A_0\Vert_\infty$, $\Vert F_{\underline
A_\epsilon}\Vert \lesssim \Vert F_{\underline A_0}\Vert_\infty$ ),
and by Lemma \ref{L5.6},
\begin{align}
  &E_3 :=\epsilon^2 \int_{B\setminus B_{2\lambda}(p)}\vert
[\underline A_\epsilon,g(I^2_{\lambda,p}-
h_{\lambda,p})g^{-1}]\vert^2\,dx
\,\lesssim\,\epsilon^2\int_{B\setminus
B_{2\lambda}(p)}|I^2_{\lambda,p}|^2\,dx\,\lesssim\,\epsilon^2\lambda^2\lesssim\epsilon^3\,;\notag\\
  &E_4= 2\epsilon\int_{B\setminus B_{2\lambda}(p)}({F_{\underline
A_\epsilon}}^{\epsilon},gF_{I^2_{\lambda,p}}g^{-1})\,dx-2\epsilon\int_{B\setminus
B_{2\lambda}(p)}({F_{\underline
A_\epsilon}}^{\epsilon},gd_{I^2_{\lambda,p}}h_{\lambda,p}g^{-1})\,dx+
O(\epsilon^3);\notag\\
  &E_5:=
 2\epsilon^2\int_{B\setminus B_{2\lambda}(p)}({F_{\underline
A_\epsilon }}^{\epsilon},[\underline A_\epsilon , g(I^2_{\lambda,p}-
h_{\lambda,p})g^{-1}])\,dx\lesssim\epsilon^2\int_{B\setminus
B_{2\lambda}(p)}|I^2_{\lambda,p}|\,dx\,\lesssim\epsilon^2\lambda^2\lesssim
\epsilon^3\,.\notag\\
\notag
\end{align}

\noindent $\bullet$ \emph{Estimate of $E_6$:}

\begin{align}
  E_6&=2\epsilon\int_{B\setminus B_{2\lambda}(p)}\bigl(g(F_{I^2_{\lambda,p}}-d_{I^2_{\lambda,p}}h_{\lambda,p}+
\frac{1}{2}[h_{\lambda,p},h_{\lambda,p}])g^{-1},[\underline A_\epsilon ,g(I^2_{\lambda,p}-h_{\lambda,p})g^{-1}]\bigr)\,dx\notag\\
&=\,2\epsilon\int_{B\setminus
B_{2\lambda}(p)}(gF_{I^2_{\lambda,p}}g^{-1},[\underline
A_\epsilon,gPI^2_{\lambda,p}g^{-1}])\,dx+O(\epsilon^3),\notag
\end{align}
\noindent since
\begin{align}
 &\epsilon\int_{B\setminus B_{2\lambda}(p)}|\underline A_\epsilon
||d_{I^2_{\lambda,p}}h_{\lambda,p}||I^2_{\lambda,p}| \,dx\lesssim
\int_{B\setminus B_{2\lambda}(p)}|\bigl(|d
h_{\lambda,p}||I^2_{\lambda,p}| + |I^2_{\lambda,p}|^2\,
|h_{\lambda,p}|\bigr)\,dx \lesssim \epsilon \lambda^4 \sim
\epsilon^3\,\notag\\
&\int_{B\setminus B_{2\lambda}(p)}|\underline A_\epsilon
|\,|d_{I^2_{\lambda,p}}h_{\lambda,p}||h_{\lambda,p}|\,dx\lesssim
\int_{B\setminus
B_{2\lambda}(p)}\bigl(|dh_{\lambda,p}||h_{\lambda,p}|+|I^2_{\lambda,p}|\,|h_{\lambda,p}|^2\bigr)\,dx
\lesssim \epsilon\lambda^4\sim \epsilon^3\,,\notag\\
&\int_{B\setminus B_{2\lambda}(p)}\bigl(|\underline A_\epsilon |\,
|I^2_{\lambda,p}||h_{\lambda,p}|^2 +
|h_{\lambda,p}|^3\bigr)\,dx\,\lesssim\,\epsilon\lambda^6\,\sim\,
\epsilon^4. \notag
\end{align}

\noindent And, adding up $E_1+\cdots + E_6,$
\begin{align}
\label{YM1-I}   &\epsilon^2 \int_{B\setminus
B_{2\lambda}(p)}|{F_{A(\q)}}^{\epsilon}|^2\,dx\,=\,
\epsilon^2\int_{B\setminus
B_{2\lambda}(p)}|{F_{\underline{A}_{\epsilon}}}^{\epsilon}|^2\,dx\,+\,
\int_{B^4\setminus B_{2\lambda}(p)}|F_{I^2_{\lambda,p}}|^2\,dx\notag\\
   &+\int_{B\setminus
B_{2\lambda}(p)}|d_{I^2_{\lambda,p}}h_{\lambda,p}|^2\,dx
-2\int_{B\setminus
B_{2\lambda}(p)}(F_{I^2_{\lambda,p}},d_{I^2_{\lambda,p}}h_{\lambda,p})\,dx+\int_{B\setminus
B_{2\lambda}(p)}(F_{I^2_{\lambda,p}},[h_{\lambda,p},h_{\lambda,p}])\,dx\notag\\
   &+2\epsilon\int_{B\setminus B_{2\lambda}(p)}({F_{\underline
A_\epsilon}}^{\epsilon},gF_{I^2_{\lambda,p}}g^{-1})\,dx-2\epsilon\int_{B\setminus
B_{2\lambda}(p)}({F_{\underline A_\epsilon}}^{\epsilon},
gd_{I^2_{\lambda,p}}h_{\lambda,p}g^{-1})\,dx\notag\\
   &+2\epsilon\int_{B\setminus
B_{2\lambda}(p)}(gF_{I^2_{\lambda,p}}g^{-1},[\underline A_\epsilon
,gPI^2_{\lambda,p}g^{-1}]))\,dx+O(\epsilon^3).
\end{align}

\noindent Quite analogously, we obtain the following estimate:
\begin{align}
\label{Tr1} &\epsilon^2\int _{B\setminus
B_{2\lambda}(p)}\tr\bigl({F_{A(\q)}}^{\epsilon}\wedge
{F_{A(\q)}}^{\epsilon})= \epsilon^2\int _{B\setminus
B_{2\lambda}(p)}\tr\bigl({F_{\underline
A_\epsilon}}^{\epsilon}\wedge {F_{\underline
A_\epsilon}}^{\epsilon}\bigr)+\int _{B\setminus
B_{2\lambda}(p)}\tr\bigl(F_{I^2_{\lambda,p}}\wedge
F_{I^2_{\lambda,p}}\bigr)+\notag\\
&+\int_{B\setminus
B_{2\lambda}(p)}\tr\bigl(d_{I^2_{\lambda,p}}h_{\lambda,p}\wedge
d_{I^2_{\lambda,p}}h_{\lambda,p}\bigr)-2 \int_{B\setminus
B_{2\lambda}(p)}\tr\bigl(F_{I^2_{\lambda,p}}\wedge
d_{I^2_{\lambda,p}}h_{\lambda,p}\bigr)+\int_{B\setminus
B_{2\lambda}}\tr\bigl(F_{I^2_{\lambda,p}}\wedge [h_{\lambda,p},h_{\lambda,p}]\bigr)\notag\\
&+2\epsilon\int_{B\setminus B_{2\lambda}(p)}\tr\bigl({F_{\underline
A_\epsilon}}^{\epsilon}\wedge
gF_{I^2_{\lambda,p}}g^{-1}\bigr)-2\epsilon\int_{B\setminus
B_{2\lambda}(p)}\tr\bigl({F_{\underline
A_\epsilon}}^{\epsilon}\wedge
gd_{I^2_{\lambda,p}}h_{\lambda,p}g^{-1}\bigr)\notag\\
& +2\epsilon\int_{B\setminus
B_{2\lambda}(p)}\tr\bigl(gF_{I^2_{\lambda,p}}g^{-1}\wedge[\underline
A_\epsilon,gPI^2_{\lambda,p}g^{-1}]\bigr)+O(\epsilon^3).
\end{align}

\noindent Summing \eqref{Tr1} to  \eqref{YM1-I}, and using
$*F_{I^2_{\lambda,p}}= F_{I^2_{\lambda,p}}$, finally yields

\begin{align}
\label{D1} &\epsilon^2\int_{B\setminus
B_{2\lambda}(p)}|{F_{A(\q)}}^{\epsilon}|^2\,dx+\epsilon^2\int
_{B\setminus B_{2\lambda}(p)}\tr\bigl({F_{A(\q)}}^{\epsilon}\wedge
{F_{A(\q)}}^{\epsilon})= \epsilon^2\int_{B\setminus
B_{2\lambda}(p)}|{F_{\underline{A}_{\epsilon}}}^{\epsilon}|^2\,dx\notag\\
&+\epsilon^2\int _{B\setminus
B_{2\lambda}(p)}\tr\bigl({F_{\underline
A_\epsilon}}^{\epsilon}\wedge {F_{\underline
A_\epsilon}}^{\epsilon}\bigr)+2\int_{B\setminus
B_{2\lambda}(p)}|d_{I^2_{\lambda,p}}h_{\lambda,p}^-|^2\,dx
-4\epsilon\int_{B\setminus B_{2\lambda}(p)}(({F_{\underline
A_\epsilon}}^{\epsilon})^-,gd_{I^2_{\lambda,p}}h_{\lambda,p}g^{-1})\,dx
+O(\epsilon^3).
\end{align}

\medskip

\noindent \textbf{(2) Estimate of
$\epsilon^2\int_{B_{2\lambda}(p)\setminus
B_{\lambda/2}(p)}|{F_{A(\q)}}^{\epsilon}|^2\,dx$}\;.

\noindent On $B_{2\lambda}(p)\setminus B_{\lambda/2}(p)$,
$A(\q)=(1-\beta_{\lambda,p})\underline{A}_{\epsilon}+\frac{1}{\epsilon}gPI^2_{\lambda,p}g^{-1}$
and
\begin{align}
\label{YM2} & \epsilon^2\int_{B_{2\lambda}(p)\setminus
B_{\lambda/2}(p)}|{F_{A(\q)}}^{\epsilon}|^2\,dx =\epsilon^2
\int_{B_{2\lambda}(p)\setminus B_{\lambda/2}(p)}
|{F_{(1-\beta_{\lambda,p})\underline{A}_{\epsilon}}}^{\epsilon}|^2 \,dx\notag\\
&+ \int_{B_{2\lambda}(p)\setminus
B_{\lambda/2}(p)}|F_{\,PI^2_{\lambda,p}}|^2\,dx
+\epsilon^2\int_{B_{2\lambda}(p)\setminus B_{\lambda/2}(p)}|[(1-\beta_{\lambda,p})\underline{A}_{\epsilon} , gPI^2_{\lambda,p}g^{-1}]|^2\,dx\notag\\
&+2\epsilon\int_{B_{2\lambda}(p)\setminus
B_{\lambda/2}(p)}({F_{(1-\beta_{\lambda,p})\underline{A}_{\epsilon}}}^{\epsilon},F_{gPI^2_{\lambda,p}g^{-1}})\,dx+
 2\epsilon^2 \int_{B_{2\lambda}(p)\setminus B_{\lambda/2}(p)}({F_{(1-\beta_{\lambda,p})\underline{A}_{\epsilon}}}^{\epsilon},[(1-\beta_{\lambda,p})\underline{A}_{\epsilon},
gPI^2_{\lambda,p}g^{-1}])\,dx\notag\\
&+2\epsilon\int_{B_{2\lambda}(p)\setminus
B_{\lambda/2}(p)}(F_{gPI^2_{\lambda,p}g^{-1}},[(1-\beta_{\lambda,p})\underline{A}_{\epsilon},
 gPI^2_{\lambda,p}g^{-1}])\,dx:= E_1 + \cdots +
 E_6\notag\\
\end{align}

\medskip

\noindent $\bullet$ \emph{Estimates of $E_1$--$E_6$:}

\noindent Since
${F_{(1-\beta_{\lambda,p})\underline{A}_{\epsilon}}}^{\epsilon}=(1-\beta_{\lambda,p})
{F_{\underline{A}_{\epsilon}}}^{\epsilon}-d\beta_{\lambda,p}\wedge\underline{A}_{\epsilon}-
\frac{\epsilon}{2}\beta_{\lambda,p}(1-\beta_{\lambda,p})[\underline{A}_{\epsilon},\underline{A}_{\epsilon}]\,,$
we find
\begin{equation}
\notag E_1 \,:=\,\epsilon^2\int_{B_{2\lambda}(p)\setminus
B_{\lambda/2}(p)}|{F_{(1-\beta_{\lambda,p})\underline{A}_{\epsilon}}}^{\epsilon}|^2\,dx\lesssim\epsilon^3\,.
\end{equation}

\noindent The estimate of $E_2$ is already done (cf. \eqref{*} with
$D= B_{2\lambda}(p)\setminus B_{\lambda/2}(p)$).

\noindent By Lemma \ref{L5.6} and Proposition \ref{P2.1},
\begin{align}
\notag & |E_3|\lesssim \epsilon^2 \Vert \underline
A_\epsilon\Vert_\infty \int_{B_{2\lambda}(p)\setminus
B_{\lambda/2}(p)} \bigl(|I^2_{\lambda,p}|^2+
|h_{\lambda,p}|^2\bigr)\,dx\lesssim
\epsilon\lambda^2\lesssim \epsilon^3;\notag\\
&E_4= \epsilon\int_{B_{2\lambda}(p)\setminus
B_{\lambda/2}(p)}({F_{(1-\beta_{\lambda,p})\underline{A}_{\epsilon}}}^{\epsilon},gF_{I^2_{\lambda,p}}g^{-1})\,dx+O(\epsilon^3);\notag\\
&|E_5|\lesssim\epsilon^2\lambda^{-1}\int_{B_{2\lambda}(p)\setminus
B_{\lambda/2}(p)}
|I^2_{\lambda,p}|\,dx+\epsilon^2\lambda^{-1}\int_{B_{2\lambda}(p)\setminus
B_{\lambda/2}(p)}|h_{\lambda,p}|\,dx
\lesssim\epsilon^2\lambda^2\lesssim\epsilon^3;\notag\\
&E_6=\epsilon\int_{B_{2\lambda}(p)\setminus
B_{\lambda/2}(p)}(gF_{I^2_{\lambda,p}}g^{-1},[(1-\beta_{\lambda,p})\underline
A_{\epsilon},gPI^2_{\lambda,p}g^{-1}])\,dx+O(\epsilon^3),\notag\\
\notag
\end{align}
\noindent since
\begin{align}
&\int_{B_{2\lambda}(p)\setminus B_{\lambda/2}(p)}|\underline
A_\epsilon|\,\bigl(
|d_{I^2_{\lambda,p}}h_{\lambda,p}||I^2_{\lambda,p}| +
|d_{I^2_{\lambda,p}}h_{\lambda,p}||h_{\lambda,p}|+|I^2_{\lambda,p}||h_{\lambda,p}|^2
+ |h_{\lambda,p}|^3\bigr)\,dx\notag\\
\lesssim&~\epsilon(\lambda^4+\lambda^4+\lambda^7+\lambda^{10})\lesssim\epsilon^3.
\notag
\end{align}

\noindent Adding up $E_1+\cdots+E_6,$ we obtain
\begin{align}
\label{YM2-I} &\epsilon^2\int_{B_{2\lambda}(p)\setminus
B_{\lambda/2}(p)}|{F_{A(\q)}}^{\epsilon}|^2\,dx
=\int_{B_{2\lambda}(p)\setminus
B_{\lambda/2}(p)}|F_{I^2_{\lambda,p}}|^2\,dx+
\int_{B_{2\lambda}(p)\setminus B_{\lambda/2}(p)}|d_{I^2_{\lambda,p}}h_{\lambda,p}|^2\,dx\notag\\
&-2\int_{B_{2\lambda}(p)\setminus
B_{\lambda/2}(p)}(F_{I^2_{\lambda,p}},d_{I^2_{\lambda,p}}h_{\lambda,p})\,dx+
2\epsilon\int_{B_{2\lambda}(p)\setminus B_{\lambda/2}(p)}({F_{(1-\beta_{\lambda,p})\underline{A}_{\epsilon}}}^{\epsilon},gF_{I^2_{\lambda,p}}g^{-1})\,dx\notag\\
&+2\epsilon\int_{B_{2\lambda}(p)\setminus
B_{\lambda/2}(p)}(gF_{I^2_{\lambda,p}}g^{-1},[(1-\beta_{\lambda,p})\underline{A}_{\epsilon},gPI^2_{\lambda,p}g^{-1}])
\,dx+O(\epsilon^3).
\end{align}

\noindent With the same calculation, one obtains the analogous
formula for trace. Adding up the two and using self duality (as done
earlier to obtain \eqref{D1}), finally yields
\begin{equation}
\label{D2} \epsilon^2\int_{B_{2\lambda}(p)\setminus
B_{\lambda/2}(p)}|{F_{A(\q)}}^{\epsilon}|^2\,dx+\epsilon^2\int
_{B\setminus B_{2\lambda}(p)}\tr\,\bigl({F_{A(\q)}}^{\epsilon}\wedge
{F_{A(\q)}}^{\epsilon})=2 \int_{B_{2\lambda}(p)\setminus
B_{\lambda/2}(p)} |d_{I^2_{\lambda,p}}h_{\lambda,p}^-|^2\,dx +
O(\epsilon^3).
\end{equation}

\medskip

\noindent\textbf{(3) Estimate of $\int_{B_{\lambda/2}(p)\setminus
B_{\lambda/4}(p)}|{F_{A(\q)}}^{\epsilon}|^2\,dx$}\;.

\noindent On $B_{\lambda/2}(p)\setminus B_{\lambda/4}(p)$,
$A(\q)=\frac{1}{\epsilon}gPI^2_{\lambda,p}g^{-1}+\frac{1}{\epsilon}
\beta_{\lambda/4,p}gh_{\lambda,p}g^{-1}$ and
\begin{align}
\label{YM3} & \epsilon^2\int_{B_{\lambda/2}(p)\setminus
B_{\lambda/4}(p)}|{F_{A(\q)}}^{\epsilon}|^2 \,dx
=\int_{B_{\lambda/2}(p)\setminus B_{\lambda/4}(p)}|F_{PI^2_{\lambda,p}}|^2 \,dx\notag\\
&+\int_{B_{\lambda/2}(p)\setminus
B_{\lambda/4}(p)}|d_{PI^2_{\lambda,p}}(\beta_{\lambda/4,p}
h_{\lambda,p}) |^2\,dx
+\frac{1}{4}\int_{B_{\lambda/2}(p)\setminus B_{\lambda/4}(p)}|[\beta_{\lambda/4,p}h_{\lambda,p}, \beta_{\lambda/4,p}h_{\lambda,p}]|^2\,dx\notag\\
&+2\int_{B_{\lambda/2}(p)\setminus B_{\lambda/4}(p)}(F_{
PI^2_{\lambda,p}},d_{ PI^2_{\lambda,p}}\beta_{\lambda/4,p}
h_{\lambda,p})\,dx+\int_{B_{\lambda/2}(p)\setminus
B_{\lambda/4}(p)}(F_{ PI^2_{\lambda,p}},[\beta_{\lambda/4,p}
h_{\lambda,p} ,\beta_{\lambda/4,p}
h_{\lambda,p}])\,dx\notag\\
& + \int_{B_{\lambda/2}(p)\setminus B_{\lambda/4}(p)}\bigl(d_{
PI^2_{\lambda,p}}(\beta_{\lambda/4,p}
h_{\lambda,p}),[\beta_{\lambda/4,p}h_{\lambda,p} ,
 \beta_{\lambda/4,p}h_{\lambda,p} ]\bigr)\,dx:= F_1 + \cdots +
 F_6.\notag\\
\end{align}
\medskip

 \noindent $\bullet$ \emph{Estimates of $F_1$--$F_6$:}

 \medskip

\noindent The term $F_1$ has already been estimated in \eqref{*}

\noindent Moreover,
\begin{align}
\notag &|F_2 |\lesssim \int_{B_{\lambda/2}(p)\setminus
B_{\lambda/4}(p)}\bigl(|dh_{\lambda,p}|^2 + |I^2_{\lambda,p}|^2
+|dh_{\lambda,p}|^2 |h_{\lambda,p}|^4\bigr)\,dx
\lesssim\lambda^8+\lambda^6+\lambda^{8}\lesssim\epsilon^3;\notag\\
&|F_3 |\lesssim \int_{B_{\lambda/2}(p)\setminus B_{\lambda/4}(p)}
|h_{\lambda,p}|^4\,dx\lesssim \lambda^{12}\sim\epsilon^6;\notag\\
&F_4 = \int_{B_{\lambda/2}(p)\setminus B_{\lambda/4}(p)}\bigl(F_{
I^2_{\lambda,p}}, d_{
PI^2_{\lambda,p}} (\beta_{\lambda/4,p} h_{\lambda,p})\bigr)\,dx+ O(\epsilon^3),\notag\\
\notag
\end{align}
\noindent since
\begin{align}
&\notag \biggl|\int_{B_{\lambda/2}(p)\setminus
B_{\lambda/4}(p)}\bigl(- 2d_{I^2_{\lambda,p}}h_{\lambda,p} + \frac
{1}{2} [h_{\lambda,p}, h_{\lambda,p}], d_{PI^2_{\lambda,p}}
(\beta_{\lambda/4,p}h_{\lambda,p})\bigr)\,dx\biggr|\notag\\
\lesssim&~\int_{B_{\lambda/2}(p)\setminus B_{\lambda/4}(p)}\bigl(
|dh_{\lambda,p}|^2 + |dh_{\lambda,p}||I^2_{\lambda,p}|
+|h_{\lambda,p}| + |dh_{\lambda,p}| |h_{\lambda,p}|^2 +
|I^2_{\lambda,p}|^2 |h_{\lambda,p}|^2 +
|I^2_{\lambda,p}||h_{\lambda,p}|^4\bigr)\,dx \lesssim
\epsilon^3;\notag\\
\notag
\end{align}
\begin{equation}
\notag F_5= \int_{B_{\lambda/2}(p)\setminus B_{\lambda/4}(p)}
\bigl(F_{I^2_{\lambda,p}},
[\beta_{\lambda/4,p}h_{\lambda,p},\beta_{\lambda/4,p}h_{\lambda,p}])\,dx
+ O(\epsilon^\frac {9}{2}\bigr),
\end{equation}

\noindent since

\begin{align*}
&\biggl|\int_{B_{\lambda/2}(p)\setminus B_{\lambda/4}(p)} \Bigl(-
d_{I^2_{\lambda,p}}h_{\lambda,p} + \frac {1}{4} [h_{\lambda,p},
h_{\lambda,p}],
[\beta_{\lambda/4,p}h_{\lambda,p},\beta_{\lambda/4,p}h_{\lambda,p}]\Big) \biggr|\\
\lesssim&~\int_{B_{\lambda/2}(p)\setminus B_{\lambda/4}(p)}
\bigl(|d_{I^2_{\lambda,p}}h_{\lambda,p}|\,|h_{\lambda,p}|^2 +
|I^2_{\lambda,p}|\,|h_{\lambda,p}|^3 + |h_{\lambda,p}|^4\bigr)\,dx
\lesssim \epsilon^\frac {9}{2};
\end{align*}

$$ |F_6 |\lesssim\int_{B_{\lambda/2}(p)\setminus
B_{\lambda/4}(p)} \bigl(|dh_{\lambda,p}||h_{\lambda,p}|^2 +
|I^2_{\lambda,p}||h_{\lambda,p}|^3 + |h_{\lambda,p}|^4\bigr)\,dx
\lesssim
\lambda^{10}+\lambda^9+\lambda^{12}\lesssim\epsilon^{\frac{9}{2}}.
$$

\medskip
\noindent Adding up $F_1+\cdots +F_6$ yields

\begin{align}
\label{YM3-I} &\epsilon^2\int_{B_{\lambda/2}(p)\setminus
B_{\lambda/4}(p)}|{F_{A(\q)}}^{\epsilon}|^2\,dx
=\int_{B_{\lambda/2}(p)\setminus
B_{\lambda/4}(p)}|F_{I^2_{\lambda,p}}|^2\,dx+
\int_{B_{\lambda/2}(p)\setminus B_{\lambda/4}(p)}|d_{I^2_{\lambda,p}}h_{\lambda,p}|^2\,dx\notag\\
&+2\int_{B_{\lambda/2}(p)\setminus
B_{\lambda/4}(p)}\bigl(F_{I^2_{\lambda,p}},-d_{I^2_{\lambda,p}}h_{\lambda,p}
+ F_{\beta_{\lambda/4,p}h_{\lambda,p}}+
[I^2_{\lambda,p},\beta_{\lambda/4,p}h_{\lambda,p}]\bigr)\,dx+O(\epsilon^3).
\end{align}

\noindent The analogous formula for trace added to the above gives

\begin{align}
\label{D3} &\epsilon^2\int_{B_{\lambda/2}(p)\setminus
B_{\lambda/4}(p)}|{F_{A(\q)}}^{\epsilon}|^2\,dx +
\,\epsilon^2\int_{B_{\lambda/2}(p)\setminus
B_{\lambda/4}(p)}\tr\,\bigl({F_{A(\q)}}^{\epsilon}\wedge
{F_{A(\q)}}^{\epsilon}\bigr)
\notag\\
=&~2\int_{B_{\lambda/2}(p)\setminus
B_{\lambda/4}(p)}|d_{I^2_{\lambda,p}}h_{\lambda,p}^-|^2\,dx+O(\epsilon^3)\;.
\end{align}

\medskip

\noindent\textbf{(4) Estimate of
$\int_{B_{\lambda/4}(p)}|{F_{A(\q)}}^{\epsilon}|^2\,dx$}\;.

\noindent On $B_{\lambda/4}(p)$, we have
$A(\q)=\frac{1}{\epsilon}gI^1_{\lambda,p}g^{-1}$ and
\begin{equation}
\label{YM4}
\epsilon^2\int_{B_{\lambda/4}(p)}|{F_{A(\q)}}^{\epsilon}|^2\,dx=\int_{B_{\lambda/4}(p)}|F_{I^1_{\lambda,p}}|^2\,dx.
\end{equation}

\noindent Thus, combining the above with the likewise formula for
trace,

\begin{equation}
\label{D4}
\epsilon^2\int_{B_{\lambda/4}(p)}|{F_{A(\q)}}^{\epsilon}|^2\,dx+
\epsilon^2\int_{B_{\lambda/4}(p)}\tr\,\bigl({F_{A(\q)}}^{\epsilon}\wedge
{F_{A(\q)}}^{\epsilon}\bigr)=0.
\end{equation}

\medskip
\noindent We are now ready to sum the contributions from the four
sub-domains, i.e. \eqref{D1}, \eqref{D2}, \eqref{D3}, \eqref{D4},
and obtain
\begin{align}
\label{expansion} \epsilon^2\int_B |{F_{A(\q)}}^{\epsilon}|^2\,dx&=-
\epsilon^2\int _{B}\tr\bigl({F_{A(\q)}}^{\epsilon}\wedge
{F_{A(\q)}}^{\epsilon})+ \epsilon^2\int_{B}\tr\bigl({F_{\underline
A_\epsilon}}^{\epsilon}\wedge
{F_{\underline A_\epsilon}}^{\epsilon}\bigr)+\epsilon^2\int_{B} |{F_{\underline{A}_{\epsilon}}}^{\epsilon}|^2\,dx\notag\\
&\quad+2\int_{B}|d_{I^2_{\lambda,p}}h_{\lambda,p}^-|^2\,dx
-4\epsilon\int_{B}(({F_{\underline
A_\epsilon}}^{\epsilon})^-,gd_{I^2_{\lambda,p}}h_{\lambda,p}g^{-1})\,dx
+O(\epsilon^3)\notag\\
= 8\pi^2+\epsilon^2&\int_{B}
|{F_{\underline{A}_{\epsilon}}}^{\epsilon}|^2\,dx+
+2\int_{B}|d_{I^2_{\lambda,p}}h_{\lambda,p}^-|^2\,dx
-4\epsilon\int_{B}(({F_{\underline
A_\epsilon}}^{\epsilon})^-,gd_{I^2_{\lambda,p}}h_{\lambda,p}g^{-1})\,dx
+O(\epsilon^3)\notag\\
= 8\pi^2+ \epsilon^2&\int_{B}
|{F_{\underline{A}_{\epsilon}}}^{\epsilon}|^2\,dx+2\lambda^4
\int_{B}|d_{I^2_{\lambda,p}}h_{p}^-|^2\,dx
-4\epsilon\int_{B}(({F_{\underline
A_0}}^0)^-,gd_{I^2_{\lambda,p}}h_{\lambda,p}g^{-1})\,dx
+O(\epsilon^3),
\end{align}

\noindent where for the first equality we have used the estimates
\begin{align}
&\biggl|\epsilon^2\int_{B_{2\lambda}}\tr\bigl({F_{\underline
A_\epsilon}}^{\epsilon}\wedge {F_{\underline
A_\epsilon}}^{\epsilon}\bigr)\biggr|\lesssim
\epsilon^2\int_{B_{2\lambda}}
|{F_{\underline{A}_{\epsilon}}}^{\epsilon}|^2\,dx\lesssim
\epsilon^2\lambda^4\lesssim \epsilon^4,\notag\\
&\biggl|\epsilon\int_{B_{2\lambda}(p)}(({F_{\underline
A_{\epsilon}}}^{\epsilon})^-,
gd_{I^2_{\lambda,p}}h_{\lambda,p}g^{-1})\,dx \biggr| \lesssim
\epsilon\int_{B_{2\lambda}(p)}\bigl(|dh_{\lambda,p}|+|I^2_{\lambda,p}||h_{\lambda,p}|\bigr)\,dx
\lesssim \epsilon\lambda^6\lesssim\epsilon^{4},\notag\\
&\int_{B_{\lambda/4}(p)}|d_{I^2_{\lambda,p}}h_{\lambda,p}^-|^2\,dx\lesssim\lambda^6\lesssim\epsilon^3,
\notag
\end{align}

\noindent and for the second equality we have used the  topological
constraint
\begin{equation}
\label{constraint}
\epsilon^2\int_{B^4}\tr({F_{A(\q)}}^{\epsilon}\wedge
{F_{A(\q)}}^{\epsilon}) -\epsilon^2\int_{B^4}\tr({F_{\underline
A_\epsilon}}^{\epsilon}\wedge {F_{\underline
A_\epsilon}}^{\epsilon})=8\pi^2,
\end{equation}

\noindent  and the estimate
\begin{equation}
\notag
\biggl|\epsilon\int_{B_{2\lambda}(p)}(({F_{\underline{A}_{\epsilon}}}^{\epsilon})^-
- ({F_{\underline A_{0}}}^0)^-,
gd_{I^2_{\lambda,p}}h_{\lambda,p}g^{-1})\,dx \biggr| \lesssim
\epsilon^2\lambda^6\lesssim \epsilon^5.
\end{equation}

\noindent Thus, \eqref{J-epsilon} holds.

\medskip

\noindent The asymptotic expansion of the derivative of
$J_{\epsilon}(\q)$ is computed similarly and we omit the
calculation. \hfill$\Box$

\subsection{Estimate of $\|\nabla\YMe(A(\q))\|_{A(\q);1,2,\ast}$}
For connections $A$ on the bundle $P:=P(\q)$ and one-forms $a\in
C^{\infty}(T^{\ast}\overline{B}^4\otimes\Ad(P))$, we define the
$L^p_k$-norm $\|a\|_{A;k,p}$ by
\begin{equation}
\label{norm}
\|a\|_{A;k,p}:=\sum_{j=1}^k\|({\nabla_A}^{\epsilon})^ja\|_p+\|a\|_p\;,
\end{equation}
where ${\nabla_A}^{\epsilon}=\nabla+\epsilon[A,\cdot]$,
$({\nabla_A}^{\epsilon})^j={\nabla_A}^{\epsilon}\cdots{\nabla_A}^{\epsilon}$
($j$-times) and $\|\cdot\|_p$ is the $L^p$-norm on $B^4$. We denote
by $L^p_{0,1}(T^{\ast}B^4\otimes\Ad(P))$ the completion of
$C^{\infty}_0(T^{\ast}B^4\otimes\Ad(P))$ with respect to the norm
above, and define the spaces
\begin{align*}
L^p_{0,k}(T^{\ast}B^4\otimes\Ad(P)):=L^p_k(T^{\ast}B^4\otimes\Ad(P))\cap
L^p_{0,1}(T^{\ast}B^4\otimes\Ad(P))\;;\\
L^p_{0^T,k}(T^{\ast}B^4\otimes\Ad(P))=\{\alpha\in
L^p_k(T^{\ast}B^4\otimes\Ad(P)):\iota^{\ast}\alpha=0\}\;.
\end{align*} Note that these are independent of the choice of
the connection $A$.

Let now the spaces $\mathcal{A}^\ast(A_0)$, $\mathcal{B}^\ast(A_0)$
and their Sobolev correspondents be defined as in $\S3.1.$ The
results in \cite{IM3} allow us to identify the tangent bundle
$T\mathcal{B}^{\ast,p}_{k,+1}(A_0)\to\mathcal{B}^{\ast,p}_{k,+1}(A_0)$
with a sub-bundle of
$\mathcal{A}^p_{k,+1}(A_0)\times_{\mathcal{G}^{\ast,p}_{k+1}}L^p_{0^T,k}(T^{\ast}B^4\otimes\Ad(P))
\to\mathcal{B}^{\ast,p}_{k,+1}(A_0)$, defined as
$$\mathcal{S}^p_{k,+1}(A_0):=S^p_{k,+1}(A_0)/\mathcal{G}^{\ast,p}_{k+1}\to\mathcal{B}^{\ast,p}_{k,+1}(A_0)\;,$$
where
$$S^p_{k,+1}(A_0)=\{(A,\alpha)\in\mathcal{A}^p_{k,+1}(A_0)\times L^p_{0^T,k}(T^{\ast}B^4\otimes\Ad(P)):\;{d_A^{\ast}}^{\epsilon}\alpha=0\}$$
and $\mathcal{G}^{\ast,p}_{k+1}$ acts diagonally on
$S^p_{k,+1}(A_0)$.

The correspondence $\mathcal{A}^p_{k,+1}(A_0)\times
L^p_{0^T,k}(T^{\ast}B^4\otimes\Ad(P))\ni (A,a)\mapsto\|a\|_{A;p,k}$
 is $\mathcal{G}^{\ast,p}_{k+1}$-invariant, therefore it descends to the
quotient.

We are interested in the case $k=1$ and $p>2$, thus we define the
gradient of the functional $\YM_{\epsilon}$ on
$\mathcal{B}^{\ast,p}_{1,+1}(A_0)$ as
\begin{equation}
\label{gradYM}
\nabla\YMe(A)(a)=2\int_{B^4}({F_A}^{\epsilon},{d_A}^{\epsilon}a)\;,
\end{equation}
for $(A,a)\in
T\mathcal{B}^{\ast,p}_{1,+1}(A_0)=\mathcal{S}^p_{1,+1}(A_0)$.

Since \eqref{gradYM} is $\mathcal{G}^{\ast,p}_2$-invariant, that is,
$\nabla\YMe(g\cdot A)(g\cdot a)=\nabla\YMe(A)(a)$ for
$g\in\mathcal{G}^{\ast,p}_2$, the functional $\nabla\YMe$ descends
to $\mathcal{S}^p_{1,+1}(A_0)$.

\noindent Observe that a connection $A$ is a solution to
$\bigl(\mathcal{D}_{\epsilon}\bigr)$ if and only if
$\nabla\YMe(A)=0$ on $\mathcal{S}^p_{1,+1}(A_0)$.

We define the dual $L^2_1$-norm of $\nabla\YMe$ as
$$\|\nabla\YMe(A)\|_{A;1,2,\ast}:=\sup\{\nabla\YMe(A)(a):a\in L^2_{0,1}(T^{\ast}B^4\otimes\Ad(P)),~\|a\|_{A;1,2}\le 1\}\;,$$
or equivalently, by the $\mathcal{G}^{\ast,p}_2$-invariance, as
\begin{align}
\label{normgradYM}
\|\nabla\YMe(A)\|_{A;1,2,\ast}=\sup\{\nabla\YMe(A)(a):\;&a\in L^2_{0,1}(T^{\ast}B^4\otimes\Ad(P))\notag\\
&\text{such that ${d_A^{\ast}}^{\epsilon}a=0$ and $\|a\|_{A;1,2}\le
1$}\}\;.
\end{align}

We are now ready to estimate
$\|\nabla\YMe(A(\q))\|_{A(\q);1,2,\ast}$.
\begin{lemma}
 \label{L3.1} For $\q\in\PP(d_0,\lambda_0;D_1,D_2;\epsilon)$, there exists a constant $C>0$ depending only on $d_0,\lambda_0,D_1$
 and $D_2$ such that the following holds for all small $\epsilon>0$:
$$\|\nabla\YMe(A(\q))\|_{A(\q);1,2,\ast}\le C\epsilon^{1/2}.$$
\end{lemma}
\textit{Proof:} Using \eqref{normgradYM}, we estimate
$\nabla\YMe(A(\q))(\alpha)$ for $\alpha\in
L^2_{0,1}(T^{\ast}B^4\otimes\Ad(P))$ with
${d_{A(\q)}^{\ast}}^{\epsilon}\alpha=0$ and
$\|\alpha\|_{A(\q);1,2}\le 1$, with $A(\q)$ represented in the same
gauge as in \eqref{A(q)}. We have
\begin{align}
\label{nabla}
\frac{1}{2}\nabla\YMe(A(\q))(\alpha)&=\int_{B^4}({F_{A(\q)}}^{\epsilon},{d_{A(\q)}}^{\epsilon}\alpha)\,dx
=\int_{B^4}({d_{A(\q)}^{\ast}}^{\epsilon}{F_{A(\q)}}^{\epsilon},\alpha)\,dx\notag\\
&=\int_{B^4\setminus
B_{2\lambda}(p)}({d_{A(\q)}^{\ast}}^{\epsilon}{F_{A(\q)}}^{\epsilon},\alpha)\,dx+
\int_{B_{2\lambda}(p)\setminus B_{\lambda/2}(p)}({d_{A(\q)}^{\ast}}^{\epsilon}{F_{A(\q)}}^{\epsilon},\alpha)\,dx\notag\\
&\quad+\int_{B_{\lambda/2}(p)\setminus
B_{\lambda/4}(p)}({d_{A(\q)}^{\ast}}^{\epsilon}{F_{A(\q)}}^{\epsilon},\alpha)\,dx
+\int_{B_{\lambda/4}(p)}({d_{A(\q)}^{\ast}}^{\epsilon}{F_{A(\q)}}^{\epsilon},\alpha)\,dx\,.
\end{align}
The last term in \eqref{nabla} vanishes since $A(\q)$ is Yang Mills
on $B_{\lambda/4}(p)$. We now proceed estimating the remaining three
terms.

\textit{Estimate on $B^4\setminus B_{2\lambda}(p)$:} On
$B^4\setminus B_{2\lambda}(p)$, we write
$A(\q)=\underline{A}_{\epsilon}+\frac{1}{\epsilon}gPI^2_{\lambda,p}g^{-1}:=A_1+A_2$
and
\begin{align}
\label{est1}
{d_{A(\q)}^{\ast}}^{\epsilon}{F_{A(\q)}}^{\epsilon}&={d_{A_1}^{\ast}}^{\epsilon}{F_{A_1}}^{\epsilon}+{d_{A_2}^{\ast}}^{\epsilon}{F_{A_2}}^{\epsilon}+\epsilon\ast
d\ast[A_1,A_2]+\epsilon\ast[A_1,\ast {F_{A_2}}^{\epsilon}]+
\epsilon\ast[A_2,\ast {F_{A_1}}^{\epsilon}]\notag\\
 &\quad+\epsilon^2\ast[A_1,\ast[A_1,A_2]]+\epsilon^2\ast[A_2,\ast[A_1,A_2]].
\end{align}

\noindent Since $A_1=\underline{A}_{\epsilon}$ is Yang Mills, we
have ${d_{A_1}^{\ast}}^{\epsilon}F_{A_1}^{\epsilon}=0$. Also, since
${F_{A_2}}^{\epsilon}=\frac{1}{\epsilon}g(F_{I^2_{\lambda,p}}-d_{I^2_{\lambda,p}}h_{\lambda,p}+\frac{1}{2}[h_{\lambda,p},h_{\lambda,p}])g^{-1}$
and $I^2_{\lambda,p}$ is Yang Mills, we have
\begin{align}
\label{distar}
{d_{A_2}^{\ast}}^{\epsilon}{F_{A_2}}^{\epsilon}&=\frac{1}{\epsilon}gd_{I^2_{\lambda,p}}^{\ast}d_{I^2_{\lambda,p}}h_{\lambda,p}g^{-1}-
\frac{1}{2\epsilon}gd_{I^2_{\lambda,p}}^{\ast}[h_{\lambda,p},h_{\lambda,p}]g^{-1}-
\frac{1}{\epsilon}\ast g[h_{\lambda,p},\ast F_{I^2_{\lambda,p}}]g^{-1}\notag\\
 &\quad+\frac{1}{\epsilon}\ast g[h_{\lambda,p},\ast d_{I^2_{\lambda,p}}h_{\lambda,p}]g^{-1}-
 \frac{1}{2\epsilon}\ast g[h_{\lambda,p},\ast[h_{\lambda,p},h_{\lambda,p}]]g^{-1}.
\end{align}
For the first term in \eqref{distar} we compute explicitly
\begin{align}
\label{ft}
d_{I^2_{\lambda,p}}^{\ast}d_{I^2_{\lambda,p}}h_{\lambda,p}&=d^{\ast}dh_{\lambda,p}+d^{\ast}[I^2_{\lambda,p},h_{\lambda,p}]+
\ast[I^2_{\lambda,p},\ast d_{\lambda,p}]+\ast[I^2_{\lambda,p},\ast[I^2_{\lambda,p},h_{\lambda,p}]]\notag\\
 &=-dd^{\ast}h_{\lambda,p}+\ast[I^2_{\lambda,p},\ast[I^2_{\lambda,p},h_{\lambda,p}]]+
d^{\ast}[I^2_{\lambda,p},h_{\lambda,p}]+\ast[I^2_{\lambda,p},\ast
dh_{\lambda,p}],
\end{align}
where we have used the harmonicity of $h_{\lambda,p}$ (i.e.,
$d^{\ast}dh_{\lambda,p}+dd^{\ast}h_{\lambda,p}=0$). Since
${d_{A(\q)}^{\ast}}^{\epsilon}\alpha=0$, we have
$d^{\ast}\alpha=-\epsilon\ast[A(\q),\ast\alpha]$ and
\begin{align}
\label{begin} \biggl|\int_{B^4\setminus
B_{2\lambda}(p)}(-dd^{\ast}h_{\lambda,p},\alpha)\,dx\biggr|
&\lesssim\biggl|\int_{B^4}(d^{\ast}h_{\lambda,p},d^{\ast}\alpha)\,dx
\biggr|+
\biggr|\int_{\partial B_{2\lambda}(p)}(d^{\ast}h_{\lambda,p},\alpha)\,dx\biggr|\notag\\
&\lesssim\int_{B^4\setminus B_{2\lambda}(p)}\epsilon|\nabla
h_{\lambda,p}||A(\q)||\alpha|\,dx+
\lambda^2\int_{\partial B_{2\lambda}(p)}|\alpha|\,dx\notag\\
&\lesssim\,\lambda^2\int_{B^4\setminus
B_{2\lambda}(p)}(\epsilon+|I^2_{\lambda,p}|)|\alpha|\,dx+
\lambda^4\biggl(\int_{\partial B_{2\lambda}(p)}|\alpha|^3\,dx\biggr)^{1/3}\notag\\
&\lesssim\,\lambda^2\epsilon\|\alpha\|_{A(\q);1,2}+
\lambda^2\biggl(\int_{B^4\setminus
B_{2\lambda}(p)}|I^2_{\lambda,p}|^{4/3}\,dx\biggr)^{3/4}
\biggl(\int_{B^4\setminus B_{2\lambda}(p)}|\alpha|^4\,dx\biggr)^{1/4}\notag\\
&\quad+\lambda^4\|\alpha\|_{A(\q);1,2}\lesssim\epsilon^2|\log\epsilon|^{3/4}\|\alpha\|_{A(q);1,2}\,,
\end{align}
where we have integrated by parts in the first line, used Lemma
\ref{L5.6} in the second line, H\"older's inequality in the third
line, and Sobolev and trace inequalities in the last line
($L^2_1\subset L^4$ in dimension  4, $\|\alpha\|_{L^3(\partial
B_{2\lambda}(p))}\le C\||\alpha|\|_{L^2_{1/2}(\partial
B_{2\lambda}(p))}\le C\||\alpha|\|_{L^2_1(B_{2\lambda}(p))}\le
C\|\alpha\|_{A(\q);1,2}).$

The remaining terms in \eqref{ft} are also estimated using H\"older
and Sobolev inequalities and Lemma \ref{L5.6} as follows:
\begin{align}
\biggl|\int_{B^4\setminus
B_{2\lambda}(p)}(\ast[I^2_{\lambda,p},\ast[I^2_{\lambda,p},h_{\lambda,p}]],\alpha)\,dx
\biggr|&\lesssim\lambda^2\biggl(\int_{B^4\setminus
B_{2\lambda}(p)}|I^2_{\lambda,p}|^{8/3}\,dx&\biggr)^{3/4}
\biggl(\int_{B^4\setminus B_{2\lambda}(p)}|\alpha|^4\,dx\biggr)^{1/4}\notag\\
 &\lesssim\epsilon^{3/2}\|\alpha\|_{A(\q);1,2}\,;
\end{align}
\begin{align}
 &\biggl|\int_{B^4\setminus
B_{2\lambda}(p)}(d^{\ast}[I^2_{\lambda,p},h_{\lambda,p}],\alpha)\,dx
\bigg|\lesssim\int_{B^4\setminus B_{2\lambda}(p)}\lambda^2(|\nabla I^2_{\lambda,p}|+|I^2_{\lambda,p}|)|\alpha|\,dx\notag\\
 &\lesssim\lambda^2\biggl[\biggl(\int_{B^4\setminus B_{2\lambda}(p)}|\nabla I^2_{\lambda,p}|^{4/3}\,dx\biggr)^{3/4}
+\biggl(\int_{B^4\setminus
B_{2\lambda}(p)}|I^2_{\lambda,p}|^{4/3}\,dx\biggr)^{3/4}\biggr]
\biggl(\int_{B^4\setminus B_{2\lambda}(p)}|\alpha|^4\,dx\biggr)^{1/4}\notag\\
&\lesssim\lambda^3\|\alpha\|_{A(\q);1,2}+\lambda^4|\log\lambda|^{3/4}\|\alpha\|_{A(\q);1,2}
\lesssim \,\epsilon^{3/2}\|\alpha\|_{A(\q);1,2}\,;
\end{align}
\begin{align}
 \biggl|\int_{B^4\setminus
B_{2\lambda}(p)}(\ast[I^2_{\lambda,p},\ast
dh_{\lambda,p}],\alpha)\,dx \biggl|
&\lesssim\lambda^2\int_{B^4\setminus B_{2\lambda}(p)}|I^2_{\lambda,p}||\alpha|\,dx\notag\\
 &\lesssim\lambda^2\biggl(\int_{B^4\setminus B_{2\lambda}(p)}|I^2_{\lambda,p}|^{4/3}\,dx\biggr)^{3/4}
\biggl(\int_{B^4\setminus B_{2\lambda}(p)}|\alpha|^4\,dx\biggr)^{1/4}\notag\\
 &\lesssim\lambda^4|\log\lambda|\|\alpha\|_{A(\q);1,2}\lesssim\epsilon^2|\log\epsilon|^{3/4}\|\alpha\|_{A(\q);1,2}\,.
\end{align}
The remaining terms in \eqref{distar} are estimated by Lemma
\ref{L5.6} as
\begin{align}
&\biggl|\int_{B^4\setminus B_{2\lambda}(p)}(gd_{I^2_{\lambda,p}}^{\ast}[h_{\lambda,p},h_{\lambda,p}]g^{-1},\alpha)\,dx\biggr|\notag\\
&\lesssim\lambda^4\int_{B^4\setminus B_{2\lambda}(p)}|\alpha|\,dx+
\lambda^4\int_{B^4\setminus B_{2\lambda}(p)}|I^2_{\lambda,p}||\alpha|\,dx\notag\\
&\lesssim\lambda^4\|\alpha\|_{A(\q);1,2}+\lambda^4\biggl(\int_{B^4\setminus
B_{2\lambda}(p)}|I^2_{\lambda,p}|^{4/3}\,dx\biggr)^{3/4}
\biggl(\int_{B^4\setminus B_{2\lambda}(p)}|\alpha|^4\,dx\biggr)^{1/4}\notag\\
&\lesssim\lambda^4\|\alpha\|_{A(\q);1,2}+\lambda^6|\log\lambda|^{3/4}\|\alpha\|_{A(\q);1,2}\lesssim\epsilon^2\|\alpha\|_{A(\q);1,2}\,;
\end{align}
\begin{align}
\biggl|\int_{B^4\setminus B_{2\lambda}(p)}(g\ast[h_{\lambda,p},\ast
F_{I^2_{\lambda,p}}]g^{-1},\alpha)\,dx\biggr|
&\lesssim\lambda^2\int_{B^4\setminus B_{2\lambda}(p)}|F_{I^2_{\lambda,p}}||\alpha|\,dx\notag\\
&\lesssim\lambda^2\biggl(\int_{B^4\setminus
B_{2\lambda}(p)}|F_{I^2_{\lambda,p}}|^{4/3}\,dx\biggr)^{3/4}
\biggl(\int_{B^4\setminus B_{2\lambda}(p)}|\alpha|^4\,dx\biggr)^{1/4}\notag\\
&\lesssim\lambda^3\|\alpha\|_{A(\q);1,2}\lesssim\epsilon^{3/2}\|\alpha\|_{A(\q);1,2}\;;
\end{align}
\begin{align}
&\biggl|\int_{B^4\setminus B_{2\lambda}(p)}(g\ast[h_{\lambda,p},\ast
d_{I^2_{\lambda,p}}h_{\lambda,p}]g^{-1},\alpha)\,dx
\biggr|\notag\\
&\lesssim\lambda^4\int_{B^4\setminus B_{2\lambda}(p)}|\alpha|\,dx+
\lambda^4\int_{B^4\setminus
B_{2\lambda}(p)}|I^2_{\lambda,p}||\alpha|\,dx\lesssim\lambda^4\|\alpha\|_{A(\q);1,2}\lesssim\epsilon^2\|\alpha\|_{A(\q);1,2}\,;
\end{align}
\begin{equation}
\label{end}
 \biggl|\int_{B^4\setminus
B_{2\lambda}(p)}(g\ast[h_{\lambda,p},\ast[h_{\lambda,p},h_{\lambda,p}]],\alpha)\,dx
\biggr|\lesssim\lambda^6\|\alpha\|_{A(\q);1,2}\lesssim\epsilon^3\|\alpha\|_{A(\q);1,2}.
\end{equation}
From \eqref{begin}--\eqref{end}, we obtain
\begin{equation}
\label{begin2} \biggl|\int_{B^4\setminus
B_{2\lambda}(p)}({d_{A_2}^{\ast}}^{\epsilon}{F_{A_2}}^{\epsilon},\alpha)\,dx
\biggr|\lesssim\epsilon^{1/2}\|\alpha\|_{A(\q);1,2}.
\end{equation}
The remaining terms of \eqref{est1} are estimated as follows:
\begin{align}
&\biggl|\int_{B^4\setminus B_{2\lambda}(p)}(\epsilon\ast
d\ast[A_1,A_2],\alpha)\,dx\biggr| \lesssim
\epsilon\int_{B^4\setminus B_{2\lambda}(p)}|\nabla
A_1||A_2||\alpha|dx+
\epsilon\int_{B^4\setminus B_{2\lambda}(p)}|A_1||\nabla A_2||\alpha|\,dx\notag\\
&\lesssim\int_{B^4\setminus B_{2\lambda}(p)}(|I^2_{\lambda,p}|+
\lambda^2)|\alpha|\,dx+\int_{B^4\setminus B_{2\lambda}(p)}(|\nabla I^2_{\lambda,p}|+\lambda^2)|\alpha|\,dx\notag\\
&\lesssim\lambda^2\|\alpha\|_{A(\q);1,2}+\biggl(\int_{B^4\setminus
B_{2\lambda}(p)}|I^2_{\lambda,p}|^{4/3}\,dx\biggr)^{3/4}
\|\alpha\|_{A(\q);1,2}+\biggl(\int_{B^4\setminus B_{2\lambda}(p)}|\nabla I^2_{\lambda,p}|^{4/3}\,dx\biggr)^{3/4}\|\alpha\|_{A(\q);1,2}\notag\\
&\lesssim
\lambda^2\|\alpha\|_{A(\q);1,2}+\lambda^2|\log\lambda|^{3/4}\|\alpha\|_{A(\q);1,2}+\lambda\|\alpha\|_{A(\q);1,2}
\lesssim\epsilon^{1/2}\|\alpha\|_{A(\q);1,2}\,;
\end{align}
\begin{align}
 &\biggl|\int_{B^4\setminus B_{2\lambda}(p)}(\epsilon\ast[A_1,\ast {F_{A_2}}^{\epsilon}],\alpha)\,dx\biggr|
\lesssim \epsilon\int_{B^4\setminus
B_{2\lambda}(p)}|{F_{A_2}}^{\epsilon}||\alpha|\,dx
\lesssim \int_{B^4\setminus B_{2\lambda}(p)}(|F_{I^2_{\lambda,p}}|+|d_{I^2_{\lambda,p}}h_{\lambda,p}|+|h_{\lambda,p}|^2)|\alpha|\,dx\notag\\
&\lesssim \biggl(\int_{B^4\setminus
B_{2\lambda}(p)}|F_{I^2_{\lambda,p}}|^{4/3}\,dx\biggr)^{3/4}\|\alpha\|_{A(\q);1,2}+
\lambda^2\|\alpha\|_{A(\q);1,2}+\lambda^2\biggl(\int_{B^4\setminus
B_{2\lambda}(p)}|I^2_{\lambda,p}|^{4/3}\,dx\biggr)^{3/4}\|\alpha\|_{A(\q);1,2}
\notag\\
 &\quad+\lambda^4\|\alpha\|_{A(\q);1,2}\lesssim\lambda\|\alpha\|_{A(\q);1,2}\lesssim\epsilon^{1/2}\|\alpha\|_{A(\q);1,2}\,;
\end{align}
\begin{align}
\biggl|\int_{B^4\setminus B_{2\lambda}(p)}(\epsilon\ast[A_2,\ast
{F_{A_1}}^{\epsilon}],\alpha)\,dx\biggr|&\lesssim\epsilon\int_{B^4\setminus
B_{2\lambda}(p)}|A_2||\alpha|\,dx
 \lesssim\int_{B^4\setminus B_{2\lambda}(p)}(|I^2_{\lambda,p}|+\lambda^2)|\alpha|\,dx\notag\\
 &\lesssim\biggl(\int_{B^4\setminus B_{2\lambda}(p)}|I^2_{\lambda,p}|^{4/3}\,dx\biggr)^{3/4}\|\alpha\|_{A(\q);1,2}+
 \lambda^2\|\alpha\|_{A(\q);1,2}\notag\\
 &\lesssim\lambda^2|\log\lambda|^{3/4}\|\alpha\|_{A(\q);1,2}\lesssim\epsilon|\log\epsilon|^{3/4}\|\alpha\|_{A(\q);1,2}\,;
\end{align}
\begin{align}
\biggl|\int_{B^4\setminus
B_{2\lambda}(p)}(\epsilon^2\ast[A_1,\ast[A_1,A_2]],\alpha)\,dx\biggr|
&\lesssim\epsilon\int_{B^4\setminus B_{2\lambda}(p)}(|I^2_{\lambda,p}|+\lambda^2)|\alpha|\,dx\notag\\
&\lesssim\epsilon\biggl(\int_{B^4\setminus
B_{2\lambda}(p)}|I^2_{\lambda,p}|^{4/3}\,dx\biggr)^{3/4}\|\alpha\|_{A(\q);1,2}+
 \epsilon\lambda^2\|\alpha\|_{A(\q);1,2}\notag\\
&\lesssim\epsilon\lambda^2|\log\lambda|^{3/4}\|\alpha\|_{A(\q);1,2}\lesssim\epsilon^2|\log\epsilon|^{3/4}\|\alpha\|_{A(\q);1,2}\,;
\end{align}
\begin{align}
\label{end2}&\biggl|\int_{B^4\setminus
B_{2\lambda}(p)}(\epsilon^2\ast[A_2,\ast[A_1,A_2]],\alpha)\,dx\biggr|
\lesssim\epsilon^2\int_{B^4\setminus B_{2\lambda}(p)}|A_2|^2|\alpha|\,dx\notag\\
&\lesssim\int_{B^4\setminus
B_{2\lambda}(p)}(|I^2_{\lambda,p}|^2+\lambda^4)|\alpha|\,dx
 \lesssim\biggl(\int_{B^4\setminus B_{2\lambda}(p)}|I^2_{\lambda,p}|^{8/3}\,dx\biggr)^{3/4}\|\alpha\|_{A(\q);1,2}+
 \lambda^4\|\alpha\|_{A(\q);1,2}\notag\\
&\lesssim\lambda\|\alpha\|_{A(\q);1,2}\lesssim\epsilon^{1/2}\|\alpha\|_{A(\q);1,2}.
\end{align}
Combining \eqref{begin2}--\eqref{end2}, we obtain
\begin{equation}
\label{est2} \biggl|\int_{B^4\setminus
B_{2\lambda}(p)}({d_{A(\q)}^{\ast}}^{\epsilon}{F_{A(\q)}}^{\epsilon},\alpha)
\biggr|\lesssim\epsilon^{1/2}\|\alpha\|_{A(\q);1,2}.
\end{equation}

\medskip

\textit{Estimate on $B_{2\lambda}(p)\setminus B_{\lambda/2}(p)$:} On
$B_{2\lambda}(p)\setminus B_{\lambda/2}(p)$, we write
$A(\q)=(1-\beta_{\lambda,p})\underline{A}_{\epsilon}+\frac{1}{\epsilon}gPI^2_{\lambda,p}g^{-1}=:A_1+A_2$,
and make use again of the expansion \eqref{est1}.

Since
${F_{A_1}}^{\epsilon}=-d\beta_{\lambda,p}\wedge\underline{A}_{\epsilon}+(1-\beta_{\lambda,p})d\underline{A}_{\epsilon}+
\frac{\epsilon}{2}(1-\beta_{\lambda,p}^2)[\underline{A}_{\epsilon},\underline{A}_{\epsilon}]$,
we easily see that
$|{d_{A_1}^{\ast}}^{\epsilon}{F_{A_1}}^{\epsilon}|\lesssim\lambda^{-2}$,
and
\begin{align}
\biggl|\int_{B_{2\lambda}(p)\setminus
B_{\lambda/2}(p)}({d_{A_1}^{\ast}}^{\epsilon}{F_{A_1}}^{\epsilon},\alpha)\,dx&\biggr|
\lesssim\lambda^{-2}\int_{B_{2\lambda}(p)\setminus B_{\lambda/2}(p)}|\alpha|\,dx\notag\\
 &\lesssim\lambda\biggl(\int_{B_{2\lambda}(p)\setminus
B_{\lambda/2}(p)}|\alpha|^4\,dx\biggr)^{1/4}\lesssim\epsilon^{1/2}\|\alpha\|_{A(\q);1,2}\,.\notag
\end{align}
The term ${d_{A_2}^{\ast}}^{\epsilon}{F_{A_2}}^{\epsilon}$ is given
by the formula \eqref{distar}, thus we may proceed as we did earlier
and obtain
\begin{equation}
\biggl|\int_{B_{2\lambda}(p)\setminus
B_{\lambda/2}(p)}({d_{A_2}^{\ast}}^{\epsilon}{F_{A_2}}^{\epsilon},\alpha)\,
\biggr|\lesssim\epsilon^{1/2}\|\alpha\|_{A(\q);1,2}\,.\notag
\end{equation}
The estimates of all the remaining terms in \eqref{est1} are quite
similar and they yield
\begin{equation}
\label{begin2bis} \biggl|\int_{B_{2\lambda}(p)\setminus
B_{\lambda/2}(p)}({d_{A(\q)}^{\ast}}^{\epsilon}{F_{A(\q)}}^{\epsilon},\alpha)\,dx
\biggr|\lesssim\epsilon^{1/2}\|\alpha\|_{A(\q);1,2}.
\end{equation}

\medskip

\textit{Estimate on $B_{\lambda/2}(p)\setminus B_{\lambda/4}(p)$:}

\noindent On $B_{\lambda/2}(p)\setminus B_{\lambda/4}(p)$, we write
$A(\q)=\frac{1}{\epsilon}gI^2_{\lambda,p}g^{-1}+
\frac{1}{\epsilon}(\beta_{\lambda,p}-1)gh_{\lambda,p}g^{-1}=:A_2+A_3$,
and use the expansion
\begin{align}
\label{3.49}
{d_{A(\q)}^{\ast}}^{\epsilon}{F_{A(\q)}}^{\epsilon}&={d_{A_2}^{\ast}}^{\epsilon}{F_{A_2}}^{\epsilon}+{d_{A_3}^{\ast}}^{\epsilon}{F_{A_3}}^{\epsilon}+\epsilon\ast
d\ast[A_2,A_3]+\epsilon\ast[A_2,\ast {F_{A_3}}^{\epsilon}]+
\epsilon\ast[A_3,\ast {F_{A_2}}^{\epsilon}]\notag\\
 &\quad+\epsilon^2\ast[A_2,\ast[A_2,A_3]]+\epsilon^2\ast[A_3,\ast[A_2,A_3]]\,,
\end{align}
where the first term above vanishes, since
${d_{A_2}^{\ast}}^{\epsilon}{F_{A_2}}^{\epsilon}=0$.

For the second term, we have
${F_{A_3}}^{\epsilon}=\frac{1}{\epsilon}d\beta_{\lambda,p}\wedge
gh_{\lambda,p}g^{-1}+\frac{1}{\epsilon}(\beta_{\lambda,p}-1)gdh_{\lambda,p}g^{-1}+\frac{1}{2\epsilon}
(\beta_{\lambda,p}^2-1)g[h_{\lambda,p},h_{\lambda,p}]g^{-1},$ and we
can easily see from Lemma \ref{L5.6} that
$|{d_{A_3}^{\ast}}^{\epsilon}{F_{A_3}}^{\epsilon}|\lesssim\epsilon^{-1}$.
Thus we have
\begin{align}
\label{begin3} \biggl|\int_{B_{\lambda/2}(p)\setminus
B_{\lambda/4}(p)}({d_{A_3}^{\ast}}^{\epsilon}{F_{A_3}}^{\epsilon},\alpha)\,dx
\biggr|&
\lesssim\epsilon^{-1}\int_{B_{\lambda/2}(p)\setminus B_{\lambda/4}(p)}|\alpha|\,dx\notag\\
&\lesssim\epsilon^{-1}\lambda^3\biggl(\int_{B_{\lambda/2}(p)\setminus
B_{\lambda/4}(p)}|\alpha|^4\,dx\biggr)^{1/4}
\lesssim\epsilon^{1/2}\|\alpha\|_{A(\q);1,2}.
\end{align}
The remaining terms of \eqref{3.49} are estimated as follows:
\begin{align}
 &\biggl|\int_{B_{\lambda/2}(p)\setminus
B_{\lambda/4}(p)}(\epsilon\ast d\ast[A_2,A_3],\alpha)\,dx\biggr|
\lesssim \epsilon^{-1}\int_{B_{\lambda/2}(p)\setminus
B_{\lambda/4}(p)}(|\nabla
I^2_{\lambda,p}|\lambda^2+|I^2_{\lambda,p}|\lambda)|\alpha|\,dx
\notag\\
&\lesssim
\epsilon^{-1}\lambda^2\biggl(\int_{B_{\lambda/2}(p)\setminus
B_{\lambda/4}(p)}|\nabla I^2_{\lambda,p}|^{4/3}\,dx\biggr)^{3/4}
\|\alpha\|_{A(\q);1,2}+\epsilon^{-1}\lambda\biggl(\int_{B_{\lambda/2}(p)\setminus
B_{\lambda/4}(p)}|I^2_{\lambda,p}|^{4/3}\,dx\biggr)^{3/4}
\|\alpha\|_{A(\q);1,2}\notag\\
&\lesssim
 \epsilon^{-1}\lambda^3\|\alpha\|_{A(\q);1,2}+\epsilon^{-1}\lambda^3\|\alpha\|_{A(\q);1,2}\lesssim\epsilon^{1/2}\|\alpha\|_{A(\q);1,2}\,;
\end{align}
\begin{align}
\biggl|\int_{B_{\lambda/2}(p)\setminus
B_{\lambda/4}(p)}(\epsilon\ast[A_2,\ast
{F_{A_3}}^{\epsilon}],\alpha)\,dx\biggr|&
\lesssim\epsilon^{-1}\lambda\int_{B_{\lambda/2}(p)\setminus B_{\lambda/4}(p)}|I^2_{\lambda,p}||\alpha|\,dx\notag\\
&\lesssim\epsilon^{-1}\lambda\biggl(\int_{B_{\lambda/2}(p)\setminus
B_{\lambda/4}(p)}|I^2_{\lambda,p}|^{4/3}\,dx\biggr)^{3/4}
\|\alpha\|_{A(\q);1,2}\notag\\
&\lesssim\epsilon^{-1}\lambda^3\|\alpha\|_{A(\q);1,2}\lesssim\epsilon^{1/2}\|\alpha\|_{A(\q);1,2}\,;
\end{align}
\begin{align}
\biggl|\int_{B_{\lambda/2}(p)\setminus
B_{\lambda/4}(p)}(\epsilon\ast[A_3,\ast
{F_{A_2}}^{\epsilon}],\alpha)\,dx\biggr|&
\lesssim\epsilon^{-1}\lambda^2\int_{B_{\lambda/2}(p)\setminus B_{\lambda/4}(p)}|F_{I^2_{\lambda,p}}||\alpha|\,dx\notag\\
&\lesssim\epsilon^{-1}\lambda^2\biggl(\int_{B_{\lambda/2}(p)\setminus
B_{\lambda/4}(p)}|F_{I^2_{\lambda,p}}|^{4/3}\,dx\biggr)^{3/4}
\|\alpha\|_{A(\q);1,2}\notag\\
&\lesssim\epsilon^{-1}\lambda^3\|\alpha\|_{A(\q);1,2}\lesssim\epsilon^{1/2}\|\alpha\|_{A(\q);1,2}\,;
\end{align}
\begin{align}
\biggl|\int_{B_{\lambda/2}(p)\setminus
B_{\lambda/4}(p)}(\epsilon^2\ast[A_2,\ast[A_2,A_3]],\alpha)\,dx\biggr|
&\lesssim\epsilon^{-1}\lambda^2\int_{B_{\lambda/2}(p)\setminus B_{\lambda/4}(p)}|I^2_{\lambda,p}|^2|\alpha|\,dx\notag\\
&\lesssim\epsilon^{-1}\lambda^2\biggl(\int_{B_{\lambda/2}(p)\setminus
B_{\lambda/4}(p)}|I^2_{\lambda,p}|^{8/3}\,dx\biggr)^{3/4}
\|\alpha\|_{A(\q);1,2}\notag\\
&\lesssim\epsilon^{-1}\lambda^3\|\alpha\|_{A(\q);1,2}\lesssim\epsilon^{1/2}\|\alpha\|_{A(\q);1,2}\,;
\end{align}
\begin{align}
\label{final3} \biggl|\int_{B_{\lambda/2}(p)\setminus
B_{\lambda/4}(p)}(\epsilon^2\ast[A_3,\ast[A_2,A_3]],\alpha)\,dx
\biggr|&\lesssim\epsilon^{-1}\lambda^4\int_{B_{\lambda/2}(p)\setminus B_{\lambda/4}(p)}|I^2_{\lambda,p}||\alpha|\,dx\notag\\
&\lesssim\epsilon^{-1}\lambda^4\biggl(\int_{B_{\lambda/2}(p)\setminus
B_{\lambda/4}(p)}|I^2_{\lambda,p}|^{4/3}\,dx\biggr)^{3/4}
\|\alpha\|_{A(\q);1,2}\notag\\
&\lesssim\epsilon^{-1}\lambda^6\|\alpha\|_{A(\q);1,2}\lesssim\epsilon^2\|\alpha\|_{A(\q);1,2}\,.
\end{align}
Combining \eqref{begin3}--\eqref{final3}, we obtain
\begin{equation}
\label{est3} \biggl|\int_{B_{\lambda/2}(p)\setminus
B_{\lambda/4}(p)}({d_{A(\q)}^{\ast}}^{\epsilon}{F_{A(\q)}}^{\epsilon},\alpha)\,dx
\biggr|\lesssim\epsilon^{1/2}\|\alpha\|_{A(\q);1,2}.
\end{equation}

\medskip

\textit{Estimate on $B_{\lambda/4}(p)$:}

\noindent On $B_{\lambda/4}(p)$,
$A(\q)=\frac{1}{\epsilon}gI^1_{\lambda,p}g^{-1}$ and
${d_{A(\q)}^{\ast}}^{\epsilon}{F_{A(\q)}}^{\epsilon}=0$. We thus
have
\begin{equation}
\label{est4}
\int_{B_{\lambda/4}(p)}({d_{A(\q)}^{\ast}}^{\epsilon}{F_{A(\q)}}^{\epsilon},\alpha)\,dx=0\;.
\end{equation}
From \eqref{3.49}, \eqref{begin2bis}, \eqref{est3}, \eqref{est4}, we
finally obtain
$$\biggl|\int_{B^4}({d_{A(\q)}^{\ast}}^{\epsilon}{F_{A(\q)}}^{\epsilon},\alpha) \biggr|\lesssim\epsilon^{1/2}\|\alpha\|_{A(\q);1,2}\,.$$
This completes the proof. \hfill$\Box$

\subsection{Estimate of the remainder $R(\q;a)$}
Let $A\in\mathcal{A}(A_0;\q)$ and $a\in
L^2_{1,0}(T^{\ast}B^4\otimes\Ad(P))$. In this section, we estimate
the dual norm of the remainder $R(\q;a)$, defined via the formula
$$R(\q;a):=\nabla\YMe(A(\q)+a)-\nabla\YMe(A(\q))-\nabla^2\YMe(A(\q))a\;,$$
where the Hessian of $\YMe$, denoted by $\nabla^2\YMe(A)$, is given
by
\begin{equation}
\label{3.58}
\langle\nabla^2\YMe(A)a,b\rangle:=2\int_{B^4}({d_A}^{\epsilon}a,{d_A}^{\epsilon}b)
+2\int_{B^4}({F_A}^{\epsilon},\epsilon[a,b])\;\hbox{ for all } a,
b\in L^2_{1,0}(T^{\ast}B^4\otimes\Ad(P))\,,
\end{equation}
 where
$\langle\cdot,\cdot\rangle$ denotes the pairing between
$L^2_{1,0}(T^{\ast}B^4\otimes\Ad(P))$ and its dual.
\begin{lemma}
\label{L3.2} For $\q\in\PP(d_0,\lambda_0;D_1,D_2;\epsilon)$, $a\in
L^2_{1,0}(T^{\ast}B^4\otimes\Ad(P))$, one has
$$\|R(\q;a)\|_{A(\q);1,2,\ast}\le C\epsilon(\|a\|_{A(\q);1,2}^2+\epsilon\|a\|_{A(\q);1,2}^3)\;,$$
where $C>0$ is a constant depending only on $d_0,\lambda_0,D_1,D_2$.
\end{lemma}
\textit{Proof:} For $a,b\in L^2_{1,0}(T^{\ast}B^4\otimes\Ad(P))$,
one has
\begin{align}
\notag
 \nabla\YMe(A(\q)+a)(b)&=2\int_{B^4}({F_{A(\q)+a}}^{\epsilon},{d_{A(\q)+a}}^{\epsilon}b)
=2\int_{B^4}({F_{A(\q)}}^{\epsilon}+{d_{A(\q)}}^{\epsilon}a+\frac{\epsilon}{2}[a,a],
{d_{A(\q)}}^{\epsilon}b+\epsilon[a,b])\,dx\notag\\
 &=\nabla\YMe(A(\q))(b)+\langle\nabla^2\YMe(A(\q))a,b\rangle+\langle
R(\q;a),b\rangle\;,
\end{align}
where
\begin{equation}
\notag \langle
R(\q;a),b\rangle=2\epsilon\int_{B^4}({d_{A(\q)}}^{\epsilon}a,[a,b])\,dx
+\epsilon\int_{B^4}([a,a],{d_{A(\q)}}^{\epsilon}b)\,dx+\epsilon^2\int_{B^4}([a,a],[a,b])\,dx.
\end{equation}
Therefore, by the H\"older's inequality
\begin{equation}
\label{remainder} |\langle R(\q;a),b\rangle|\le
C\epsilon\|{d_{A(\q)}}^{\epsilon}a\|_2\|a\|_4\|b\|_4
+C\epsilon\|a\|^2_4\|{d_{A(\q)}}^{\epsilon}b\|_2+C\epsilon^2\|a\|_4^3\|b\|_4
\end{equation}
Applying the Weitzenb\"ock formula,
${\Delta_A}^{\epsilon}:={d_A}^{\epsilon}{d_A^{\ast}}^{\epsilon}+{d_A^{\ast}}^{\epsilon}{d_A}^{\epsilon}={\nabla_A^{\ast}}^{\epsilon}{\nabla_A}^{\epsilon}+\epsilon\{{F_A}^{\epsilon},\cdot\}$
(cf.~\cite{DK},~\cite{FU}) (here, we only need to know that
$\{\cdot,\cdot\}$ is a bilinear form) and the Sobolev inequality
$\|a\|_4\le C\|a\|_{A;1,2}$, we obtain
\begin{equation}
\label{dA}
\|{d_{A(\q)}}^{\epsilon}a\|_2^2\le\|{\nabla_{A(\q)}}^{\epsilon}a\|_2^2+\epsilon\langle\{{F_{A(\q)}}^{\epsilon},a\},a\rangle
\le\|{\nabla_{A(\q)}}^{\epsilon}a\|^2_2+C\epsilon\|{F_{A(\q)}}^{\epsilon}\|_2\|a\|_4^2\le
C\|a\|_{A(\q);1,2}^2\;.
\end{equation}
The assertion of the lemma now follows from \eqref{remainder},
\eqref{dA}. \hfill$\Box$

\subsection{Estimate of the modified Hessian}
For $\q\in\PP(d_0,\lambda_0;D_1,D_2;\epsilon)$ and
$A\in\mathcal{A}_{+1}(A_0)$, the modified Hessian, $\mathcal{H}_A$,
is defined as a bilinear form on
$L^2_{1,0}(T^{\ast}B^4\otimes\Ad(P))$ as follows: for $a,b\in
L^2_{1,0}(T^{\ast}B^4\otimes\Ad(P))$,
\begin{equation}
\label{HA}
\mathcal{H}_A(a,b):=\frac{1}{2}\langle\nabla^2\YMe(A)a,b\rangle+({d_A^{\ast}}^{\epsilon}a,{d_A^{\ast}}^{\epsilon}b)_{L^2(B^4)}\;.
\end{equation}
With this definition, $\mathcal{H}_A$ is continuous. The following
positivity result holds for the modified Hessian:

\begin{lemma}
\label{L3.3} For $\q\in\PP(d_0,\lambda_0;D_1,D_2;\epsilon)$, there
exists a constant $C$ depending only on $d_0,\lambda_0, D_1, D_2$
such that for small $\epsilon>0$ and $a\in
L^2_{1,0}(T^{\ast}B^4\otimes\Ad(P))\cap
T_{A(\q)}\mathcal{N}(d_0,\lambda_0)^{\perp}$, there holds
$$\mathcal{H}_{A(\q)}(a,a)\ge C\|a\|_{A(\q);1,2}^2.$$
Here, $ L^2_{1,0}(T^{\ast}B^4\otimes\Ad(P))\cap
T_{A(\q)}\mathcal{N}(d_0,\lambda_0)^{\perp}$ is the orthogonal
complement of $T_{A(\q)}\mathcal{N}(d_0,\lambda_0)$ in
$L^2_{1,0}(T^{\ast}B^4\otimes\Ad(P))$.
\end{lemma}

\medskip
\noindent To prove the lemma above, one needs to introduce further
notation and some auxiliary lemmas.

\noindent We define the following family of $SU(2)_{\epsilon}$-Yang
Mills connections on $\R^4$: for $\q\in\PP(d_0,\lambda_0)$,
\begin{equation}
\label{tildeA(q)} \tilde{A}(\q)=\left\{\begin{array}{ll}
 \frac{1}{\epsilon}gI^2_{\lambda,p}g^{-1}&\quad\mbox{in } \R^4\setminus\{p\}\\
 \frac{1}{\epsilon}gI^1_{\lambda,p}g^{-1}&\quad\mbox{in }
 B^4_{\lambda/4}(p)\,.
 \end{array}\right.
 \end{equation}
 These connections (absolute minimizers for the Yang Mills functional)
 live on the bundles $\tilde{P}(\q)$, defined by the
 data
 \begin{equation}
 \label{Ptilde(q)}
 \biggl(\{B_{\lambda/4}(p), \;\R^4\setminus\{p\}\}, \;\{gg_{12,p}g^{-1}\}\biggr)\,.
 \end{equation}
 Note that these bundles are
 extensions to $\R^4$ of the bundles $P(\q)$ defined in $\S3.1$ ($\tilde{P}(\q)|_{B^4}=P(\q)$) and that $\tilde{A}(\q)=A(\q)$ on $B_{\lambda/4}(p).$
Similarly to what has been done in $\S3.1$,  we set \begin{equation}
\label{tildeN} \mathcal{\tilde
 N}(d_0,\lambda_0):=\{\tilde{A}(\q):\q\in\PP(d_0,\lambda_0)\}\,,
 \end{equation}
 and apply the
convention that everything is pulled back to the
 bundle $\tilde P(\q_0)$ by the bundle isomorphisms
 $\tilde{\varphi}(\q): \tilde P(\q_0)\overset{\sim}{\to}\tilde
 P(\q)$. To the bundles and connections just defined over $\R^4$, correspond bundles and
 connections on $S^4$ (by pull back under the stereographic
 projection from the north pole), which we denote by $P(\q)_{S^4}$ and
 $A(\q)_{S^4}$, respectively.

 \medskip

To prove Lemma \ref{L3.3}, we describe the tangent space of
$\mathcal N(d_0,\lambda_0)$ at $A(\q)$ by
\begin{equation}
\label{TN}
T_{A(\q)}\mathcal{N}(d_0,\lambda_0):=\text{span}\Big\langle\frac{\partial
A}{\partial p_i}(\q),\frac{\partial A}{\partial\xi_j[g]}(\q),
\frac{\partial A}{\partial\lambda}(\q)\Big\rangle_{1\le i\le 4,1\le
j\le 3}\subset L^2_1(T^{\ast}B^4\otimes\Ad(P))\;.
\end{equation}
This is a finite dimensional (at most 8 dimensional) space. Indeed,
for small $\epsilon>0$, it is exactly $8$-dimensional (cf. Lemma 3.2
in \cite{IM3}).

Likewise, the tangent space of $\tilde{N}(d_0,\lambda_0)$ at
$\tilde{A}(\q)$ is described by
\begin{equation}
\label{TtildeN}
T_{\tilde{A}(\q)}\tilde{N}(d_0,\lambda_0):=\text{span}\Big\langle\frac{\partial}{\partial
p_i}\tilde{A}(\q),\frac{\partial}
{\partial\xi_j[g]}\tilde{A}(\q),\frac{\partial}{\partial\lambda}\tilde{A}(\q)\Big\rangle_{1\le
i\le 4,1\le j\le 3} \subset
L^2_{1;\tilde{A}(\q)}(T^{\ast}\R^4\otimes\Ad(\tilde{P}))\;.
\end{equation}

For $\q:=(p,[g],\lambda)\in\PP(d_0,\lambda_0;D_1,D_2;\epsilon)$, we
introduce the function space
\begin{equation*}
L^2_{1;\tilde{A}(\q)}(T^{\ast}\R^4\otimes\Ad(\tilde{P})) =\big\{a\in
L^1_{\text{loc}}(T^{\ast}\R^4\otimes\Ad(\tilde{P})):\int_{\R^4}|{\nabla_{\tilde{A}(\q)}}^{\epsilon}a|^2+
\frac{|a|^2}{(1+|x|^2)^2}\,dx<+\infty\big\}.
\end{equation*}
Thus $a\in L^2_{1;\tilde{A}(\q)}(T^{\ast}\R^4\otimes\Ad(\tilde{P}))$
if and only if its pull-back $\pi^{\ast}a$ is in
$L^2_1(T^{\ast}S^4)$, where $\pi:S^4\to\R^4\cup\{\infty\}$ is the
stereographic projection from the north pole. For technical reasons,
we define the weighted inner product on
$L^2_{1;\tilde{A}(\q)}(T^{\ast}\R^4\otimes\Ad(\tilde{P}))$ by
\begin{equation}
\label{wip}
(\alpha,\beta)_{\tilde{A}(\q);1,2;\R^4}:=\int_{\R^4}({\nabla_{\tilde{A}(\q)}}^{\epsilon}\alpha,{\nabla_{\tilde{A}(\q)}}^{\epsilon}\beta)+
w(x)(\alpha,\beta)\,dx,
\end{equation}
where $w(x)=1$ for $|x|\le 1$, and $w(x)=1/(1+|x|^2)^2$ for $|x|>
1$.

\noindent In the following proofs, we denote by
$T_{\tilde{A}(\q)}\tilde{\mathcal N}(d_0,\lambda_0)^{\perp}$ the
orthogonal complement of $T_{\tilde{A}(\q)}\tilde{\mathcal
N}(d_0,\lambda_0)$ in
$L^2_{1;\tilde{A}(\q)}(T^{\ast}\R^4\otimes\Ad(\tilde{P}))$ and
define
\begin{equation}
\label{AA}
\tilde{\mathcal{A}}^p_{k,+1}:=\{A\in\tilde{P}:\text{$\pi^*A$ is an
$L^p_k$-connection on $P(\q)_{S^4}$}\},
\end{equation}
where $\pi$ is the stereographic projection. We denote by
$\tilde{\mathcal{G}}^p_{k+1}$ the group of gauge transformations on
$\tilde{P}$ which come from $L^p_{k+1}$-gauge transformations on the
bundle $P(\q)_{S^4}$, and we define
\begin{equation}
\label{Btilde}
\tilde{\mathcal{B}}^p_{k,+1}:=\tilde{A}^p_{k,+1}/\tilde{\mathcal{G}}^p_{k+1}\,.
\end{equation}
In order to prove Lemma \ref{L3.3} we need the following lemma:

\begin{lemma}
\label{L3.4} For $\q\in\PP(d_0,\lambda_0;D_1,D_2;\epsilon)$, there
exists $C>0$ such that for $a\in T_{\tilde{A}(\q)} \tilde{\mathcal
N}(d_0,\lambda_0)^{\perp}$, one has
\begin{equation}
\notag
\int_{\R^4}|{d_{\tilde{A}(\q)}}^{\epsilon}a|^2\,dx+\int_{\R^4}({F_{\tilde{A}(\q)}}^{\epsilon},\epsilon[a,a])\,dx+
\int_{\R^4}|{d_{\tilde{A}(\q)}^{\ast}}^{\epsilon}a|^2\,dx \ge
C\Big(\|{\nabla_{\tilde{A}(\q)}}^{\epsilon}a\|^2_{L^2(\R^4)}+\Big\|\frac{a}{1+|x|^2}\Big\|^2_{L^2(\R^4)}\Big).
\end{equation}
\end{lemma}
\textit{Proof:} We only prove the assertion for the case
$\epsilon=1$. The general case follows by the Lie algebra
isomorphism
$\phi_{\epsilon}:\mathfrak{su}_{\epsilon}(2)\to\mathfrak{su}(2)$.

Recall that the instanton $\tilde{A}(\q)$ is action minimizing,
therefore the Hessian of $\YM$ at $\tilde{A}(\q)$ is non-negative,
i.e.,
\begin{equation}
\label{nonneg}
\int_{\R^4}|d_{\tilde{A}(\q)}a|^2\,dx+\int_{\R^4}(F_{\tilde{A}(\q)},[a,a])\,dx\ge
0,
\end{equation}
for all $a\in
L^2_{1;\tilde{A}(\q)}(T^{\ast}\R^4\otimes\Ad(\tilde{P}))$.

We claim:
\newtheorem{claim}{Claim}[section]
\begin{claim}
\label{Cl3.1} For $a\in
T_{\tilde{A}(\q)}\tilde{N}(d_0,\lambda_0)^{\perp}$ with $a\not\equiv
0$,
$$\int_{\R^4}|d_{\tilde{A}(\q)}a|^2\,dx+\int_{\R^4}(F_{\tilde{A}(\q)},[a,a])\,dx+
\int_{\R^4}|d_{\tilde{A}(\q)}^{\ast}a|^2\,dx>0\,.$$
\end{claim}
\textit{Proof of Claim \ref{Cl3.1}:} We already know that this is
non-negative by \eqref{nonneg}. By contradiction, assume that there
exists $a\in T_{\tilde{A}(\q)}\tilde{N}(d_0,\lambda_0)^{\perp}$,
$a\not\equiv 0$ such that
\begin{equation}
\label{equal}
\int_{\R^4}|d_{\tilde{A}(\q)}a|^2\,dx+\int_{\R^4}(F_{\tilde{A}(\q)},[a,a])\,dx+\int_{\R^4}|d_{\tilde{A}(\q)}^{\ast}a|^2\,dx=0.
\end{equation}
From \eqref{nonneg} and \eqref{equal}, it follows that
$d_{\tilde{A}(\q)}^{\ast}a=0$ and
\begin{equation}
\notag
\int_{\R^4}|d_{\tilde{A}(\q)}a|^2\,dx+\int_{\R^4}(F_{\tilde{A}(\q)},[a,a])\,dx
=0\,.
\end{equation}
It follows that $a$ minimizes the quadratic functional
$\alpha\mapsto\int_{\R^4}|d_{\tilde{A}(\q)}\alpha|^2\,dx+\int_{\R^4}(F_{\tilde{A}(\q)},[\alpha,\alpha])\,dx$
in $L^2_{1;\tilde{A}(\q)}(T^{\ast}\R^4\otimes\Ad(\tilde{P}))$, thus
its first variation computed at $a$ is zero, that is, for all
$\varphi\in
L^2_{1;\tilde{A}(\q)}(T^{\ast}\R^4\otimes\Ad(\tilde{P}))$,
\begin{align*}
&
\frac{d}{d\epsilon}\Big|_{\epsilon=0}\int_{\R^4}|d_{\tilde{A}(\q)}(a+\epsilon\varphi)|^2\,dx+
\int_{\R^4}(F_{\tilde{A}(\q)},[a+\epsilon\varphi,a+\epsilon\varphi])\,dx+
\int_{\R^4}|d_{\tilde{A}(\q)}(a+\epsilon\varphi)|^2\,dx\\
&=2\int_{\R^4}(d_{\tilde{A}(\q)}a,d_{\tilde{A}(\q)}\varphi)\,dx+
\int_{\R^4}(F_{\tilde{A}(\q)},[a,\varphi]+[\varphi,a])\,dx+
2\int_{\R^4}(d_{\tilde{A}(\q)}^{\ast}a,d_{\tilde{A}(\q)}^{\ast}\varphi)\,dx=0\,,
\end{align*}
or, equivalently (since $d_{\tilde{A}(\q)}^{\ast}a=0$),
\begin{equation}
\label{equal2} \nabla^2\YM(\tilde{A}(\q))(a)=0.
\end{equation}
It follows from the elliptic regularity theory that $a\in
C^{\infty}(T^{\ast}\R^4\otimes\Ad(\tilde{P}))\cap
L^2_{1;\tilde{A}(\q)}(T^{\ast}\R^4\otimes\Ad(\tilde{P}))$. It is
well-known that the set of all solutions of \eqref{equal2}, which
satisfy $d_{\tilde{A}(\q)}^{\ast}a=0$ constitute the tangent space
of the $1$-instanton moduli space $\mathcal{M}_{+1}(S^4)$ over
$\R^4\cup\{\infty\}=S^4$. One has
$T_{\tilde{A}(\q)}\mathcal{M}_{+1}(S^4)\subset
T_{\tilde{A}(\q)}\mathcal{B}^p_{1,+1}\subset\tilde{S}^p_{1,+1;\tilde{A}(\q)}:=\{a\in
L^2_{1;\tilde{A}(\q)}(T^{\ast}\R^4\otimes\Ad(\tilde{P})):d_{\tilde{A}(\q)}^{\ast}a=0\}$.
Since $T_{\tilde{A}(\q)}\tilde{N}(d_0,\lambda_0)$ coincides with the
tangent space of $\mathcal{M}_{+1}(S^4)$ at $\tilde{A}(\q)$ (up to
infinitesimal gauge transformations), this contradicts the fact that
$a$ is orthogonal to the tangent space of $\tilde{N}(d_0,\lambda_0)$
at $\tilde{A}(\q)$. This completes the proof of Claim \ref{Cl3.1}.
\hfill$\Box$

\medskip

\textit{Completion of the proof of Lemma \ref{L3.4}:} By
contradiction, assume that there exists a sequence  $\{a_n\}\subset
T_{\tilde{A}(\q)}\tilde{\mathcal N}(d_0,\lambda)^{\perp}$ such that
\begin{equation}
\label{sequence}
\int_{\R^4}|\nabla_{\tilde{A}(\q)}a_n|^2\,dx+\int_{\R^4}w(x)|a_n|^2\,dx=1\quad
\forall n,
\end{equation}
and
\begin{equation}
\label{to0}
\int_{\R^4}|d_{\tilde{A}(\q)}a_n|^2\,dx+\int_{\R^4}(F_{\tilde{A}(\q)},[a_n,a_n])\,dx+
\int_{\R^4}|d_{\tilde{A}(\q)}^{\ast}a_n|^2\,dx\to 0\;,\;
\text{as}\;n\to\infty.
\end{equation}

By \eqref{sequence}, passing to a subsequence, we may assume that
$a_n\rightharpoonup a$ weakly in
$L^2_{1;\tilde{A}(\q)}(T^{\ast}\R^4\otimes\Ad(\tilde{P}))$ as
$n\to\infty,$ for some $a\in
L^2_{1;\tilde{A}(\q)}(T^{\ast}\R^4\otimes\Ad(\tilde{P}))$. The first
and the third integrals above are clearly lower semi-continuous with
respect to the weak convergence in
$L^2_{1;\tilde{A}(\q)}(T^{\ast}\R^4\otimes\Ad(\tilde{P}))$. We
assert that the second integral
$a\mapsto\int_{\R^4}(F_{\tilde{A}(\q)},[a,a]) $ is continuous with
respect to the weak convergence in
$L^2_{1;\tilde{A}(\q)}(T^{\ast}\R^4\otimes\Ad(\tilde{P}))$. To see
this, we write $a_n=a+b_n$ with $b_n\rightharpoonup 0$ weakly in
$L^2_{1;\tilde{A}(\q)}(T^{\ast}\R^4\otimes\Ad(\tilde{P}))$. We have
\begin{equation}
\label{98}
\int_{\R^4}(F_{\tilde{A}(\q)},[a_n,a_n])\,dx=\int_{\R^4}(F_{\tilde{A}(\q)},[b_n,b_n])\,dx+
2\int_{\R^4}(F_{\tilde{A}(\q)},[b_n,a])\,dx+\int_{\R^4}(F_{\tilde{A}(\q)},[a,a])\,dx.
\end{equation}
By the Sobolev embedding, we have (modulo passing to a subsequence)
$b_n\to 0$ in $L^p_{\loc}(\R^4)$ for any $p<4$. Fixing an arbitrary
$R>0$, we have
\begin{equation}
\notag \int_{\R^4}(F_{\tilde{A}(\q)},[b_n,b_n])\,dx=\int_{|x|\le
R}(F_{\tilde{A}(\q)},[b_n,b_n])\,dx+\int_{|x|>R}(F_{\tilde{A}(\q)},[b_n,b_n])\,dx.
\end{equation}
Both integrals on the right hand side go to $0$ as $n\to\infty$ (by
H\"older's inequality, the second one is bounded above by
$C\big(\int_{|x|>R}|F_{\tilde{A}(\q)}|^2\,dx\big)^{1/2}\|b_n\|_{1,2;\tilde{A}(\q)}^2$,
which goes to $0$ as $R\to\infty$). Thus the first integral in
\eqref{98} goes to $0$ as $n\to\infty$. With similar arguments, one
shows that the second integral in \eqref{98} also goes to zero and
the assertion is proved.

Inequality \eqref{nonneg} and \eqref{to0} then imply
\begin{equation}
\notag
\int_{\R^4}|d_{\tilde{A}(\q)}a|^2\,dx+\int_{\R^4}(F_{\tilde{A}(\q)},[a,a])
+\int_{\R^4}|d_{\tilde{A}(\q)}^{\ast}a|^2\,dx=0\;,
\end{equation}
and, finally, since $a\in T_{\tilde{A}(\q)}\tilde{\mathcal
N}(d_0,\lambda_0)^{\perp}$, applying Claim \ref{Cl3.1}, one obtains
that $a=0$.

We then have, by \eqref{to0},
\begin{equation}
\label{101}
\int_{\R^4}|d_{\tilde{A}(\q)}a_n|^2\,dx+\int_{\R^4}|d_{\tilde{A}(\q)}^{\ast}a_n|^2\,dx\to
0\quad\text{for $n\to\infty$},
\end{equation}
and by the Weitzenb\"ock formula,
\begin{align}
\label{102}
\int_{\R^4}|d_{\tilde{A}(\q)}a_n|^2\,dx+\int_{\R^4}|d_{\tilde{A}(\q)}^{\ast}a_n|^2\,dx
&=
\int_{\R^4}|\nabla_{\tilde{A}(\q)}a_n|^2\,dx+(\{F_{\tilde{A}(\q)},a_n\},a_n)\,dx\notag\\
&=\int_{\R^4}|\nabla_{\tilde{A}(\q)}a_n|^2\,dx+o(1)\;\;
 (n\to\infty)\,.
\end{align}

Combining \eqref{101}, \eqref{102}, we obtain
$\int_{\R^4}|\nabla_{\tilde{A}(\q)}a_n|^2\,dx\to0$ as $n\to\infty$,
and $\int_{\R^4}w(x)|a_n|^2\,dx\to 0$, by the Sobolev embedding.
This contradicts \eqref{sequence} and completes the proof of Lemma
\ref{L3.4}. \hfill$\Box$

\medskip

To prove Lemma \ref{L3.3}, we also need to estimate the difference
between the bilinear forms $\mathcal{H}_{A(\q)}$ and
$\mathcal{H}_{\tilde{A}(\q)}$, where
$\mathcal{H}_{\tilde{A}(\q)}(a,b):=\frac{1}{2}\langle\nabla^2\YM_{\epsilon}(\tilde{A}(\q))a,b\rangle+
({d_{\tilde{A}(\q)}^{\ast}}^{\epsilon}a,{d_{\tilde{A}(\q)}^{\ast}}^{\epsilon}b)_{L^2(\R^4)}$
for $a,b\in
L^2_{1;\tilde{A}(\q)}(T^{\ast}\R^4\otimes\Ad(\tilde{P}))$. This is
the content of the following lemma.
\begin{lemma} \label{L3.5} Let $\q\in\PP(d_0,\lambda_0;D_1,D_2;\epsilon)$.
The bilinear forms $\mathcal{H}_{A}(\q)$ and
$\mathcal{H}_{\tilde{A}(\q)}$ on $L^2_{1,0}(T^{\ast}B^4\otimes\Ad
(P))\times L^2_{1,0}(T^{\ast}B^4\otimes\Ad(P))$ satisfy
$$\|\mathcal{H}_{A(\q)}-\mathcal{H}_{\tilde{A}(\q)}\|_{A(\q);1,2,\ast}\lesssim\epsilon\;.$$
(Here, $(\cdot,\cdot)_2$ denotes the $L^2$-inner product).
\end{lemma}

\textit{Proof:} Set $b=b(\q)=A(\q)-\tilde{A}(\q)$. For
$\alpha,\beta\in L^2_{1;0}(T^{\ast}B^4\otimes\Ad(P))$, a simple
computation yields
\begin{align}
\label{103}
 &(\mathcal{H}_{A(\q)}-\mathcal{H}_{\tilde{A}(\q)})(\alpha,\beta)
=\epsilon({d_{\tilde{A}(\q)}}^{\epsilon}\alpha,[b,\beta])_2+\epsilon([b,\alpha],{d_{\tilde{A}(\q)}}^{\epsilon}\beta)_2+
\epsilon(d_{A(\q)}b,[\alpha,\beta])_2+
\notag\\
&+\frac{\epsilon^2}{2}([b,b],[\alpha,\beta])_2-\epsilon({d_{\tilde{A}(\q)}^{\ast}}^{\epsilon}\alpha,\ast[b,\beta])_2-
\epsilon(\ast[b,\alpha],{d_{\tilde{A}(\q)}^{\ast}}^{\epsilon}\beta)_2+
2\epsilon^2([b,\alpha],[b,\beta])_2.
\end{align}
Since
$b=b(\q)=(1-\beta_{\lambda,p})\underline{A}_{\epsilon}+\frac{1}{\epsilon}(\beta_{\lambda/4,p}-1)gh_{\lambda,p}g^{-1}$
satisfies $\epsilon\|b\|_{\infty}\lesssim\epsilon$, integrating by
parts the third addend at the right hand side of \eqref{103}, one
obtains

\begin{equation}
\notag
\|(\mathcal{H}_{A(\q)}-\mathcal{H}_{\tilde{A}(\q)})(\alpha,\beta)\|\lesssim\epsilon
\|\alpha\|_{A(\q);1,2}\|\beta\|_{A(\q);1,2}
\end{equation}
This completes the proof of Lemma \ref{L3.5}. \hfill$\Box$

\medskip

For the next steps, we use the orthonormal basis
$\langle\A_1(\q),\A_2(\q),\ldots,\A_8(\q)\rangle$ of
$T_{A(\q)}\mathcal{N}(d_0,\lambda_0)$ given by
$\A_i(\q)=A_{\q_i}(\q)$, where the vector fields $\q_i$ ($1\le i\le
8$) on $\PP(d_0,\lambda_0)$ are defined in \cite{IM3} and
$A_{\q_i}(\q)$ denotes the directional derivative of $A(\q)$ in the
direction $\q_i(\q)$, and the basis
$\langle\tilde\A_1(\q),\tilde\A_2(\q),\ldots,\tilde\A_8(\q)\rangle$
of $T_{\tilde A(\q)}\mathcal{\tilde N}(d_0,\lambda_0)$, defined via
$\tilde{\A}_i(\q):=\tilde{A}_{\q_i}(\q)$. We also need the
orthonormal basis
$\langle\hat{\A}_1(\q),\ldots,\hat{\A}_8(\q)\rangle$ of
$T_{\tilde{A}(\q)}\tilde{N}(d_0,\lambda_0)$ constructed via the
Gram-Schmidt's orthogonalization procedure applied to
$\tilde{\A}_1(\q),\ldots,\tilde{\A}_8(\q)$ (for details, see
\cite{IM3}). One needs the following technical lemmas, proved in
\cite{IM3}.
\begin{lemma}
\label{L3.6} For $\q\in\PP(d_0,\lambda_0;D_1,D_2;\epsilon)$,
\begin{align*}\|\A_i(\q)-\tilde{\A}_i(\q)\|_{A(\q);1,2;B^4}\lesssim\epsilon^{3/2}\quad(1\le
i\le 4),\\
\|\A_i(\q)-\tilde{\A}_i(\q)\|_{A(\q);1,2;B^4}\lesssim\epsilon\quad(5\le
i\le 8)\,.\end{align*}
\end{lemma}
\begin{lemma}
\label{L3.7} For $\q\in\PP(d_0,\lambda_0;D_1,D_2;\epsilon)$,
\begin{align*}\|(\nabla^2\YMe(A(\q))-\nabla^2\YMe(\tilde{A}(\q)))\A_i(\q)\|_{A(\q);1,2,\ast}&
\lesssim\epsilon^{3/2}\quad(1\le
i\le 8),\\
\|({d_{A(\q)}}^{\epsilon}{d_{A(\q)}^{\ast}}^{\epsilon}-{d_{\tilde{A}(\q)}}^{\epsilon}{d_{\tilde{A}(\q)}^{\ast}}^{\epsilon})\A_i(\q)\|_{A(\q);1,2,\ast}&
\lesssim\epsilon^{3/2}\quad(1\le i\le 8)\,.
\end{align*}
\end{lemma}
Lemma \ref{L3.6} and the results in \cite{IM3} yield the estimate
\begin{equation}
\label{105}
\|\hat{\A}_i(\q)-\A_i(\q)\|_{A(\q);1,2}\le\|\hat{\A}_i(\q)-\tilde{\A}_i(\q)\|_{A(\q);1,2}+
\|\tilde{\A}_i(\q)-\A_i(\q)\|_{A(\q);1,2} \lesssim\epsilon\; (1\le
i\le 8)\,.
\end{equation}
In order to prove Lemma \ref{L3.3}, we need to define the
``topological projections"
\begin{equation}
\label{tpr} Q : L^2_{1,0}(T^{\ast}B^4\otimes\Ad(P))\to
L^2_{1,0}(T^{\ast}B^4\otimes\Ad(P))\cap
T_{A(\q)}\mathcal{N}(d_0,\lambda_0)^{\perp},
\end{equation}
\begin{equation}
\label{tpr1} \tilde Q :
L^2_{1;\tilde{A}(\q)}(T^{\ast}\R^4\otimes\Ad(\tilde{P}))\to
L^2_{1;\tilde{A}(\q)}(T^{\ast}\R^4\otimes\Ad(\tilde{P}))\cap
T_{\tilde{A}(\q)}\tilde{N}(d_0,\lambda_0)^{\perp}.
\end{equation}
(Recall that the orthogonal complements
$T_{A(\q)}{\mathcal{N}(d_0,\lambda_0)}^{\perp}$, $T_{\tilde
A(\q)}{\mathcal{\tilde N}(d_0,\lambda_0)}^{\perp}$ are calculated
with respect to the inner product $(\cdot,\cdot)_{A(\q);1,2;B^4}$
and $(\cdot,\cdot)_{\tilde{A}(\q);1,2;\R^4}$, respectively).

\medskip

\textit{Proof of Lemma \ref{L3.3}:} Let $a\in
L^2_{1,0}(T^{\ast}B^4\otimes\Ad(P))$. By extending $a$ trivially to
$\R^4\setminus B^4$ by $0$, we may also regard it as an element of
$L^2_{1;\tilde{A}(\q)}(T^{\ast}\R^4\otimes\Ad(\tilde{P}))$ and we
define the components
\begin{align}
\label{comp} &a^{\intercal}:=Pa\,,\quad \mbox{   with }
P:=\text{Id}-Q\,,\qquad
a^{\perp}:=Qa\,,\\
&a^{\ttilde}:=\tilde{P}a\,,\quad \mbox{   with } \tilde{P} :=
\text{Id}-\tilde{Q}\,,\qquad a^{\tperp}:=\tilde{Q}a\,.
\end{align}

We now prepare to estimate $\mathcal{H}_{A(\q)}(a,a)$. We have
\begin{align}
\label{106}
\mathcal{H}_{A(\q)}(a,a)=\mathcal{H}_{A(\q)}(a^{\intercal}+a^{\perp},a^{\intercal}+a^{\perp})=
\mathcal{H}_{A(\q)}(a^{\intercal},a^{\intercal})+2\mathcal{H}_{A(\q)}(a^{\intercal},a^{\perp})+
\mathcal{H}_{A(\q)}(a^{\perp},a^{\perp})\,.
\end{align}
The last term of \eqref{106} can be estimated by Lemma \ref{L3.5} as
follows:
\begin{align}
\label{107} \mathcal{H}_{A(\q)}(a^{\perp},a^{\perp})
=\mathcal{H}_{\tilde{A}(\q)}(a^{\perp},a^{\perp})+(\mathcal{H}_{A(\q)}-
\mathcal{H}_{\tilde{A}(\q)})(a^{\perp},a^{\perp})
\ge\mathcal{H}_{\tilde{A}(\q)}(a^{\perp},a^{\perp})-C\epsilon\|a^{\perp}\|^2_{A(\q);1,2}\,.
\end{align}
To estimate the first term of \eqref{107}, we write
$\alpha^{\perp}=(\alpha^{\perp})^{\ttilde}+(\alpha^{\perp})^{\tperp}$.
The first addend expressed in terms of the orthonormal basis
$\langle \hat{\A}_1(\q),\ldots,\hat{\A}_8(\q)\rangle$ of
$T_{\tilde{A}(\q)}\tilde{N}(d_0,\lambda_0)$ is given by
$$(a^{\perp})^{\ttilde}=\sum_{i=1}^8(a^{\perp},\hat{\A}_i(\q))_{\tilde{A}(\q);1,2;\R^4}\hat{\A}_i(\q)\,.$$
These components satisfy
\begin{align}
\label{108}
(a^{\perp},\hat{\A}_i(\q))_{\tilde{A}(\q);1,2;\R^4}&=(a^{\perp},\hat{\A}_i(\q))_{A(\q);1,2;B^4}+
((a^{\perp},\hat{\A}_i(\q))_{\tilde{A}(\q);1,2;\R^4}-(a^{\perp},\hat{\A}_i(\q))_{A(\q);1,2;B^4})\notag\\
&=(a^{\perp},\hat{\A}_i(\q)-\A_i(\q))_{A(\q);1,2;B^4}+((a^{\perp},\hat{\A}_i(\q))_{\tilde{A}(\q);1,2;\R^4}-
(a^{\perp},\hat{\A}_i(\q))_{A(\q);1,2;B^4}).
\end{align}
Here, we have (recall that $\text{supp}\,a^{\perp}\subset B^4$):
\begin{align}
\label{109}
&(a^{\perp},\hat{\A}_i(\q))_{\tilde{A}(\q);1,2;\R^4}-(a^{\perp},\hat{\A}_i(\q))_{A(\q);1,2;B^4}\notag\\
&=\int_{B^4}({\nabla_{\tilde{A}(\q)}}^{\epsilon}a^{\perp},{\nabla_{\tilde{A}(\q)}}^{\epsilon}\hat{\A}_i(\q)) -\int_{B^4}({\nabla_{A(\q)}}^{\epsilon}a^{\perp},{\nabla_{A(q)}}^{\epsilon}\hat{\A}_i(\q)) \notag\\
&=\epsilon\int_{B^4}({\nabla_{A(\q)}}^{\epsilon}a^{\perp},[b(\q),\hat{\A}_i(\q)]) +\epsilon\int_{B^4}([b(\q),a^{\perp}],{\nabla_{A(\q)}}^{\epsilon}\hat{\A}_i(\q)) \notag\\
&\quad+\epsilon^2\int_{B^4}([b(\q),a^{\perp}],[b(\q),\hat{\A}_i(\q)])
\lesssim\epsilon\|a^{\perp}\|_{A(\q);1,2;B^4}\,,
\end{align}
where $b(\q)=\tilde{A}(\q)-A(\q)$ and we have used
$\epsilon\|b\|_{\infty}\lesssim\epsilon$.

By \eqref{105}, \eqref{108}, \eqref{109}, one obtains
\begin{equation}
\label{110}
|(a^{\perp},\hat{\A}_i(\q))_{\tilde{A}(\q);1,2;\R^4}|\lesssim\epsilon\|a^{\perp}\|_{A(\q);1,2;B^4},\;1\le
i\le 8,
\end{equation}
and therefore,
\begin{equation}
 \label{111}
\|(a^{\perp})^{\ttilde}\|_{\tilde{A}(\q);1,2;\R^4}\lesssim\epsilon\|a^{\perp}\|_{A(\q);1,2;B^4}.
\end{equation}
We thus obtain
\begin{align}
\label{112}
&\mathcal{H}_{\tilde{A}(\q)}(a^{\perp},a^{\perp})=\mathcal{H}_{\tilde{A}(\q)}((a^{\perp})^{\tperp},(a^{\perp})^{\tperp})
+2\mathcal{H}_{\tilde{A}(\q)}((a^{\perp})^{\tperp},(a^{\perp})^{\ttilde})+
\mathcal{H}_{\tilde{A}(\q)}((a^{\perp})^{\ttilde},(a^{\perp})^{\ttilde})\notag\\
&\ge\mathcal{H}_{\tilde{A}(\q)}((a^{\perp})^{\tperp},(a^{\perp})^{\tperp})-
C\|(a^{\perp})^{\tperp}\|_{\tilde{A}(\q);1,2;\R^4}\|(a^{\perp})^{\ttilde}\|_{\tilde{A}(\q);1,2;\R^4}-
C\|(a^{\perp})^{\ttilde}\|_{\tilde{A}(\q);1,2;\R^4}^2\notag\\
&\ge C\|(a^{\perp})^{\tperp}\|^2_{\tilde{A}(\q);1,2;\R^4}-
C\epsilon\|(a^{\perp})^{\tperp}\|_{\tilde{A}(\q);1,2;\R^4}\|a^{\perp}\|_{A(\q);1,2;B^4}-
C\epsilon^2\|a^{\perp}\|_{A(\q);1,2;B^4}^2,
\end{align}
where we used Lemma \ref{L3.4} to estimate
$\mathcal{H}_{\tilde{A}(\q)}((a^{\perp})^{\tperp},(a^{\perp})^{\tperp})$.

From $(a^{\perp})^{\tperp}=a^{\perp}-(a^{\perp})^{\ttilde}$,
estimate \eqref{111} and
$\|a^{\perp}\|_{\tilde{A}(\q);1,2;\R^4}\approx\|a^{\perp}\|_{A(\q);1,2;B^4}$
(since $\text{supp}\,a^{\perp}\subset B^4$), one obtains
\begin{equation}
\label{113}
 C^{-1}\|a^{\perp}\|_{A(\q);1,2;B^4}\le
\|(a^{\perp})^{\tperp}\|_{\tilde{A}(\q);1,2;\R^4}\le
C\|a^{\perp}\|_{A(\q);1,2;B^4}.
\end{equation}
Finally, combining \eqref{107}, \eqref{112}, \eqref{113},
\begin{equation}
\label{114} \mathcal{H}_{A(\q)}(a^{\perp},a^{\perp})\ge
C\|a^{\perp}\|_{A(\q);1,2;B^4}^2
\end{equation}
for all small $\epsilon>0$, which is `almost' the assertion of Lemma
\ref{L3.3}. It is left to prove that the constant $C$ can be taken
independent of $\q\in\PP(d_0,\lambda_0;D_1,D_2;\epsilon)$. To this
purpose, we first observe that $\mathcal{H}_{\tilde{A}}(a,a)$ and
$\|{\nabla_{\tilde{A}}}^{\epsilon}a\|_{2;\R^4}$ are conformally
invariant. Hence, the inequality
$\mathcal{H}_{\tilde{A}(\q)}((a^{\perp})^{\tperp},(a^{\perp})^{\tperp})\ge
C\|{\nabla_{\tilde{A}(\q)}}^{\epsilon}(a^{\perp})^{\tperp}\|_{2;\R^4}^2$
holds for some $C>0$ independent of $\q$. By the Poincar\'e
inequality, we have $\|a\|_2\le C\|\nabla |a|\|_2\le
C\|{\nabla_{\tilde{A}(\q)}}^{\epsilon}a\|_2$ for $a\in
L^2_{1,0}(T^{\ast}B^4\otimes\Ad(P))$ and some $C>0$ independent of
$\q$. Thus
\begin{align*}
\|{\nabla_{\tilde{A}(\q)}}^{\epsilon}(a^{\perp})^{\tperp}\|_{2;\R^4}&\ge\|{\nabla_{\tilde{A}(\q)}}^{\epsilon}a^{\perp}\|_{2;\R^4}-
\|(a^{\perp})^{\ttilde}\|_{\tilde{A}(\q);1,2;\R^4}\\
&\ge C\|a^{\perp}\|_{A(\q);1,2;B^4}-C\epsilon\|a^{\perp}\|_{A(\q);1,2;B^4}\\
&\ge C\|a^{\perp}\|_{A(\q);1,2;B^4}
\end{align*}
for some $C>0$ independent of $\q$. The assertion follows.

\noindent Note that  \eqref{114}  yields the estimate
\begin{equation}
\label{115} \mathcal{H}_{A(\q)}(a,a)\ge
C\|a^{\perp}\|_{A(\q);1,2;B^4}^2-C\|a^{\intercal}\|_{A(\q);1,2;B^4}^2,
\end{equation}
which is used in the next section. \hfill$\Box$

\medskip

\subsection{The auxiliary equation}
In this section we solve the equation in
$T_{A(\q)}\mathcal{N}(d_0,\lambda_0)^{\perp}$ (i.e., essentially
orthogonally to the kernel of the Hessian of the $\epsilon$-Yang
Mills functional, which is the obstruction to the direct application
of the implicit function theorem). Thus we introduce the following
\textit{auxiliary equation}, associated to the Yang Mills equation
$\nabla\YMe(A(\q)+a)=0$:
\begin{equation}
\label{aux}
Q\Big(\frac{1}{2}\nabla\YMe(A(\q)+a)+{d_{A(\q)+a}}^{\epsilon}{d_{A(\q)+a}^{\ast}}^{\epsilon}a\Big)=0\,,
\end{equation}
where $Q$ is the topological projection defined in \eqref{tpr}. We
shall solve \eqref{aux} for $a\in
L^p_{1,0}(T^{\ast}B^4\otimes\Ad(P))\cap
T_{A(\q)}\mathcal{N}(d_0,\lambda_0)^{\perp}$.

 For $2<p<4$, we define the following duality paring:
 \begin{align}
 &L^p_{1,0}(T^{\ast}B^4\otimes\Ad(P))\times L^{p'}_{1,0}(T^{\ast}B^4\otimes\Ad(P))\to \R\notag\\
 &(a,b)\to \langle a,b\rangle :=\int_{B^4}({\nabla_{A(\q)}}^{\epsilon}a,{\nabla_{A(\q)}}^{\epsilon}b)\,dx+\int_{B^4}(a,b)\,dx\,,\notag
 \end{align}
 for $a\in L^p_{1,0}(T^{\ast}B^4\otimes\Ad(P))$, $b\in L^{p'}_{1,0}(T^{\ast}B^4\otimes\Ad(P))$, with $1/p+1/p'=1$.

 \noindent For $a,b\in C^{\infty}_0(T^{\ast}B^4\otimes\Ad(P))$, we define
  \begin{align}
&\Big\langle\frac{1}{2}\nabla\YMe(A(\q)+a)+{d_{A(\q)+a}}^{\epsilon}{d_{A(\q)+a}^{\ast}}^{\epsilon}a,b\Big\rangle\notag\\
&=\int_{B^4}({F_{A(\q)+a}}^{\epsilon},{d_{A(\q)+a}}^{\epsilon}b)\,dx+\int_{B^4}({d_{A(\q)+a}^{\ast}}^{\epsilon}a,{d_{A(\q)+a}^{\ast}}^{\epsilon}b)\,dx\notag\\
&=\int_{B^4}({F_{A(\q)}}^{\epsilon}+{d_{A(\q)}}^{\epsilon}a+\frac{\epsilon}{2}[a,a],{d_{A(\q)}}^{\epsilon}b+\epsilon[a,b])\,dx
+\int_{B^4}({d_{A(\q)}^{\ast}}^{\epsilon}a+\epsilon\ast[a,\ast
a],{d_{A(\q)}^{\ast}}^{\epsilon}b+\epsilon\ast[a,\ast b])\,dx.\notag
 \end{align}
By the Sobolev embeddings $L^p_1(B^4)\subset L^{4p/(4-p)}(B^4)$,
$L^{p'}_1(B^4)\subset L^{4p/(3p-4)}(B^4)$, one obtains
 \begin{align}
 \label{smoothH}
&\Big\langle\frac{1}{2}\nabla\YMe(A(\q)+a)+{d_{A(\q)+a}}^{\epsilon}{d_{A(\q)+a}^{\ast}}^{\epsilon}a,b\Big\rangle\notag\\
\le &\,
(\|{F_{A(\q)}}^{\epsilon}\|_{\infty}+C\|a\|_{A(\q);1,p}+C\|a\|_{A(\q);1,p}^2+C\|a\|_{A(\q);1,p}^3)\|b\|_{A(\q);1,p'}\,,
 \end{align}
 for some constant $C>0$ independent of $a$ and $b$. It follows
 that
 $$ \sup\Big\{\Big\langle\frac{1}{2}\nabla\YMe(A(\q)+a)+{d_{A(\q)+a}}^{\epsilon}{d_{A(\q)+a}^{\ast}}^{\epsilon}a,b\Big\rangle:\|b\|_{A(\q);1,p'}\le1\Big\}<\infty\,,$$
 for $a\in L^p_{1,0}(T^{\ast}B^4\otimes\Ad(P))$, thus
 $\frac{1}{2}\nabla\YMe(A(\q)+a)+{d_{A(\q)+a}}^{\epsilon}{d_{A(\q)+a}^{\ast}}^{\epsilon}a\in L^p_{1,0}(T^{\ast}B^4\otimes\Ad(P)).$

\noindent We obtain the following existence lemma. The proof is
essentially an application of the contraction mapping principle
together with uniform estimates of the Hessian as given in Lemma
3.4.
\begin{lemma} \label{L3.8} There exist $2<p_0<4$ and $\epsilon_0>0$ such that for all $2<p<p_0$, $0<\epsilon<\epsilon_0$ and
$\q\in\PP(d_0,\lambda_0;D_1,D_2;\epsilon)$, there exists $\delta>0$
such that the auxiliary equation \eqref{aux} has a unique solution
$a=a(\q)\in L^p_{1,0}(T^{\ast}B^4\otimes\Ad(P))\cap
T_{A(\q)}\mathcal{N}(d_0,\lambda_0)^{\perp}$ which satisfies
$\|a(\q)\|_{A(\q);1,p}<\delta$ and
\begin{equation}
\label{120} \|a(\q)\|_{A(\q);1,2;B^4}\le
C\|\nabla\YMe(A(\q))\|_{A(\q);1,2,\ast}
\end{equation}
for some $C>0$ depending only on $d_0,\lambda_0,D_1$ and $D_2$.
\end{lemma}
\textit{Proof:} For $\q\in\PP(d_0,\lambda_0;D_1,D_2;\epsilon)$, we
define the functional $$F:L^p_{1,0}(T^{\ast}B^4\otimes\Ad(P))\cap
T_{A(\q)}\mathcal{N}(d_0,\lambda_0)^{\perp} \to
L^p_{1,0}(T^{\ast}B^4\otimes\Ad(P))\cap
T_{A(\q)}\mathcal{N}(d_0,\lambda_0)^{\perp}\,,$$ by
$$F(a):=Q\Big(\frac{1}{2}\nabla\YMe(A(\q)+a)+
{d_{A(\q)+a}}^{\epsilon}{d_{A(\q)+a}^{\ast}}^{\epsilon}a\Big).$$
We show that
\begin{align}\notag
F'(0)&=Q\Big(\frac{1}{2}\nabla^2\YMe(A(\q))+d_{A(\q)}d_{A(\q)}^{\ast}\Big):\notag\\
&L^p_{1,0}(T^{\ast}B^4\otimes\Ad(P))\cap
T_{A(\q)}\mathcal{N}(d_0,\lambda_0)^{\perp}\to
L^p_{1,0}(T^{\ast}\otimes\Ad(P))\cap
T_{A(\q)}\mathcal{N}(d_0,\lambda_0)^{\perp}\notag
\end{align} is an isomorphism. Moreover, for our purpose, we shall give an estimate of the inverse norm for each $\q\in\PP(d_0,\lambda_0;D_1,D_2;\epsilon)$ which depends only on $d_0$, $D_1$, $D_1$ and $\epsilon$. To prove the invertibility, suppose
that $F'(0)a=0$ for $a\in L^p_{1,0}(T^{\ast}B^4\otimes\Ad(P))\cap
T_{A(\q)}\mathcal{N}(d_0,\lambda_0)^{\perp}$. We then have
\begin{equation}
\notag 0=\langle F'(0)a,a\rangle
=\frac{1}{2}\langle\nabla^2\YMe(A(\q))a,a\rangle+\langle
{d_{A(\q)}^{\ast}}^{\epsilon}a,{d_{A(\q)}^{\ast}}^{\epsilon}a\rangle
\ge C\|a\|_{A(\q);1,2;B^4}^2,
\end{equation}
where the last inequality comes from Lemma \ref{L3.3}. Therefore,
$a=0$ and $F'(0)$ is one to one.

To show that $F'(0)$ is onto, one needs to solve the equation
$$F'(0)a=Q\big(\frac{1}{2}\nabla^2\YMe(A(\q))a+{d_{A(\q)}}^{\epsilon}{d_{A(\q)}^{\ast}}^{\epsilon}a\big)=b$$
for arbitrary $b\in L^p_{1,0}(T^{\ast}\otimes\Ad(P))\cap
T_{A(\q)}\mathcal{N}(d_0,\lambda_0)^{\perp}$. This is equivalent to
\begin{equation}
\label{122} \Big\langle
Q\Big(\frac{1}{2}\nabla^2\YMe(A(\q))a+{d_{A(\q)}}^{\epsilon}{d_{A(\q)}^{\ast}}^{\epsilon}a\Big),\varphi\Big\rangle=\langle
b,\varphi\rangle :=(b,\varphi)_{A(\q);1,2;B^4}
\end{equation}
for all $\varphi\in L^{p'}_{1,0}(T^{\ast}B^4\otimes\Ad(P))$.

By a density argument, it is sufficient to show the existence of
$a\in L^p_{1,0}(T^{\ast}B^4\otimes\Ad(P))\cap
T_{A(\q)}\mathcal{N}(d_0,\lambda_0)^{\perp}$ such that \eqref{122}
holds for $\varphi\in L^2_{1,0}(T^{\ast}\otimes\Ad(P))$. Since
\eqref{122} is always satisfied for $\varphi\in
T_{A(\q)}\mathcal{N}(d_0,\lambda_0)$, we may also assume that
$\varphi\in L^2_{1,0}(T^{\ast}B^4\otimes\Ad(P))\cap
T_{A(\q)}\mathcal{N}(d_0,\lambda_0)^{\perp}$. Solving \eqref{122} is
equivalent to finding a critical point for the functional $a\mapsto
\frac{1}{4}\langle\nabla^2\YMe(A(\q))a,a\rangle+\frac{1}{2}\langle
{d_{A(\q)}}^{\epsilon}{d_{A(\q)}^{\ast}}^{\epsilon}a,a\rangle-\langle
b,a\rangle$ in $L^p_{1,0}(T^{\ast}B^4\otimes\Ad(P))\cap
T_{A(\q)}\mathcal{N}(d_0,\lambda_0)^{\perp}$. By Lemma \ref{L3.3},
there is a unique critical point $a\in
L^2_{1,0}(T^{\ast}B^4\otimes\Ad(P))\cap
T_{A(\q)}\mathcal{N}(d_0,\lambda_0)^{\perp}$, which also satisfies
the estimate
\begin{equation}
\label{123} \|a\|_{A(\q);1,2;B^4}\le C\|b\|_{A(\q);1,2;B^4}\le
C\|b\|_{A(\q);1,p;B^4}.
\end{equation}
We show that this solution is actually in
$L^p_{1,0}(T^{\ast}B^4\otimes\Ad(P))$. By the Weitzenb\"ock formula,
\eqref{122} can be written as
\begin{equation}
\notag
({\nabla_{A(\q)}}^{\epsilon}a,{\nabla_{A(\q)}}^{\epsilon}\varphi)_{2;B^4}=-\epsilon({F_{A(\q)}}^{\epsilon},[a,\varphi])_{2;B^4}+
\epsilon(\{{F_{A(\q)}}^{\epsilon},\varphi\},a)_{2;B^4}+ \langle
b,\varphi\rangle\,.
\end{equation}
H\"older's inequality and Sobolev embeddings yield
\begin{equation}
\label{125}
|({\nabla_{A(\q)}}^{\epsilon}a,{\nabla_{A(\q)}}^{\epsilon}\varphi)_{2;B^4}|\le
C\|\epsilon{F_{A(\q)}}^{\epsilon}\|_{p;B^4}\|a\|_{A(\q);1,2;B^4}\|
\varphi\|_{A(\q);1,p';B^4}+\|b\|_{A(\q);1,p;B^4}\|\varphi\|_{A(\q);1,p';B^4},
\end{equation}
where $C>0$ depends only on $d_0,\lambda_0,D_1$ and $D_2$.

On the other hand, by H\"older's inequality and Sobolev embedding,
\begin{equation}
\label{126} |(a,\varphi)_{2;B^4}|\le
C\|a\|_{A(\q);1,2;B^4}\|\varphi\|_{A(\q);1,p';B^4}
\end{equation}
for $2<p\le 4$, where $C>0$ is an absolute constant.

From \eqref{125}, \eqref{126}, it follows that
\begin{equation}
\label{127}
\biggl|\int_{B^4}({\nabla_{A(\q)}}^{\epsilon}a,{\nabla_{A(\q)}}^{\epsilon}\varphi)+(a,\varphi)\,dx
\biggr|\le
C(1+\|\epsilon{F_{A(\q)}}^{\epsilon}\|_{p;B^4})(\|a\|_{A(\q);1,2;B^4}+\|b\|_{A(\q);1,p;B^4})\|\varphi\|_{A(\q);1,p';B^4}\,,
\end{equation}
thus (by \eqref{123}),
$({\nabla_{A(\q)}^{\ast}}^{\epsilon}{\nabla_{A(\q)}}^{\epsilon}+1)a\in
L^p_{-1}(T^{\ast}B^4\otimes\Ad(P))$ with
$\|({\nabla_{A(\q)}^{\ast}}^{\epsilon}{\nabla_{A(\q)}}^{\epsilon}+1)a\|_{A(\q);-1,p;B^4}\le
C(1+\|\epsilon{F_{A(\q)}}^{\epsilon}\|_{p;B^4})\|b\|_{A(\q);1,p;B^4}.$
Since we have the estimate
$\|\epsilon{F_{A(\q)}}^{\epsilon}\|_{p;B^4}\le
C\epsilon^{\frac{2}{p}-1}$ for some constant $C>0$ depending only on
$d_0$, $D_1$ and $D_2$, it finally follows that
$\|a\|_{A(\q);1,p;B^4}\le
C(1+\epsilon^{\frac{2}{p}-1})\|b\|_{A(\q);1,p;B^4}$, where $C>0$
depends only on $d_0,\lambda_0,D_1, D_2$. This completes the proof
of the invertibility of $F'(0)$ and the estimate of the norm of
$F'(0)^{-1}$ for $\q\in\PP(d_0,\lambda_0;D_1,D_2;\epsilon)$.

From this, the existence of $a=a(\q)\in
L^p_{1,0}(T^{\ast}B^4\otimes\Ad(P))\cap
T_{A(\q)}\mathcal{N}(d_0,\lambda_0)^{\perp}$ satisfying \eqref{aux}
follows directly from the contraction mapping theorem. In fact,
$F(a)=0$ can be written as
\begin{equation}
\label{128} QF'(0)a=-
Q\Big(\frac{1}{2}\nabla\YMe(A(\q))\Big)-Q\Big(\frac{1}{2}R(\q;a)+R_2(\q;a)\Big)\,,
\end{equation}
where
$R_2(\q;a)={d_{A(\q)+a}}^{\epsilon}{d_{A(\q)+a}^{\ast}}^{\epsilon}a-{d_{A(\q)}}^{\epsilon}{d_{A(\q)}^{\ast}}^{\epsilon}a$,
and, since $QF'(0): L^p_{1,0}(T^{\ast}B^4\otimes\Ad(P))\cap
T_{A(\q)}\mathcal{N}(d_0,\lambda_0)^{\perp}\to
L^p_{1,0}(T^{\ast}B^4\otimes\Ad(P))\cap
T_{A(\q)}\mathcal{N}(d_0,\lambda_0)^{\perp}$ admits an inverse,
\eqref{128} is equivalent to
\begin{equation}
\label{129}
 a=T(\q;a)
:=-(QF'(0))^{-1}Q\Big(\frac{1}{2}\nabla\YMe(A(\q))\Big)-(QF'(0))^{-1}Q\Big(\frac{1}{2}R(\q;a)+R_2(\q;a)\Big)\,,
\end{equation}
where $T(\q;\cdot)$ is a contraction from the ball of radius
$\delta$ in $L^p_{1,0}(T^{\ast}B^4\otimes\Ad(P))\cap
T_{A(\q)}\mathcal{N}(d_0,\lambda_0)^{\perp}$ into itself, if $p>2$
is close to $2$ and $\delta$ is sufficiently small. To show this,
for $b\in L^{p'}_{1,0}(T^{\ast}B^4\otimes\Ad(P))$ we write
\begin{align}
&\qquad\qquad\Big\langle\frac{1}{2}R(\q;a_1)+R_2(\q;a_1)-\frac{1}{2}R(\q;a_2)-R_2(\q;a_2),b\Big\rangle=\notag\\
&\epsilon\int_{B^4}\bigl(({d_{A(\q)}}^{\epsilon}a_1,[a_1,b])-({d_{A(\q)}}^{\epsilon}a_2,[a_2,b])\bigr)dx+
\frac{\epsilon}{2}\int_{B^4}\bigl(([a_1,a_1],{d_{A(\q)}}^{\epsilon}b)-([a_2,a_2],{d_{A(\q)}}^{\epsilon}b)\bigr)dx\,+\notag\\
&\frac{\epsilon^2}{2}\int_{B^4}\bigl(([a_1,a_1],[a_1,b])-([a_2,a_2],[a_2,b])\bigr)dx+
\epsilon\int_{B^4}\bigl(({d_{A(\q)}^{\ast}}^{\epsilon}a_1,\ast[a_1,\ast b])-({d_{A(\q)}^{\ast}}^{\epsilon}a_2,\ast[a_2,\ast b])\bigr)dx\,+\notag\\
&\epsilon\int_{B^4}\bigl((\ast[a_1,\ast
a_1],{d_{A(\q)}^{\ast}}^{\epsilon}b)-(\ast[a_2,\ast
a_2],{d_{A(\q)}^{\ast}}^{\epsilon}b)\bigr)dx
+\epsilon^2\int_{B^4}\bigl(([a_1,\ast a_1],[a_1,\ast b])-([a_2,\ast
a_2],[a_2,\ast b])\bigr)dx\,.\notag
\end{align}
By the Sobolev embeddings $L^p_1\subset L^{4p/(4-p)}$,
$L^{p'}_1\subset L^{4p/(3p-4)}$, and H\"older's inequality, one
obtains
\begin{align}
&\Big|\Big\langle\frac{1}{2}R(\q;a_1)+R_2(\q;a_1)-\frac{1}{2}R(\q;a_2)-R_2(\q;a_2),b\Big\rangle\Big|\notag\\
&\le
C\epsilon(\|a_1\|_{1,p}+\|a_2\|_{1,p}+\epsilon\|a_1\|_{1,p}^2+\epsilon\|a_2\|_{1,p}^2)\|a_1-a_2\|_{1,p}\|b\|_{1,p'}
\notag\end{align} for some $C>0$ independent of $a_1$, $a_2$ and
$\epsilon$. Combining the estimate of the norm of $F'(0)^{-1}$ as
given above, the operator  $T(\q;\cdot)$ satisfies
$$\|T(\q;a_1)-T(\q;a_2)\|_{1,p}\le C(1+\epsilon^{\frac{2}{p}-1})\epsilon(\|a_1\|_{1,p}+\|a_2\|_{1,p}+\epsilon\|a_1\|_{1,p}^2+\epsilon\|a_2\|_{1,p}^2)\|a_1-a_2\|_{1,p}$$
and
\begin{equation}
\notag \|T(\q;a)\|_{1,p}\le
C(1+\epsilon^{\frac{2}{p}-1})(\|\nabla\YMe(A(\q))\|_{1,p}+\epsilon\|a\|_{1,p}^2+\epsilon^2\|a\|_{1,p}^3)\,,
\end{equation}
where we have the estimate of the form
$\|\nabla\YMe(A(\q))\|_{1,p}\le C\epsilon^{\alpha(p)}$ as
$\epsilon\to 0$, and $\alpha(p)\to1/2$ as $p\to 2$ (c.f. Lemma 3.2
and its proof). It follows that there exists $\delta>0$, $2<p_0<4$
and $\epsilon_0>0$ such that for all $0<\epsilon<\epsilon_0$ and
$2<p<p_0$, $T(\q;\cdot)$ is a contraction of the ball of radius
$\delta$ in $L^p_{1,0}(T^{\ast}B^4\otimes\Ad(P))\cap
T_{A(\q)}\mathcal{N}(d_0,\lambda_0)^{\perp}$. By the contraction
mapping theorem, for $0<\epsilon<\epsilon_0$ there exists a unique
solution of \eqref{129} in the ball of radius $\delta$. Note that
$a(\q)$ depends smoothly on $\q$, by the implicit function theorem.

\noindent To prove the estimate \eqref{120}, we observe that the
auxiliary equation \eqref{aux} yields
\begin{equation}
\label{134}
\Big\langle\frac{1}{2}\nabla\YMe(A(\q)+a),a\Big\rangle+({d_{A(\q)+a}^{\ast}}^{\epsilon}a,{d_{A(\q)+a}^{\ast}}^{\epsilon}a)_{2;B^4}=0\,,
\end{equation}
where
\begin{align}
\label{135}
({d_{A(\q)+a}^{\ast}}^{\epsilon}a,{d_{A(\q)+a}^{\ast}}^{\epsilon}a)_{2;B^4}&=({d_{A(\q)}^{\ast}}^{\epsilon}a,{d_{A(\q)}^{\ast}}^{\epsilon}a)_{2;B^4}+
2\epsilon({d_{A(\q)}^{\ast}}^{\epsilon}a,\ast[a,\ast a])_{2;B^4}
+\epsilon^2([a,\ast a],[a,\ast a])_{2;B^4}\notag\\
&\ge({d_{A(\q)}^{\ast}}^{\epsilon}a,{d_{A(\q)}^{\ast}}^{\epsilon}a)_{2;B^4}-C\epsilon(\|a\|_{A(\q);1,2;B^4}^3+\epsilon\|a\|_{A(\q);1,2;B^4}^4)\,.
\end{align}
Combining \eqref{134}, \eqref{135}, Lemma \ref{L3.2} and Lemma
\ref{L3.3}, we finally obtain
\begin{align}
C\|a\|_{A(\q);1,2;B^4}^2&\le
-\langle\nabla\YMe(A(\q)),a\rangle-\langle R(\q;a),a\rangle
+C\epsilon(\|a\|_{A(\q);1,2;B^4}^3+\epsilon\|a\|_{A(\q),1,2;B^4}^4)\notag\\
&\le \|\nabla\YMe(A(\q))\|_{A(\q);1,2,\ast}\|a\|_{A(\q);1,2;B^4}
+C\epsilon(\|a\|_{A(\q);1,2;B^4}^3+\epsilon\|a\|_{A(\q);1,2;B^4}^4),\notag
\end{align}
thus the required estimate \eqref{120}, since
$\|a\|_{A(\q);1,2;B^4}\le C\|a\|_{A(\q);1,p;B^4}$ is small.

This completes the proof. \hfill$\Box$

\subsection{Estimates for $\|a_{\q_i}(\q)\|_{A(\q);1,2;B^4}$, $1\leq i\leq 8$}

Let $a=a(\q)$ be as in Lemma \ref{L3.8}. We now estimate the
directional derivatives of $a(\q)$ in the direction $\q_i$, denoted
by $a_{\q_i}(\q)$. These estimates are needed to prove Proposition
\ref{P3.2}, which allows us to regard the problem of finding
multiple solutions to $(\mathcal D_\epsilon)$ as a finite
dimensional problem, and are also needed in \cite{IM2} where the
latter is solved.
\begin{lemma}
\label{L3.9} The following estimates hold:
\begin{align}
\|a_{\q_i}(\q)\|_{A(\q);1,2;B^4}&\lesssim\epsilon^{3/2}\quad \mbox{
for } 1\le
i\le 4\,,\\
\|a_{\q_i}(\q)\|_{A(\q);1,2;B^4}&\lesssim\epsilon\quad \mbox { for }
5\le i\le 8\,.
\end{align}
\end{lemma}
\textit {Proof:} To prove this lemma,  we write
$a_{\q_i}(\q)=a_{\q_i}(\q)^{\intercal}+a_{\q_i}(\q)^{\perp}$, where
$a_{\q_i}(\q)^{\intercal}=Pa_{\q_i}(\q)$ and
$a_{\q_i}(\q)^{\perp}=Qa_{\q_i}(\q)$ are defined in the course of
the proof of Lemma \ref{L3.3} (cf. \eqref{comp}), and estimate these
components separately.

\medskip\noindent \textit{Estimate of $\|a_{\q_j}(\q)^{\intercal}\|_{A(\q);1,2;B^4}$:}
For this estimate we need the following lemma, proved in \cite{IM3}.
\begin{lemma} \label{L3.10} Let $\A_i(\q)$ be the element of the orthonormal basis constructed in Lemma 3.2 in \cite{IM3}.
The following estimates hold:
\begin{equation}
\label{140}
\|{\A_i}_{\q_j}(\q)^{\perp}\|_{A(\q);1,2;B^4}\lesssim\epsilon
\end{equation}
for $1\le i,j\le 8$, where ${\A_i}_{\q_j}(\q)$ denotes the
directional derivative of $\A_i(\q)$ in the direction $\q_j$.
\end{lemma}
Since $a(\q)\in L^2_{1,0}(T^{\ast}B^4\otimes\Ad(P))\cap
T_{A(\q)}\mathcal{N}(d_0,\lambda_0)^{\perp}$, one has
$(a(\q),\A_i(\q))_{A(\q);1,2;B^4}=0$ for $1\le i\le 8$.
Differentiating this with respect to $\q_j$, one obtains
\begin{align}
\label{141}
&(a_{\q_j}(\q),\A_i(\q))_{A(\q);1,2;B^4}+(a(\q),{\A_i}_{\q_j}(\q))_{A(\q);1,2;B^4}\notag\\
&+\epsilon([\A_j(\q),a(\q)],{\nabla_{A(\q)}}^{\epsilon}\A_i(\q))_{2;B^4}+
\epsilon({\nabla_{A(\q)}}^{\epsilon}a(\q),[\A_j(\q),\A_i(\q)])_{2;B^4}=0\,.
\end{align}
By Lemma \ref{L3.1}, \eqref{120}, and Lemma \ref{L3.10}, the second
term of \eqref{141} is estimated as
\begin{equation}
|(a(\q),{\A_i}_{\q_j}(\q))_{A(\q);1,2;B^4}|\le\|a(\q)\|_{A(\q);1,2;B^4}\|{\A_i}_{\q_j}(\q)^{\perp}\|_{A(\q);1,2;B^4}
\lesssim\epsilon^{3/2},
\end{equation}
while the third and fourth terms are estimated as
\begin{align}
\epsilon|([\A_j(\q),a(\q)],&{\nabla_{A(\q)}}^{\epsilon}\A_i(\q))_{2;B^4}|\le C\epsilon\|\A_j(\q)\|_{4;B^4}\|a(\q)\|_{4;B^4}\|{\nabla_{A(\q)}}^{\epsilon}\A_i(\q)\|_{2;B^4}\notag\\
&\le
C\epsilon\|\A_j(\q)\|_{A(\q);1,2;B^4}\|a(\q)\|_{A(\q);1,2;B^4}\|\A_i(\q)\|_{A(\q);1,2;B^4}
\lesssim\epsilon^{3/2}
\end{align}
and
\begin{equation}
\label{144}
\epsilon|({\nabla_{A(\q)}}^{\epsilon}a(\q),[\A_j(\q),\A_i(\q)])_{2;B^4}|\lesssim\epsilon^{3/2}.
\end{equation}
From \eqref{141}--\eqref{144}, it follows
\begin{equation}
|(a_{\q_j}(\q),\A_i(\q))_{A(\q);1,2;B^4}|\lesssim\epsilon^{3/2}\quad\text{for
$1\le i\le 8$}
\end{equation}
 which yields finally
\begin{equation}
\label{146}
\|a_{\q_j}(\q)^{\intercal}\|_{A(\q);1,2;B^4}\lesssim\epsilon^{3/2}\quad\text{for
$1\le j\le 8$}.
\end{equation}

\medskip
\noindent\textit{Estimate of
$\|a_{\q_i}(\q)^{\perp}\|_{A(\q);1,2;B^4}$:} To estimate
$\|a_{\q_j}(\q)^{\perp}\|_{A(\q);1,2;B^4}$, recall that $a(\q)$
satisfies the auxiliary equation \eqref{aux}, thus there exist
$c_i(\q)\in\R$ ($1\le i\le 8$) such that the following holds:
\begin{equation}
\label{147}
\frac{1}{2}\nabla\YMe(A(\q)+a(\q))+{d_{A(\q)+a(\q)}}^{\epsilon}{d_{A(\q)+a(\q)}^{\ast}}^{\epsilon}a(\q)=\sum_{i=1}^8c_i(\q)\A_i(\q).
\end{equation}
From now on, we simply write $A=A(\q)$, $a=a(\q)$, $\A_i=\A_i(\q)$
and $c_i=c_i(\q)$.

We now estimate $|c_i|$ for $1\le i\le 8$.
\begin{lemma}
\label{L3.11} We have $|c_i(\q)|\lesssim\epsilon^{1/2}$ for $1\le
i\le 8$.
\end{lemma}
\textit{Proof of Lemma \ref{L3.11}:} By \eqref{147}, \eqref{120} and
Lemmas \ref{L3.1}, \ref{L3.2},
\begin{align}
|c_i|&=\Big|\Big\langle\frac{1}{2}\nabla\YMe(A+a),\A_i\Big\rangle+({d_{A+a}^{\ast}}^{\epsilon}a,{d_{A+a}^{\ast}}^{\epsilon}\A_i)_{2;B^4}\Big|\notag\\
&\le C(\|\nabla\YMe(A)\|_{A;1,2,\ast}+\|\nabla^2\YMe(A)\|_{A;1,2,\ast}\|a\|_{A;1,2;B^4}+\epsilon\|a\|_{A;1,2,B^4}^2\notag\\
&\quad+\epsilon^2\|a\|_{A;1,2;B^4}^3+\|a\|_{A;1,2;B^4})\lesssim~\epsilon^{1/2}.
\end{align}
This completes the proof of Lemma \ref{L3.11}. \hfill$\Box$

\smallskip
By differentiating \eqref{147} in the direction of $\q_i$ and taking
the pairing with $a_{\q_i}^{\perp}$, one obtains
\begin{align}
\label{149}
&\Big\langle\frac{1}{2}\nabla^2\YMe(A+a)(\A_i+a_{\q_i}),a_{\q_i}^{\perp}\Big\rangle+\langle
{d_{A+a}}^{\epsilon}{d_{A+a}^{\ast}}^{\epsilon}a_{\q_i},a_{\q_i}^{\perp}\rangle
+\epsilon\langle[\A_i+a_{\q_i},{d_{A+a}^{\ast}}^{\epsilon}a],a_{\q_i}^{\perp}\rangle
\notag\\
& +\epsilon\langle {d_{A+a}}^{\epsilon}{\ast}[\A_i+a_{\q_i},\ast
a],a_{\q_i}^{\perp}\rangle
=\sum_{j=1}^8c_j({\A_j}_{\q_i},a_{\q_i}^{\perp})_{A;1,2;B^4}\,,
\end{align}
thus, by expanding the left hand side of \eqref{149},
\begin{align}
\label{150}
&\mathcal{H}_A(a_{\q_i},a_{\q_i}^{\perp})+\frac{1}{2}\langle\nabla^2\YMe(A)\A_i,a_{\q_i}^{\perp}\rangle+
\epsilon({d_A^{\ast}}^{\epsilon}a_{\q_i},\ast[a,\ast a_{\q_i}^{\perp}])_{2;B^4}\notag\\
&+\epsilon(\ast[a,\ast
a_{\q_i}],{d_A^{\ast}}^{\epsilon}a_{\q_i}^{\perp})_{2;B^4}+\epsilon^2([a,\ast
a_{\q_i}],[a,\ast a_{\q_i}])_{2;B^4}
 +\frac{1}{2}\langle(\nabla^2\YMe(A+a)-\nabla^2\YMe(A))(\A_i+a_{\q_i}),a_{\q_i}^{\perp}\rangle\notag\\
 & +\epsilon([\A_i+a_{\q_i},{d_{A+a}^{\ast}}^{\epsilon}a],a_{\q_i}^{\perp})_{2;B^4}+\epsilon(\ast[\A_i+a_{\q_i},\ast a],{d_{A+a}^{\ast}}^{\epsilon}a_{\q_i}^{\perp})_{2;B^4}
=\sum_{j=1}^8c_j({\A_j}_{\q_i},a_{\q_i}^{\perp})_{A;1,2;B^4}\,.
\end{align}
By \eqref{115}, there holds
\begin{equation}
\mathcal{H}_A(a_{\q_i},a_{\q_i}^{\perp})=\mathcal{H}_A(a_{\q_i},a_{\q_i})-
\mathcal{H}_A(a_{\q_i},a_{\q_i}^{\intercal}) \ge
C\|a_{\q_i}^{\perp}\|_{A;1,2;B^4}^2-C\|a_{\q_i}^{\intercal}\|_{A;1,2;B^4}^2-\mathcal{H}_A(a_{\q_i},a_{\q_i}^{\intercal}).
\end{equation}
Combining the above with \eqref{150}, one obtains
\begin{align}
\label{152} &C\|a_{\q_i}^{\perp}\|_{A;1,2;B^4}^2\le
\mathcal{H}_A(a_{\q_i},a_{\q_i}^{\intercal})+C\|a_{\q_i}^{\intercal}\|_{A;1,2;B^4}^2+
\frac{1}{2}|\langle\nabla^2\YMe(A)\A_i,a_{\q_i}^{\perp}\rangle|\notag\\
&\quad +\epsilon|({d_A^{\ast}}^{\epsilon}a_{\q_i},\ast[a,\ast a_{\q_i}^{\perp}])_{2;B^4}|+\epsilon|(\ast[a,\ast a_{\q_i}],{d_A^{\ast}}^{\epsilon}a_{\q_i}^{\perp})_{2;B^4}|\notag\\
&\quad+\epsilon^2|([a,\ast a_{\q_i}],[a,\ast
a_{\q_i}^{\perp}])_{2;B^4}+\frac{1}{2}|\langle(\nabla^2\YMe(A+a)-
\nabla^2\YMe(A))(\A_i+a_{\q_i}),a_{\q_i}^{\perp}\rangle|\notag\\
&\quad+\epsilon|([\A_i+a_{\q_i},{d_{A+a}^{\ast}}^{\epsilon}a],a_{\q_i}^{\perp})_{2;B^4}|+
\epsilon|(\ast[\A_i+a_{\q_i},\ast
a],{d_{A+a}^{\ast}}^{\epsilon}a_{\q_i}^{\perp})_{2;B^4}|
 +\sum_{j=1}^8|c_j||({\A_j}_{\q_i},a_{\q_i}^{\perp})_{A;1,2;B^4}|.
\end{align}

We now estimate each term in \eqref{152}.

\noindent Applying \eqref{146}, we have
\begin{equation}
\label{153}
|\mathcal{H}_A(a_{\q_i},a_{\q_i}^{\intercal})|\lesssim\|a_{\q_i}\|_{A;1,2;B^4}\|a_{\q_i}^{\intercal}\|_{A;1,2;B^4}
\lesssim\epsilon^{3/2}\|a_{\q_i}^{\perp}\|_{A;1,2;B^4}+\epsilon^3.
\end{equation}
The third term in \eqref{152} is estimated as (we write
$\tilde{A}=\tilde{A}(\q)$)
\begin{align}
\label{154}  &|\langle\nabla^2\YMe(A)\A_i,a_{\q_i}^{\perp}\rangle|
\le|\langle(\nabla^2\YMe(A)-\nabla^2\YMe(\tilde{A}))\A_i,a_{\q_i}^{\perp}\rangle|+
|\langle\nabla^2\YMe(\tilde{A})\A_i,a_{\q_i}^{\perp}\rangle|\notag\\
&\le\|(\nabla^2\YMe(A)-\nabla^2\YMe(\tilde{A}))\A_i\|_{A;1,2,\ast}\|a_{\q_i}^{\perp}\|_{A;1,2;B^4}+
|\langle\nabla^2\YMe(\tilde{A})(\A_i-\tilde{\A}_i),a_{\q_i}^{\perp}\rangle|\notag\\
&\lesssim\epsilon^{3/2}\|a_{\q_i}^{\perp}\|_{A;1,2;B^4}+\|\A_i-
\tilde{\A}_i\|_{A;1,2;B^4}\|a_{\q_i}^{\perp}\|_{A;1,2;B^4},
\end{align}
where we have used $\nabla^2\YMe(\tilde{A})\tilde{\A}_i=0$ and Lemma
\ref{L3.7}.

By \eqref{120}, \eqref{146} and Lemma \ref{L3.1}, the fourth term is
estimated as
\begin{align}
\label{155}  &\epsilon|( {d_A^{\ast}}^{\epsilon}a_{\q_i},\ast[a,\ast
a_{\q_i}^{\perp}])_{2;B^4}|
\lesssim\epsilon\|a_{\q_i}\|_{A;1,2;B^4}\|a\|_{A;1,2;B^4}\|a_{\q_i}^{\perp}\|_{A;1,2;B^4}
\notag\\
 &\lesssim\epsilon^{3/2}\|a_{\q_i}\|_{A;1,2;B^4}\|a_{\q_i}^{\perp}\|_{A;1,2;B^4}
\lesssim\epsilon^3\|a_{\q_i}^{\perp}\|_{A;1,2;B^4}+\epsilon^{3/2}\|a_{\q_i}^{\perp}\|_{A;1,2;B^4}^2.
\end{align}
The fifth and the sixth terms are estimated similarly:
\begin{equation}
\label{156} \epsilon|(\ast[a,\ast
a_{\q_i}],{d_A^{\ast}}^{\epsilon}a_{\q_i}^{\perp})_{2;B^4}|\lesssim
\epsilon^3\|a_{\q_i}^{\perp}\|_{A;1,2;B^4}+\epsilon^{3/2}\|a_{\q_i}^{\perp}\|_{A;1,2;B^4}^2\,
\end{equation}
and
\begin{equation}
\label{157} \epsilon^2|([a,\ast a_{\q_i}],[a,\ast
a_{\q_i}])_{2;B^4}|\lesssim\epsilon^{9/2}\|a_{\q_i}^{\perp}\|_{A;1,2;B^4}+\epsilon^3\|a_{\q_i}^{\perp}\|_{A;1,2;B^4}^2\,.
\end{equation}
To estimate the seventh term, we first observe that
\begin{align}
\label{158}
&\frac{1}{2}\langle(\nabla^2\YMe(A+a)-\nabla^2\YMe(A))\alpha,\beta\rangle
=\int_{B^4}({d_A}^{\epsilon}\alpha,[a,\beta])+\epsilon\int_{B^4}([a,\alpha],{d_A}^{\epsilon}\beta)\notag\\
&\quad+\epsilon^2\int_{B^4}([a,\alpha],[a,\beta])+\epsilon\int_{B^4}({d_A}^{\epsilon}a,[\alpha,\beta])
+\frac{\epsilon^2}{2}\int_{B^4}([a,a],[\alpha,\beta])\notag\\
&\lesssim\epsilon\|\alpha\|_{A;1,2;B^4}\|\beta\|_{A;1,2;B^4}\|a\|_{A;1,2;B^4}+\epsilon^2\|\alpha\|_{A;1,2;B^4}\|\beta\|_{A;1,2;B^4}\|a\|_{A;1,2;B^4}^2\notag\\
&\lesssim\epsilon^{3/2}\|\alpha\|_{A;1,2;B^4}\|\beta\|_{A;1,2;B^4}
\end{align}
and
\begin{equation}
\label{159}
\|\nabla^2\YMe(A+a)-\nabla^2\YMe(A)\|_{A;1,2,\ast}\lesssim\epsilon^{3/2}.
\end{equation}
Thus,
\begin{align}
\label{160}
 |\langle(\nabla^2\YMe(A+a)-\nabla^2\YMe(A))(\A_i+a_{\q_i}),a_{\q_i}^{\perp}\rangle|
&\lesssim\epsilon^{3/2}(1+\|a_{\q_i}\|_{A;1,2;B^4})\|a_{\q_i}^{\perp}\|_{A;1,2;B^4}\notag\\
&\lesssim\epsilon^{3/2}\|a_{\q_i}^{\perp}\|_{A;1,2;B^4}+\epsilon^{3/2}\|a_{\q_i}^{\perp}\|_{A;1,2;B^4}^2\,.
\end{align}
The eighth and the ninth terms are estimated similarly:
\begin{equation}
\label{161}
\epsilon|([\A_i+a_{\q_i},{d_{A+a}^{\ast}}^{\epsilon}a],a_{\q_i}^{\perp})_{2;B^4}|
\lesssim\epsilon^{3/2}\|a_{\q_i}^{\perp}\|_{A;1,2;B^4}+\epsilon^{3/2}\|a_{\q_i}^{\perp}\|_{A;1,2;B^4}^2
\end{equation}
and
\begin{equation}
\label{162} \epsilon|(\ast[\A_i+a_{\q_i},\ast
a],{d_{A+a}^{\ast}}^{\epsilon}a_{\q_i}^{\perp})_{2;B^4}|
\lesssim\epsilon^{3/2}\|a_{\q_i}^{\perp}\|_{A;1,2;B^4}+\epsilon^{3/2}\|a_{\q_i}^{\perp}\|_{A;1,2;B^4}^2\,.
\end{equation}
Finally,  by Lemma 3.12 and \eqref{140}, the last term is estimated
as
\begin{equation}
\label{163}
\sum_{j=1}^8|c_j||({\A_j}_{\q_i},a_{\q_i}^{\perp})_{A;1,2;B^4}|=
\sum_{j=1}^8|c_j||({\A_j}_{\q_i}^{\perp},a_{\q_i}^{\perp})_{A;1,2;B^4}|
\lesssim\epsilon^{3/2}\|a_{\q_i}^{\perp}\|_{A;1,2;B^4}.
\end{equation}
From \eqref{152}--\eqref{163}, we obtain
\begin{equation}
\label{164}
\|a_{\q_i}^{\perp}\|_{A;1,2;B^4}^2\lesssim\epsilon^{3/2}\|a_{\q_i}^{\perp}\|_{A;1,2;B^4}+
\epsilon^3+\|\A_i-\tilde{\A}_i\|_{A;1,2;B^4}\|a_{\q_i}^{\perp}\|_{A;1,2;B^4}+
\epsilon^{3/2}\|a_{\q_i}^{\perp}\|_{A;1,2;B^4}^2\,.
\end{equation}
Thus, for small  $\epsilon>0$,
\begin{equation}
\label{165}
\|a_{\q_i}^{\perp}\|_{A;1,2;B^4}\lesssim\epsilon^{3/2}+\|\A_i-\tilde{\A}_i\|_{A;1,2;B^4}\,.
\end{equation}
From \eqref{165} and Lemma \ref{L3.6}, one obtains the estimates
\begin{equation}
\label{166}
\|a_{\q_i}^{\perp}\|_{A;1,2;B^4}\lesssim\epsilon^{3/2}\quad\text{for
$1\le i\le 4$}
\end{equation}
and
\begin{equation}
\label{167}
\|a_{\q_i}^{\perp}\|_{A;1,2;B^4}\lesssim\epsilon\quad\text{for $5\le
i\le 8$}.
\end{equation}
Finally, by combining \eqref{146}, \eqref{166}, \eqref{167}, we
complete the proof of Lemma \ref{L3.9}. \hfill$\Box$

\subsection{Natural constraints}

In this section, we prove that the manifold
$\{A(\q)+a(\q):\q\in\PP(d_0,\lambda_0;D_1,D_2;\epsilon)\}$ is a
natural constraint for $\YMe$ if $\epsilon>0$ is small, more
precisely, we prove the following proposition which allows us to
transform the $\epsilon$-Dirichlet problem for the Yang Mills
functional into a finite dimensional problem.

\begin{proposition}
\label{P3.2} There exists $\epsilon_0>0$ such that for
$0<\epsilon<\epsilon_0$, $\q\in\PP(d_0,\lambda_0;D_1,D_2;\epsilon)$
is a critical point for the function
$$\Je(\q):=\epsilon^2\YMe(A(\q)+a(\q))$$
on $\PP(d_0,\lambda_0;D_1,D_2;\epsilon)$ if and only if
$A(\q)+a(\q)$ is a Yang Mills connection, where $a(\q)$ is given by
Lemma \ref{L3.8}.
\end{proposition}
\textit{Proof:} It is obvious that if $A(\q)+a(\q)$ is Yang Mills,
then $\q$ is a critical point for $\Je$. Assume on the other hand
that $\q\in\PP(d_0,\lambda_0;D_1,D_2;\epsilon)$ is a critical point
for $\Je$. Then,
\begin{equation}
\label{168} 0=\langle\Je'(\q),\q_j\rangle
=\epsilon^2\nabla\YMe(A(\q)+a(\q))(\A_i+a_{\q_i}(\q))\quad\text{for
$1\le j\le 8$}.
\end{equation}
By \eqref{147}, we have
$\frac{1}{2}\nabla\YMe(A(\q)+a(\q))=\sum_{i=1}^8c_i(\q)\A_i$ on
$T_{A(\q)+a(\q)}\mathcal{B}^p_{k,+1}(A_0)$.  Thus, by \eqref{168}
and Lemma \ref{L3.9},
\begin{align}
\label{171}
 0&=\sum_{i=1}^8c_i(\q)(\A_i(\q),\A_j(\q)+a_{\q_j}(\q))_{A(\q)+a(\q);1,2;B^4}\notag\\
 &=\sum_{i=1}^8c_i(\q)\Big((\A_i(\q),\A_j(\q)+a_{\q_j}(\q))_{A(\q);1,2;B^4}+
({\nabla_{A(\q)}}^{\epsilon}\A_i(\q),\epsilon[a(\q),\A_j(\q)+a_{\q_j}(\q)])_{2;B^4}
\notag\\
&\quad+(\epsilon[a(\q),\A_i(q)],{\nabla_{A(\q)}}^{\epsilon}(\A_j(\q)+a_{\q_j}(\q)))_{2;B^4}+
(\epsilon[a(\q),\A_i(\q)],\epsilon[a(\q),\A_j(\q)])_{2;B^4}\Big)\notag\\
&=\sum_{i=1}^8c_i(\q)\delta_{ij}+o(|c(\q)|)
=c_j(\q)+o(|c(\q)|)\quad\text{for $1\le j\le 8$}.\notag\\
\end{align}
This implies $c_j(\q)=0$ for $1\le j\le 8$, if $\epsilon>0$ is
small, and $\nabla\YMe(A(\q)+a(\q))=0$ in
$T_{A(\q)+a(\q)}\mathcal{B}^p_{k,+1}(A_0)$, thus $A(\q)+a(\q)$ is
Yang Mills.

This completes the proof. \hfill$\Box$

\end{document}